\input amstex
\documentstyle{amsppt}
\magnification=1200
\hoffset=-0.5pc
\vsize=57.2truepc
\hsize=38truepc
\nologo
\spaceskip=.5em plus.25em minus.20em

\define\almemoli{1}
\define\bieriboo{2}
\define\bierieck{3}
\define\borhiron{4}
\define\boredmod{5}
\define\kbrownfo{6}
\define\cartanei{7}
\define\chemlone{8}
\define\derhaboo{9}
\define\evluwein{10}
\define\ginzwein{11}
\define\godebook{12}
\define\gugenhtw{13}
\define\hartshor{14}
\define\poiscoho{15}
\define\bv{16}
\define\extensta{17}
\define\koszulon{18}
\define\mackbook{19}
\define\maclaboo{20}
\define\rinehart{21}
\define\weinsfte{22}
\define\xuone{23}
\define\Bobb{\Bbb}
\define\fra{\frak}

\topmatter
\title 
Duality for Lie-Rinehart algebras and the modular class
\endtitle
\author Johannes Huebschmann
\endauthor
\affil 
Universit\'e des Sciences et Technologies de Lille
\\
UFR de Math\'ematiques
\\
59 655 VILLENEUVE D'ASCQ C\'edex
\\
Johannes.Huebschmann\@univ-lille1.fr
\endaffil
\date {October 18, 1998}
\enddate
\abstract{
A notion
of homological duality 
for Lie-Rinehart algebras
is 
studied which generalizes
the ordinary 
duality in Lie algebra (co)homology
and in the (co)homology of smooth manifolds.
The duality isomorphisms can be given
by cap products with suitable fundamental
classes and hence may be taken to be natural in any 
reasonable sense.
A precise notion of Poincar\'e duality,
meaning that certain bilinear pairings over appropriate
ground rings are nondegenerate,
is then introduced, and
various examples of Lie-Rinehart
algebras are shown to satisfy Poincar\'e duality.
Thereafter
a certain intrinsic
module 
introduced by Evens, Lu, and Weinstein
for Lie algebroids
is generalized to Lie-Rinehart algebras satisfying duality.
This module
determines a characteristic class,
called the modular class of the Lie-Rinehart algebra;
this class lies in
an appropriately defined Picard group 
generalizing the abelian group of
flat line bundles on a smooth manifold.
A Poisson algebra
having the requisite regularity properties
determines a corresponding module
for its Lie-Rinehart algebra
and hence modular class
whose square yields
the module and characteristic class
for its Lie-Rinehart algebra mentioned before.
These concepts arise from abstraction from the notion
of modular vector field 
for a smooth Poisson manifold. Finally,
it is shown that 
the Poisson cohomology 
of certain Poisson algebras
satisfies Poincar\'e duality.}
\endabstract
\endkeywords
\endsubjclass
\dedicatory
Beno Eckmann anl\"a{\ss}lich seines 80.Geburtstags in Dankbarkeit zugeeignet
\enddedicatory
\endtopmatter
\document
\leftheadtext{Johannes Huebschmann}
\rightheadtext{Duality for Lie-Rinehart algebras}
\medskip\noindent {\bf Introduction}
\smallskip\noindent
There are two ways of phrasing Poincar\'e duality
in the real cohomology of a compact orientable 
smooth  $n$-dimensional manifold
$W$, by means of the natural bilinear (cup) pairing
$$
\roman H^j(W) \otimes _{\Bobb R}  \roman H^{n-j}(W)
@>>>
\roman H^n(W),
\quad 0 \leq j \leq n,
\tag0.1
$$
and via
the natural map
$$
\roman H^j(W) @>>> \roman H_j(W),
\quad 0 \leq j \leq n,
\tag0.2
$$
in singular homology
which is obtained as the cap product with the fundamental class
of $W$.
Poincar\'e duality then amounts to the nondegeneracy of
(0.1) or, equivalently,
to (0.2) being an isomorphism.
Both notions of Poincar\'e duality  are equivalent
since the universal coefficient map
from
$\roman H^j(W)$ to
$\roman{Hom}(\roman H_j(W),\Bobb R)$
is an isomorphism.
Over more general rings, e.~g. the integers,
the appropriate way
of phrasing Poincar\'e duality
is in fact by means of the
isomorphism (0.2) but this will not be important here.
For our purposes, real cohomology of a smooth manifold
will simply mean de Rham cohomology.
\smallskip
In the (co)homology
of a finite dimensional Lie algebra 
over a field,
Poincar\'e duality may likewise be phrased
in two ways,
which in fact formally correspond
to the nondegeneracy of a pairing
similar to (0.1)
and to an isomorphism
of the kind (0.2).
\smallskip
In this paper we shall 
abstract from both de Rham cohomology and Lie algebra (co)homology
and study 
duality for
Lie-Rinehart algebras
(a definition will be given shortly).
There are examples which do neither
arise from the de Rham cohomology
of a smooth manifold nor from ordinary Lie algebras,
for example the Lie-Rinehart algebra
of a foliation or that corresponding to a Poisson algebra,
and these justify the general approach we offer here.
With the appropriate notion of homology,
the isomorphism (0.2) generalizes
to arbitrary Lie-Rinehart algebras.
We refer to this generalization simply as 
{\it duality\/}
(or sometimes
as {\it naive duality\/}).
The bilinear pairing (0.1) 
may be generalized as well
but 
the question whether the resulting
pairing is  nondegenerate
(over an appropriate ring which does not necessarily
coincide with the naive ground ring, see what is said below)
is then rather delicate,
and we do not know if
the answer is always positive.
We expect that,
for a general Lie-Rinehart algebra,
nondegeneracy 
reflects a certain regularity property, though,
as it does for topological spaces
where Poincar\'e duality singles out
the {\it Poincar\'e complexes\/}
(that is, cell complexes which satisfy the statement of the
Poincar\'e duality theorem with respect to some
fundamental class)
and in particular manifolds.
When nondegeneracy holds, we shall  say
that the Lie-Rinehart algebra satisfies
{\it Poincar\'e duality\/}.
Part of the paper is devoted to developing the necessary
language to make this precise.
\smallskip
We now give an outline of the paper.
Let $R$ be a commutative ring and $A$ a commutative $R$-algebra.
An $(R,A)$-{\it Lie algebra\/} \cite\rinehart\ 
 is a Lie algebra $L$ over $R$ 
which acts on (the left of) $A$ (by derivations) and
is also
an $A$-module
satisfying suitable compatibility conditions
which generalize the customary
properties of 
the Lie algebra
of smooth vector fields on a smooth manifold
viewed as a module over its ring of smooth functions;
these 
conditions read
$$
\aligned
[\alpha,a\beta] &=  \alpha(a) \beta + a [\alpha,\beta]
\\
(a \alpha) (b) &=  a (\alpha (b))
\endaligned
\tag0.3
$$
where  $a,b \in A$ and
$\alpha,\beta \in L$.
See our paper \cite\poiscoho\ 
for more details and history.
When the ground ring $R$ is not specified, we refer to 
the pair $(A,L)$ as a {\it Lie-Rinehart\/} algebra
to honor the fact that
Rinehart \cite\rinehart\  
proved early non-trivial results about these objects;
in particular, he established
a kind of Poincar\'e-Birkhoff-Witt theorem
(which we reproduce in the proof of Theorem 2.10 below).
Any
$(R,A)$-Lie algebra $L$
gives rise to a complex $\roman{Alt}_A(L,A)$
of alternating forms with the customary
Cartan-Chevalley-Eilenberg differential
\cite\rinehart\ 
which generalizes the de Rham complex
of a manifold {\it and\/} at 
the same time the  complex
computing 
Chevalley-Eilenberg
Lie algebra cohomology,
cf. 
e.~g. \cite\cartanei.
More generally,
there are notions of left- and 
right $(A,L)$-modules;
given left- and 
right $(A,L)$-modules $M$ and 
$N$, respectively,
the cohomology
$\roman H^*(L,M)$
and
homology $\roman H_*(L,N)$
are defined
in the usual way
by $\roman{Ext}_U^*(A,M)$ and
$\roman{Tor}^U_*(N,A)$, respectively,
where $U = U(A,L)$ is the universal algebra for $(A,L)$
(see Section 2 for details;
it is the algebra of differential operators when $A$ is the ring of
smooth functions on a smooth manifold and $L$ the Lie algebra of 
smooth vector fields).
When $L$ is projective as an $A$-module,
the cohomology
$\roman H^*(L,M)$
is computed
by the complex 
$\roman{Alt}_A(L,M)$. 
When $A$ and $L$ are the algebra of smooth functions
and Lie algebra of smooth vector fields, respectively,
on a smooth manifold,
the cohomology $\roman H^*(L,A)$
coincides with the de Rham cohomology of the manifold.
\smallskip
A Lie-Rinehart algebra
$(A,L)$ 
where, as an $A$-module,  $L$ 
is finitely generated projective 
of constant rank $n$ 
will be referred to 
as a {\it duality
Lie-Rinehart algebra of rank\/} $n$.
For a general
Lie-Rinehart algebra $(A,L)$,
the universal algebra
$U = U(A,L)$
inherits structures of left- and right
$(A,L)$-modules in the usual way.
Suppose now that
$(A,L)$ is
a  duality
Lie-Rinehart algebra of rank $n$.
 Consider
the highest non-zero cohomology group
$C_L = \roman H^{\roman{top}}(L,U)=\roman H^n(L,U)$,
with reference to the left
$(A,L)$-module structure on $U$.
Under the construction of $C_L$,
the
right $(A,L)$-module structure on $U$ remains free and induces
a right $(A,L)$-module structure on
$C_L$. In Section 2 we shall show the following:
\newline\noindent
(i) 
as right $(A,L)$-modules,
$C_L$ and $\roman{Hom}_A(\Lambda_A^n L,A)$
are isomorphic, the requisite right
$(A,L)$-module structure on 
$\roman{Hom}_A(\Lambda_A^n L,A)$ being given by
the negative of the Lie derivative;
 there are  
\newline\noindent
(ii) natural isomorphisms
$$
\roman H^k(L,M)
\cong
\roman H_{n-k}(L,C_L \otimes_A M)
\tag0.4
$$
for all non-negative integers $k$ and all
left $(A,L)$-modules $M$ and,  furthermore,
\newline\noindent
(iii) natural isomorphisms
$$
\roman H_k(L,N)
\cong
\roman H^{n-k}(L,\roman{Hom}_A(C_L,N))
\tag0.5
$$
for all non-negative integers $k$ and all
right $(A,L)$-modules $N$.
See (2.11) below for details;
the requisite
left- and right
$(A,L)$-module structures
are there explained as well.
A version of the isomorphism (0.4) may be found in
\cite{\chemlone\ (5.2.2)}.
\smallskip
In Section 3 we shall show that
the duality isomorphisms
(0.4) and (0.5)
can be given by a cap product 
with a certain fundamental class
in $\roman H_n(L,C_L)$;
hence they are or may be taken to be natural
in any reasonable sense.
These versions of
the isomorphisms (0.4) and (0.5)
are the direct generalization of the isomorphism (0.2)
spelled out above.
For a general duality Lie-Rinehart algebra,
they entail
the existence of certain bilinear cohomology pairings
which are not necessarily nondegenerate, though;
see (3.10) and (7.13.1) below.
Special cases of these pairings 
arising from certain Lie algebroids 
are given in
\cite\evluwein.
\smallskip
In Section 4 we shall introduce a notion of Poincar\'e duality
for a general duality Lie-Rinehart algebra,
the terminology \lq\lq Poincar\'e\rq\rq\ 
referring to the nondegeneracy of certain cohomology pairings
generalizing (0.1) above.
A crucial ingredient
is the notion of what we call a {\it trace\/}.
For example,
given 
a smooth $n$-dimensional manifold $W$,
a trace 
for the Lie-Rinehart algebra
$(A,L) = (C^{\infty}(W), \roman{Vect}(W))$
consists of the space
$\Cal O$ of (compactly supported) sections of the orientation bundle
of $W$ and an
isomorphism
$t$ from
$\roman H^n(L,\Cal O) \ (\cong \roman H^n(W,\Bobb R_t)$
where $\Bobb R_t$ refers to the local system defined by
the orientation bundle)
to $\Bobb R$
(which is in fact given by integration).
Our  notion of \lq trace\rq\  is similar to
a corresponding one
in the theory of Serre duality for the cohomology
of projective sheaves on a projective scheme,
cf. e.~g. \cite\hartshor.
We shall prove that various duality Lie-Rinehart algebras,
including
the Lie-Rinehart algebra
$(A,L) = (C^{\infty}(W), \roman{Vect}(W))$
where $W$ is a smooth manifold
and, more generally,
the Lie-Rinehart algebra
describing the infinitesimal
structure
of a fibre bundle with compact fibre,
satisfy Poincar\'e duality for appropriate
modules.
\smallskip
In Section 5 we
shall show that, in a suitable sense, Poincar\'e duality
is preserved under extensions of
Lie-Rinehart algebras; see (5.3) for details.
We then deduce that,
cf. (5.4),
every Lie-Rinehart algebra arising
from a transitive Lie algebroid on a smooth manifold
satisfies Poincar\'e duality.
\smallskip
In Section 6 we shall introduce a certain
Picard group 
$\roman{Pic}^{\roman{flat}}(L,A)$
of 
projective rank one $A$-modules with a 
left $(A,L)$-module structure.
This group could be thought of as a group of flat line bundles,
the group structure
being induced by the operation
of tensor product over $A$.
Thereafter we introduce a certain 
left $(A,L)$-module
$Q_L$,
for the special case
where the algebra $A$ is regular in the sense that,
as an $A$-module,
the $(R,A)$-Lie algebra
$\roman{Der}(A)$
is finitely generated and projective.
In this case, we denote the $A$-dual
$\roman {Hom}_A(\Lambda_A^n \roman{Der}(A), A)$
of the highest non-zero exterior power 
$\Lambda_A^n \roman{Der}(A)$
of 
$\roman{Der}(A)$
by $\omega_A$.
For example, the algebra of smooth functions
on a smooth manifold will be regular in this sense.
The $A$-module
$\omega_A$
may be identified with
the highest non-zero exterior power of the module 
of formal differentials of $A$ or of a suitable descendant thereof. 
In algebraic geometry, 
a closely related
object is called
the \lq\lq canonical sheaf\rq\rq\ 
and denoted by $\omega$ with an appropriate subscript
whence the notation.
When
$A$ is the algebra of smooth functions
on a smooth manifold,
$\omega_A$
is the space of sections
of the highest non-zero exterior power of the
cotangent bundle.
In the general case,
the negative of the operation of Lie derivative endows
$\omega_A$
with 
a right
$(A,L)$-module structure, and
the module $Q_L$ is then defined
by $Q_L = \roman{Hom}_A(C_L,\omega_A)$;
the customary diagonal action,
with reference to the right
$(A,L)$-module structures on
$C_L$  and $\omega_A$,
endows $Q_L$ with 
a {\it left\/}
$(A,L)$-module structure.
When $A$ is the algebra
of smooth real functions on a smooth
manifold and $L$ arises from a foliation,
the
left $(A,L)$-module structure on
$Q_L$
amounts to the Bott connection for this foliation.
For the special case of
a 
Lie-Rinehart algebra arising from a
Lie algebroid,
the
left $(A,L)$-module 
$Q_L$
has been introduced in \cite\evluwein.
For a general duality Lie-Rinehart algebra 
$(A,L)$, the 
class
$[Q_L] \in \roman{Pic}^{\roman{flat}}(L,A)$
of the left $(A,L)$-module 
$Q_L$ is then {\it characteristic\/} for 
$(A,L)$.
We refer to $[Q_L]$  as the {\it modular class\/}
of $(A,L)$.
Our construction 
of this modular class
is related with a similiar construction
in \cite\evluwein\ 
but substantially differs from
the one in that paper
and is more general.
See Section 6 below for details.
Our more formal approach
involving Lie-Rinehart algebras 
over an arbitrary commutative algebra
and the notion of duality
explained in Sections 1 and 2 below
and that of Poincar\'e duality
given in Section 4
place the results of \cite\evluwein\  in their proper context,
accordingly simplify the approach in that paper,
and
clarify some of the points related with 
existence of the modular class and
the question whether or not
certain cohomology pairings are nondegenerate,
which have been left open there.
\smallskip
Finally, in Section 7 we shall apply our ideas and
constructions to a Poisson algebra $A$
which, as a commutative algebra, is assumed to be regular.
Here the requisite $(R,A)$-Lie algebra
$D_{\{\cdot,\cdot\}}$
is that arising from the formal differentials
of $A$, introduced and studied in our paper
\cite\poiscoho,
and the theory developed earlier applies.
In particular,
we have
the corresponding
left $(A,D_{\{\cdot,\cdot\}})$-module
$Q_{D_{\{\cdot,\cdot\}}}$ at our disposal.
From its Poisson structure,
$A$ also inherits
a  {\it right\/} 
$(A,D_{\{\cdot,\cdot\}})$-module
structure;  we denote the resulting
right
$(A,D_{\{\cdot,\cdot\}})$-module
by 
$A_{\{\cdot,\cdot\}}$.
This observation greatly simplifies the relevant calculations
in \cite\evluwein\ 
and allows for further generalization
of the corresponding results in
\cite\evluwein\ and
\cite\weinsfte.
In fact, 
$C_{\{\cdot,\cdot\}}$
denoting 
the dualizing module
for
$D_{\{\cdot,\cdot\}}$  with its
right
$(A,D_{\{\cdot,\cdot\}})$-module structure, 
the diagonal action
(cf. (2.3) below)
endows the $A$-module 
$\roman{Hom}_A(C_{\{\cdot,\cdot\}},A_{\{\cdot,\cdot\}})$
with a  {\it left\/}
$(A,D_{\{\cdot,\cdot\}})$-module structure.
The 
left
$(A,D_{\{\cdot,\cdot\}})$-module structure
on
$\roman{Hom}_A(C_{\{\cdot,\cdot\}},A_{\{\cdot,\cdot\}})$
generalizes
the modular vector field for a smooth Poisson manifold 
studied in \cite\weinsfte\ and elsewhere
and, furthermore, offers some additional formal insight
into this concept;
some comments will be given after Lemma 7.8 below.
Moreover, 
Theorem 7.9
will say that
$Q_{D_{\{\cdot,\cdot\}}}$
and
$\roman{Hom}_A(C_{\{\cdot,\cdot\}},A_{\{\cdot,\cdot\}})\otimes_A
\roman{Hom}_A(C_{\{\cdot,\cdot\}},A_{\{\cdot,\cdot\}})$
are isomorphic as
left
$(A,D_{\{\cdot,\cdot\}})$-modules.
This generalizes a result in \cite\evluwein.
See (7.10) and (7.14) below.
Moreover,
the 
class
$[\roman{Hom}_A(C_{\{\cdot,\cdot\}},A_{\{\cdot,\cdot\}})] 
\in \roman{Pic}^{\roman{flat}}(D_{\{\cdot,\cdot\}},A)$
of the left $(A,D_{\{\cdot,\cdot\}})$-module 
$\roman{Hom}_A(C_{\{\cdot,\cdot\}},A_{\{\cdot,\cdot\}})$ 
is then {\it characteristic\/} for 
the Poisson algebra
$A$.
We refer 
to
it as the {\it modular class\/} of $A$.
In view of 
the relationship between
$Q_{D_{\{\cdot,\cdot\}}}$
and
$\roman{Hom}_A(C_{\{\cdot,\cdot\}},A_{\{\cdot,\cdot\}})$ 
which we just mentioned,
the modular class of the Lie-Rinehart algebra
$(A,D_{\{\cdot,\cdot\}})$
is twice the modular class 
(or its square when we think of the group structure as 
being multiplicative)
of the Poisson algebra $A$.
This generalizes the corresponding observation
spelled out in \cite\evluwein.
A related but different
modular class for a Poisson manifold
has been introduced in \cite\weinsfte.
In the light of the notion of Poincar\'e duality
to be introduced in Section 4,
we then study the question of nondegeneracy
of certain pairings in Poisson cohomology.
Our approach differs substantially from that in
\cite\evluwein.
Even in the special circumstances
of \cite\evluwein,
our pairings will not in general boil down to those
in that paper.
The requisite trace will involve
the  {\it left\/}
$(A,D_{\{\cdot,\cdot\}})$-module 
$\roman{Hom}_A(C_{\{\cdot,\cdot\}},A_{\{\cdot,\cdot\}})$
mentioned before.
A special case, not so much interesting in itself but important
for formal reasons, is that of the trivial Poisson
structure: 
Given a smooth $n$-dimensional manifold,
consider its algebra $A$ of smooth functions,
endowed with the trivial Poisson structure.
The corresponding Poisson cohomology
is just the exterior $A$-algebra
$\Lambda_AV$ over the $A$-module $V$ of smooth vector fields,
and the corresponding pairing
$$
\Lambda_A^k V \otimes_A \Lambda_A^{n-k} V
@>>> 
\Lambda_A^n V
$$
is perfect over $A$; the algebra $A$  here actually coincides with
the algebra of Casimir functions.
The nondegeneracy of this pairing cannot be phrased over
the reals as ground ring.
Our theory of Poincar\'e duality includes this example,
see (7.15) below.
Furthermore, at the end of Section 7,
we give two non-trivial examples
of nondegenerate bilinear pairings in Poisson cohomology,
the pairings being defined and nondegenerate
over the subring of Casimir elements
(and NOT over the naive ground ring, $R$).
Again nondegeneracy seems to reflect a certain
regularity property.
\smallskip
I am much indebted
to A. Weinstein for having drawn my attention
to his paper
\cite\weinsfte\ 
and to
\cite\evluwein;
in fact, the latter prompted 
me to write the present paper.
I am grateful to J. Stasheff
for a number of comments on  earlier versions
which helped improve the exposition,
and to K. Mackenzie for discussions.
Many thanks are due to the referee, for his careful reading
of an earlier version of the paper,
and for having prodded me to work out in detail
the notion of Poincar\'e duality
in terms of nondegenerate bilinear pairings.
\smallskip
It is a pleasure to dedicate this paper to my former teacher,
Beno Eckmann, at the occasion of his 80th birthday in 1997
(when this paper was written).
In my mathematical youth, when I was a student at the
ETH (Z\"urich), duality groups were a fashionable topic,
and 
the idea of duality which I then learnt from
B. Eckmann (and R. Bieri),
cf. e.~g. \cite\bierieck,
merged into the present paper;
cf. also Remark 1.6 below.

\medskip\noindent{\bf 1. Duality}\smallskip\noindent  
Let $R$ be a commutative ring,
$U$ an $R$-algebra,
and let $A$ be a left $U$-module.
We shall consider the ordinary Ext- and Tor groups 
in the category of $U$-modules.
Our convention is that
(i) $\roman {Tor}_*^U(N,M)$ is defined for
a right $U$-module $N$ and a left $U$-module $M$, and that
(ii) $\roman {Ext}^*_U(M_1,M_2)$,
without further comment,
refers to two left
$U$-modules $M_1$ and $M_2$;
occasionally we shall also consider
$\roman {Ext}^*_U(N_1,N_2)$ for two
right $U$-modules
$N_1$ and $N_2$ but this will then always be indicated explicitly.
We shall say that $U$ and $A$ {\it satisfy duality in dimension \/} $n$,
with {\it dualizing module\/} $C$, 
provided there is a {\it right\/} $U$-module
$C$ such that one has {\it natural\/} isomorphisms
$$
\roman {Ext}^{k}_U(A,M)
\cong
\roman {Tor}_{n-k}^U(C,M)
\tag1.1
$$
for all non-negative integers $k$ and all
left $U$-modules $M$.
Likewise we shall say that $U$ and $A$ {\it satisfy 
inverse
duality in dimension \/} $n$,
with {\it dualizing module\/} $D$, 
provided there is a {\it right\/} $U$-module
$D$ such that one has {\it natural\/} isomorphisms
$$
\roman {Tor}_{k}^U(N,A)
\cong
\roman {Ext}^{n-k}_U(D,N)
\tag1.2
$$
for all non-negative integers $k$ and all
right $U$-modules $N$,
where
$\roman {Ext}^{*}_U(\cdot,\cdot)$ is taken in the category
of right $U$-modules.
Here the convention is that 
$\roman {Ext}^{*}(\cdot,\cdot)$
and
$\roman {Tor}_{*}(\cdot,\cdot)$
are zero in negative degrees.
\smallskip
The algebra $U$ is a bimodule over itself in the usual way, and
the groups
$\roman {Ext}^{*}_U(A,U)$
(taken in the category of left $U$-modules)
inherit right $U$-module structures from the right $U$-module structure
on the second variable $U$ which remains free when
$\roman {Ext}^{*}_U(A,U)$ is taken.
Recall that a left $U$-module $M$ is said to have 
($U$-){\it dimension\/} $n$
provided
$\roman {Ext}^{n}_U(M,M')$
is non-zero for some left $U$-module $M'$
and
$\roman {Ext}^{k+n}_U(M,M'')$
is zero for every $k>0$
and every left $U$-module $M''$;
accordingly we can talk about the ($U$-)dimension
of a right $U$-module.
We shall say that
a projective resolution
$\varepsilon \colon K \to M$ 
of a
left $U$-module $M$ 
in the category of left $U$-modules
is {\it finite of length\/} $n$
provided
each $K_j$ is a finitely generated $U$-module
and $K_j$ is zero for $j > n$;
likewise we can talk about
{\it finite length\/} $n$ projective resolutions
in the category of right $U$-modules.

\proclaim{Proposition 1.3}
Suppose that $U$ and $A$ satisfy duality
in dimension  $n$,
with dualizing module $C$.
Then the following hold.
\roster
\item"(i)"
$\roman {Ext}^{n}_U(A,U) \cong C$
as right $U$-modules, and
$\roman {Ext}^{k}_U(A,U) = 0$
for $k \ne n, \ k \geq 0$.
\item"(ii)"
As a left $U$-module, $A$ has dimension $n$.
\item"(iii)"
As a left $U$-module, $A$ has a finite projective resolution
of length $n$.
\endroster
\endproclaim

\demo{Proof}
Taking $M = U$
in (1.1), we find
$$
\roman{Ext}_U^k(A,U) \cong
\roman {Tor}_{n-k}^U(C,U)
= \cases 0,& \quad k \ne n
\\
C,& \quad k =n .
\endcases
$$
Hence (i) holds, and (ii) is clearly also true.
Finally, we note that
the functor
$\roman {Tor}_{*}^U(C,\cdot)$
commutes with direct limits whence,
by duality,
the functor $\roman{Ext}_U^k(A, \cdot)$
commutes with direct limits
for every $k$.
In view of the corollary of Theorem 1 of \cite\kbrownfo,
since the functor $\roman{Ext}_U^k(A, \cdot)$
commutes with direct limits
and since $A$ has $U$-dimension $n$,
we conclude  that
$A$ has
a finite projective resolution of length $n$. \qed
\enddemo

\proclaim{Proposition 1.4}
Suppose that the left $U$-module $A$ has dimension $n$
in such a way that
{\rm (i)}
it has a finite projective resolution
of length $n$, and that
{\rm (ii)}
$\roman {Ext}^{k}_U(A,U) = 0$
for $0 \leq k <n$.
Then
$U$ and $A$ satisfy duality
and inverse duality
in dimension  $n$,
with (the same)  dualizing module 
$C=\roman {Ext}^n_U(A,U)$.
\endproclaim

\demo{Proof}
Let $\varepsilon \colon K \to A$
be a finite projective resolution of $A$
of length $n$ in the category of left $U$-modules.
Then $K^* = \roman{Hom}_U(K,U)$ is a finite projective resolution
of 
$C=\roman {Ext}^{n}_U(A,U)$
in the category of right $U$-modules.
Hence, for every 
left $U$-module $M$,
we have an isomorphism
$$
\Phi \colon K^* \otimes_U M @>>> \roman {Hom}_U(K,M)
\tag1.4.1
$$
which, for every non-negative integer $k$,
induces
$$
\phi\colon
\roman {Tor}_{n-k}^U(C,M)
@>>>
\roman {Ext}^{k}_U(A,M),
$$
naturally in $M$.
Likewise, for
every
right $U$-module $N$, we have
an isomorphism
$$
\Psi
\colon
N \otimes_U K @>>> \roman {Hom}_U(K^*,N)
\tag1.4.2
$$
which,  for every non-negative integer $k$,
yields
$$
\psi
\colon
\roman {Tor}_{k}^U(N,A)
@>>>
\roman {Ext}^{n-k}_U(C,N),
$$
naturally in $N$. \qed
\enddemo

We note that, when duality holds in dimension $n>0$,
the dualizing module is obviously non-zero
since otherwise the $U$-module $A$
could not have stricly positive projective dimension;
when duality holds in dimension $n=0$,
it amounts to
$\roman{Hom}_U(A,M) \cong C \otimes _UM$,
in particular,
$\roman{Hom}_U(A,A) \cong C \otimes _UA$.
Since
$\roman{Hom}_U(A,A)$
contains at least a copy of the ground ring $R$,
the dualizing module $C$  cannot then be zero either.

\smallskip\noindent
{\smc Remark 1.5.}
Let $\fra g$ be a Lie algebra which,
as a module over the
ground ring $R$,
is finitely
generated and projective
of (constant) rank $n$,
let $U $ be its universal algebra $U\fra g$,
and let $A=R$, with trivial 
$\fra g$-module and hence
trivial (unital) $U$-module structure.
Then the ordinary duality isomorphisms in the homology and cohomology
of $\fra g$,
cf. \cite\cartanei\ (Exercise 15 on p. 288),
are precisely of the kind (1.1) and (1.2),
and the dualizing module $C = \roman H^n(\fra g,U)$
has the form $\Lambda^n \fra g^*$.
See Example 4.5 below for details.
In the next section we shall show that more generally,
certain Lie-Rinehart algebras
have similar duality properties.
\smallskip\noindent
{\smc Remark 1.6.}
Let  $G$ be a discrete group, $A=R$, and $U = RG$,
  the group ring of $G$. The group $G$
is called a {\it duality group\/} or {\it inverse
duality group\/}
provided there are 
natural isomorphisms of the kind (1.1) and (1.2), respectively,
with $C$ flat and $D$ projective as $R$-modules.
When $G$ is a duality group,
$C$ is isomorphic to
the highest non-zero cohomology group $\roman H^{\roman{top}}(G,RG)$ and
when it is an inverse duality group,
$D \cong \roman H^{\roman{top}}(G,RG)$.
See \cite\bieriboo\ and \cite\bierieck\ 
for details.
\smallskip\noindent
{\smc Remark 1.7.}
It may be shown that,
when  $U$ and $A$ satisfy inverse duality
in dimension  $n$,
with dualizing module $D$,
the statements (i)-(iii)
of Proposition 1.3 hold as well, with $C$ being replaced by $D$.
We refrain from spelling out details.

\medskip\noindent{\bf 2. Lie-Rinehart algebras}
\smallskip\noindent
We now suppose that $A$ is a commutative $R$-algebra.
Let $L$ be
an $(R,A)$-Lie algebra.
Recall that
its {\it universal algebra\/}
${(U(A,L),\iota_L,\iota_A)}$
is an $R$-algebra $U(A,L)$ together with a morphism
${
\iota_A
\colon 
A
\longrightarrow
U(A,L)
}$
of $R$-algebras
and
a morphism
${
\iota_L
\colon 
L
\longrightarrow
U(A,L)
}$
of Lie algebras over $R$
having the properties
$$
\iota_A(a)\iota_L(\alpha) 
= \iota_L(a\,\alpha),\quad
\iota_L(\alpha)\iota_A(a) - \iota_A(a)\iota_L(\alpha) 
= \iota_A(\alpha(a)),
$$
and
${(U(A,L),\iota_L,\iota_A)}$
is {\it universal\/} among triples
${(B,\phi_L,\phi_A)}$
having these properties.
More precisely:
{\sl Given 
\newline\noindent
{\rm (i)}\phantom{ii} another
 $R$-algebra $B$, viewed at the same time as a Lie algebra over
$R$,
\newline
\noindent
{\rm (ii)}\phantom{i} a morphism
$
\phi_L 
\colon
L
\longrightarrow
B
$ of Lie algebras over $R$, and
\newline
\noindent
{\rm (iii)} a morphism
$
\phi_A 
\colon
A
\longrightarrow
B
$
of $R$-algebras,
\newline\noindent
\noindent
so that, for ${\alpha \in L, a \in A}$,
$$
\align
\phi_A(a)\phi_L(\alpha) &= \phi_L(a\,\alpha),
\\
\phi_L(\alpha)\phi_A(a) - \phi_A(a)\phi_L(\alpha) &= \phi_A(\alpha(a)),
\endalign
$$
there is a unique morphism
${
\Phi 
\colon
U(A,L)
\longrightarrow
B
}$
of $R$-algebras
so that
$
\Phi\,\iota_A = \phi_B
$ 
and
$
\Phi\,\iota_L = \phi_L$.
}
Often we shall simply write
$U$ instead of 
${(U(A,L),\iota_L,\iota_A)}$.
See (1.6) of \cite\poiscoho\ 
for more details.
We only mention that $U$ is determined up to isomorphism
by the universal property as usual.
For example,
when
$A$ is the algebra of smooth functions on a smooth
manifold $B$ and  $L$  the Lie algebra of 
smooth vector fields
on $B$, then $U(A,L)$ is the {\it algebra of 
(globally defined)
differential operators
on\/} $B$.
\smallskip
Suppose that $A$ is endowed with the obvious
left $U$-module structure induced by the left $L$-action on $A$.
Our aim is to study the notions of duality and inverse duality
for this special case.
\smallskip
Recall that an $A$-module $M$ which is also a left $L$-module
is called a {\it left\/} $(A,L)$-module
provided
$$
\align
\alpha (a x) &=  \alpha(a) x + a \alpha (x)
\tag2.0.1
\\
(a \alpha) (x) &=  a (\alpha (x))
\tag2.0.2
\endalign
$$
where  $a \in A, \ x \in M,\ \alpha \in L$.
An $A$-module $N$ which is also a right $L$-module,
the action being written
$(x,\alpha) \mapsto x \alpha$,
is called a {\it right\/} $(A,L)$-module
provided
$$
\align
 (a x) \alpha &=  a  (x \alpha)  - \alpha(a) x  
\tag2.0.3
\\
 x (a \alpha) &=  a  (x \alpha)  - \alpha(a) x
\tag2.0.4
\endalign
$$
where  $a \in A, \ x \in M,\ \alpha \in L$.
Then left- and right $(A,L)$-module structures
correspond to
left and right $U(A,L)$-module structures, respectively,
and vice versa.
The formula (2.0.4) might look somewhat puzzling at first glance;
it is {\it not\/} the expected replica of (2.0.2). In view of the
associativity law, there is only one consistent way to evaluate the
expression $x(a\alpha)$ on the left-hand side of (2.0.4), though.
\smallskip
We now list a few facts
related with the behaviour of
$(A,L)$-modules under the operations \lq\lq $\cdot \otimes _A \cdot$\rq\rq\ 
and
\lq\lq $\roman{Hom}_A(\cdot,\cdot)$\rq\rq.
\smallskip
\roster
\item"(2.1)"
Given two left
$(A,L)$-modules
$M_1$ and $M_2$, 
the customary formula 
$$
\alpha (m_1 \otimes m_2) =
(\alpha m_1) \otimes m_2 + m_1 \otimes (\alpha m_2),
\quad
m_1 \in M_1, \ m_2 \in M_2, \ \alpha \in L,
$$
endows
their tensor product
$M_1 \otimes _A M_2$
with a  left
$(A,L)$-module structure.
\item"(2.2)"
Given two left
$(A,L)$-modules
$M_1$ and $M_2$, 
the customary formula 
$$
(\alpha \phi)(m)
=
\alpha (\phi m) - \phi(\alpha m),
\quad
m \in M_1,\  \alpha \in L,\ \phi \in \roman{Hom}_A(M_1, M_2),
$$
endows
the $A$-module $\roman{Hom}_A(M_1, M_2)$
with  a left
$(A,L)$-module structure.
\item"(2.3)"
Given two right
$(A,L)$-modules
$N_1$ and $N_2$, 
the customary formula 
$$
(\alpha \phi)(n)
=
\phi(n \alpha) - (\phi n) \alpha,
\quad
n \in N_1,\  \alpha \in L,\ \phi \in \roman{Hom}_A(N_1, N_2),
$$
endows
the $A$-module $\roman{Hom}_A(N_1, N_2)$
with a  left (!)
$(A,L)$-module structure.
\item"(2.4)"
Given left and 
right
$(A,L)$-modules
$M$ and $N$, respectively,
the  formula 
$$
(n \otimes m) \alpha=
(n \alpha) \otimes m - n  \otimes (\alpha m),
\quad
m \in M, \ n \in N, \ \alpha \in L,
$$
endows
their tensor product
$N \otimes _A M$
with a   right (!)
$(A,L)$-module structure.
\item"(2.5)"
Given left and 
right
$(A,L)$-modules
$M$ and $N$, respectively,
the  formula 
$$
(\phi\alpha)m=
(\phi m)\alpha -  \phi(\alpha m),
\quad
m \in M, \ \alpha \in L, \ \phi \in \roman{Hom}_A(M,N),
$$
endows
$\roman{Hom}_A(M,N)$
with a  right (!)
$(A,L)$-module structure.
\endroster
\smallskip
We shall say that $L$ {\it satisfies duality in dimension \/} $n$,
with {\it dualizing module\/} $C$, 
provided  $U$ and $A$
satisfy duality in dimension $n$,
with dualizing module $C$; likewise 
we shall say that $L$ {\it satisfies inverse duality in dimension \/} $n$,
with {\it dualizing module\/} $D$, 
provided  $U$ and $A$
satisfy inverse duality in dimension $n$,
with  dualizing module $D$. 
\smallskip
Recall \cite\rinehart\ that
the classical
complex 
of
multilinear alternating forms
extends 
from Lie algebras to
Lie-Rinehart algebras; for a left $(A,L)$-module $M$,
we denote by
$$
(\roman{Alt}_A(L,M),d)
\tag2.6
$$
the resulting chain complex (over $R$)
of $M$-valued $A$-multilinear alternating forms,
where $d$ is the corresponding 
Cartan-Chevalley-Eilenberg operator.
We only note that the fact that
this operator sends $A$-multilinear forms
to $A$-multilinear forms
is not obvious and requires proof;
see \cite\rinehart\ and what is said below.
Further, the ordinary constructions of 
{\it contraction\/}
$$
i \colon
L \otimes _R \roman{Alt}_A(L,M)
@>>>
\roman{Alt}_A(L,M)
\tag2.6.1
$$
and 
{\it Lie derivative\/}
$$
\lambda \colon
L \otimes _R (\roman{Alt}_A(L,M),d)
@>>>
(\roman{Alt}_A(L,M),d)
\tag2.6.2
$$
carry over as well \cite\rinehart;
for   $\alpha \in L$,
we  write 
$$
i_{\alpha}
\colon
\roman{Alt}_A(L,M)
@>>>
\roman{Alt}_A(L,M)
$$
and
$$
\lambda_\alpha
\colon
(\roman{Alt}_A(L,M),d)
@>>>
(\roman{Alt}_A(L,M),d),
$$
respectively.
The operation of Lie derivative endows
$(\roman{Alt}_A(L,M),d)$
with   a {\it left\/} $L$-module
structure
({\it not\/} with that of a left $(A,L)$-module) and,
as is common for ordinary Lie algebras,
the negative
of the Lie derivative yields a  {\it right\/}
$L$-module
structure on $(\roman{Alt}_A(L,M),d)$
({\it not\/} that of a right $(A,L)$-module).
These operations
satisfy the customary formulas
$$
\align
\lambda_{a \alpha}\omega &= a \lambda_{\alpha}\omega + (da) \cup (i_\alpha 
\omega),
\tag2.7.1
\\
\lambda_{\alpha} &= d i_\alpha +   i_\alpha d,
\tag2.7.2
\\
i_\alpha(\omega_1 \cup \omega_2)
&=
i_\alpha(\omega_1) \cup \omega_2
+ (-1)^{\roman {deg}(\omega_1)}
\omega_1 \cup i_\alpha(\omega_2),
\tag2.7.3
\endalign
$$
where $\alpha \in L,\
a \in A, \ \omega,\omega_1,\omega_2 \in
\roman{Alt}_A(L,M)$.
\smallskip

Recall that, for a finitely generated projective $A$-module,
its rank to be constant means that it is the same for every
prime ideal of $A$.
For example, when $A$ is the algebra
of smooth functions
on a smooth manifold $W$
and $L$ the $(\Bobb R,A)$-Lie algebra of smooth vector fields
on $W$,
the  requirement that $L$ be
of constant rank
amounts to the connectedness of $W$.

\proclaim{Proposition 2.8} 
Suppose that, as an $A$-module,
the $(R,A)$-Lie algebra $L$
is finitely generated projective of constant rank $n$.
For every left $(A,L)$-module $M$,
the formula
$$
\phi \alpha = -\lambda_\alpha(\phi),
\quad
\phi \in
\roman{Hom}_A(\Lambda^n_A L,M),
\ \alpha \in L,
\tag2.8.1
$$
then endows
$\roman{Hom}_A(\Lambda^n_A L,M)$
with  a right {\rm (!)}
$(A,L)$-module structure.
\endproclaim

\demo{Proof}
For
$a \in A,\  \alpha \in L,\  \phi \in
\roman{Hom}_A(\Lambda^n_A L,M) =
\roman{Alt}^n_A(L,M)$, in view of (2.7.1) and (2.7.3), we have
$$
\lambda_{a \alpha}\phi= a \lambda_{\alpha}\omega + (da) 
\cup (i_\alpha \phi)
$$
and
$$
i_\alpha ((da) \cup \phi)
=
(i_\alpha (da))\cup \phi
-
(da)\cup (i_\alpha \phi)
=
\alpha (a) \phi
-
(da)\cup (i_\alpha \phi).
$$
Since $\phi$ is in the top dimension,
$(da) \cup \phi$ is zero whence
$$
(da)\cup (i_\alpha \phi)=\alpha (a) \phi.
$$
Thus
$
\lambda_{a \alpha}\phi= a \lambda_{\alpha}\omega + 
\alpha (a) \phi$
and thence
$
\phi (a\alpha) = a (\phi \alpha) -\alpha(a) \phi.
$
Consequently (2.8.1)
yields a  right
$(A,L)$-module structure
on
$\roman{Hom}_A(\Lambda^n_A L,M)$
as asserted. \qed
\enddemo
A special case of this proposition is known from
the theory of $D$-modules,
cf. (VI.3.2) and (VI.3.3) on pp. 226/227 of \cite\boredmod.
\smallskip
Under appropriate circumstances,
the chain complex  (2.6)
computes the cohomology of 
the $(R,A)$-Lie algebra
$L$
with values in $M$.
To recall what this means,
we reproduce briefly the 
Rinehart complex for $(A,L)$;
in the next section,
we shall further exploit this
Rinehart complex and in particular deduce descriptions of
the duality isomorphisms 
(to be given in (2.11) below)
in terms of appropriate
cup- and cap products:
Let
$\Lambda_AL$
be the exterior $A$-algebra of $L$, and consider the graded
left $U(A,L)$-module
$U(A,L) \otimes _A \Lambda_AL$
where $A$ acts on  $U(A,L)$ from
the right by means of the canonical map 
from $A$ to $U(A,L)$.
For $u \in U(A,L)$ and $\alpha_1, \dots, \alpha_n \in L$, let
$$
\aligned
&d(u\otimes_A (\alpha_1\wedge \dots\wedge \alpha_n))
\\
&
=
\quad\sum_{i=1}^n (-1)^{(i-1)}u\alpha_i\otimes_A 
(\alpha_1\wedge \dots\widehat{\alpha_i}\ldots\wedge \alpha_n )
\\
&\phantom{=}+\quad
\sum_{j<k} (-1)^{(j+k)}u\otimes_A 
(\lbrack \alpha_j,\alpha_k \rbrack\wedge
\alpha_1\wedge \dots\widehat{\alpha_j}\dots
\widehat{\alpha_k}\ldots\wedge \alpha_n).
\endaligned
\tag2.9.1
$$
Rinehart \cite\rinehart\ has proved that
this yields
an 
${U(A,L)}$-linear differential
$$
d
\colon
U(A,L) \otimes _A \Lambda_AL 
\longrightarrow
U(A,L) \otimes _A \Lambda_AL,
$$
that is,
$dd =0$.
The non-trivial fact to be verified here is that,
for every $u \in U(A,L)$,\ $a \in A$, and $\alpha_1, \dots, \alpha_n \in L$, 
$$
d((ua)\otimes (\alpha_1\wedge \dots\wedge \alpha_n))
=
d(u\otimes ((a\alpha_1)\wedge \dots\wedge \alpha_n))
=
\dots
=
d(u\otimes (\alpha_1\wedge \dots\wedge (a\alpha_n) )).
$$
We only mention that,
separately,
the two summands on the right-hand side of
(2.9.1) are {\it not\/}
even well defined, and only their sum 
is $A$-multilinear in the variables
$u, \alpha_1,\dots,\alpha_n$.
We refer to the resulting chain complex
$$
K(A,L) = (U(A,L) \otimes _A \Lambda_AL,d)
\tag2.9.2
$$
as the
{\it Rinehart complex\/} for $(A,L)$.
Rinehart has also proved that,
when $L$ is projective as an $A$-module,
$K(A,L)$ is in fact a projective resolution of
$A$ in the category of left $U(A,L)$-modules.
\smallskip
We now recall that, for a general
$(R,A)$-Lie algebra $L$ and
for left and right $(A,L)$-modules $M$ and $N$, respectively,
the homology and cohomology of $L$
are defined by
$$
\roman H^*(L,M) = \roman {Ext}^{*}_U(A,M)
\quad
\text{and}
\quad
\roman H_*(L,N) = \roman {Tor}_{*}^U(N,A).
\tag2.9.3
$$
When $L$ is projective as an $A$-module,
its homology and cohomology
may thus be computed from
$K(A,L)$;
hence,
given a left $(A,L)$-module $M$,
the cohomology
$\roman H^*(L,M)$
may be obtained as that of
the chain complex 
$$
(\roman{Alt}_A(L,M),d)
=
\roman{Hom}_U(K(A,L),M)
\tag2.9.4
$$
and a similar statement can be made
for the
homology
$\roman H_*(L,N)$
of $L$ with values in a right $(A,L)$-module $N$.
Notice that the description (2.9.4) of
$(\roman{Alt}_A(L,M),d)$
justifies in particular
the earlier claim
(which may also be verified directly)
that
the Cartan-Chevalley-Eilenberg operator sends $A$-multilinear forms
to $A$-multilinear forms.
(The notation \lq\lq $d$\rq\rq\ is slightly abused here.)
For an
ordinary Lie algebra $\fra g$ over $R$,
viewed as an
$(R,A)$-Lie algebra with trivial $L$-action on $A=R$,
the Rinehart complex $K(R,\fra g)$ comes down to the standard {\it Koszul
complex\/} which we denote by $K_R\fra g$
(or Koszul resolution of $R$ when $\fra g$ is projective
as an $R$-module)
familiar in ordinary Lie algebra (co)homology.
\smallskip
When the $(R,A)$-Lie algebra $L$ is finitely generated 
and projective of constant rank $n$ as an $A$-module,
its Rinehart complexe $K(A,L)$ 
is manifestly a finite projective resolution of $A$
of length 
$n$
in the category of left $U(A,L)$-modules,
and the highest non-zero term
$
K_n(A,L)
$
of the resolution has the form
$
U(A,L)
\otimes _A
\Lambda_A^n L$.

\proclaim {Theorem 2.10}
An $(R,A)$-Lie algebra $L$ 
which,
as an $A$-module, is finitely generated and projective
of constant rank $n$,
satisfies
duality 
and inverse duality in dimension  $n$,
with {\it dualizing module\/} 
$$
C_L
=\Lambda_A^n L^* 
= \roman{Hom_A}(\Lambda_A^n L, A),
$$
the requisite isomorphism between $C_L$ and
$\roman H^n(L,U)$ being induced by 
the injection
of
$\roman{Hom_A}(\Lambda_A^n L, U)$
into
$$
\roman{Hom_A}(\Lambda_A^n L, U)
\cong
\roman{Hom_U}(K_n(A,L), U)
$$
coming from the canonical inclusion
of $A$ into $U=U(A,L)$;
here
$\roman{Hom_A}(\Lambda_A^n L, A)$
is viewed
endowed
with the  right $(A,L)$-module structure {\rm (2.8.1)}.
\endproclaim

\demo{Proof}
We show that the conditions of Proposition 1.4 
are met.
To this end,
using the fact that $L$ is projective as an $A$-module,
we
exploit the Poincar\'e-Birkhoff-Witt theorem
for Lie-Rinehart algebras
\cite\rinehart\ 
and
consider the customary filtration of $U=U(A,L)$ 
coming from powers of $L$ and
having as associated graded
algebra
$\roman E^0(U)$ 
the symmetric algebra $S_AL$ on $L$ in the category of $A$-modules.
This filtration induces a spectral sequence
having
$\roman E_1 =\roman {Ext}^{*}_{S_AL}(A,S_AL)$
and converging to
$\roman {Ext}^{*}_{U}(A,U)$.
But
$\roman {Ext}^{k}_{S_AL}(A,{S_AL}) = 0$
for $0 \leq k <n$,
and
the injection
of
$\roman{Hom_A}(\Lambda_A^n L, A)$
into
$\roman{Hom_A}(\Lambda_A^n L, S_AL)$
coming from the canonical inclusion
of $A$ into $S_AL$
induces an isomorphism 
from
$\roman{Hom_A}(\Lambda_A^n L, A)$
onto
$\roman {Ext}^{n}_{S_AL}(A,{S_AL})$
whence the spectral sequence
has 
$\roman E_1 =\roman E_{\infty}$.
Hence,
for $0 \leq k <n$,
$\roman {Ext}^{k}_U(A,U) = 0$,
and
the injection
of
$\roman{Hom_A}(\Lambda_A^n L, A)$
into
$\roman{Hom_A}(\Lambda_A^n L, U)$
coming from the canonical inclusion
of $A$ into $U$
induces an isomorphism
$$
C_L
=
\roman{Hom_A}(\Lambda_A^n L, A)
@>>>
\roman H^n(L,U).
\tag2.10.1
$$
To see that this is an isomorphism
of right
$(A,L)$-modules
of the asserted kind,
let $\alpha,\alpha_1,\dots,\alpha_n \in L$
and
$\phi \in \roman{Hom_A}(\Lambda_A^n L, A)$
so that
$\phi(\alpha_1,\dots,\alpha_n) \in A$;
a straightforward calculation shows that, then,
$$
\align
(\phi \alpha)(\alpha_1,\dots,\alpha_n)
+
(\lambda_\alpha\phi)
(\alpha_1,\dots,\alpha_n)
&=
-\alpha (\phi(\alpha_1,\dots,\alpha_n))
+ \alpha \cdot (\phi(\alpha_1,\dots,\alpha_n))
\\
&=
(\phi(\alpha_1,\dots,\alpha_n)) \cdot \alpha;
\endalign
$$
here \lq\lq\,$\cdot$\,\rq\rq\  refers to the product in $U = U(A,L)$,
$\phi \alpha 
\in \roman{Hom_A}(\Lambda_A^n L, A)$
denotes the result of
the operation
on $\phi \in \roman{Hom_A}(\Lambda_A^n L, A)$
with $\alpha \in L$ from the right,
$\alpha (\phi(\alpha_1,\dots,\alpha_n)) \in A$
that of
the operation
on 
$\phi(\alpha_1,\dots,\alpha_n) \in A$
with $\alpha$ from the left,
and the equality sign means \lq identity in
$\roman{Hom_A}(\Lambda_A^n L, U)$\rq,
$\roman{Hom_A}(\Lambda_A^n L, A)$
being viewed as
a subspace thereof
via the obvious injection of $A$ into $U$.
Since
we are in the top dimension,
the operation $\lambda_\alpha$ equals $di_\alpha$ whence
$$
\lambda_\alpha\phi
=
di_\alpha \phi
$$
is a coboundary.
Thus
the cocycles
$\phi \alpha$
and
$$
(\alpha_1,\dots,\alpha_n) \mapsto 
(\phi(\alpha_1,\dots,\alpha_n)) \cdot \alpha
$$
are cohomologous.
Hence (2.10.1)
is even an isomorphism of
right
$(A,L)$-modules.
In view of  Proposition 1.4,
this proves the claim. \qed 
\enddemo

Theorem 2.10 suggests the following definition:
An $(R,A)$-Lie algebra $L$ which, as an $A$-module,
is finitely generated and projective
of constant rank, will henceforth be referred to as a
{\it duality\/} $(R,A)$-{\it Lie algebra\/} and 
the pair $(A,L)$ will be called a 
{\it duality Lie-Rinehart algebra\/}.
The rank of
a duality $(R,A)$-Lie algebra $L$,
viewed as an $A$-module,
will be referred to as the
{\it rank\/} of $L$ and, likewise,
we shall talk about the {\it rank\/} of
a duality Lie-Rinehart algebra.

\proclaim{Corollary 2.11}
Given a duality
$(R,A)$-Lie algebra
$L$
of rank $n$,
there are  natural isomorphisms
$$
\phi
\colon
\roman H_k(L, C_L\otimes_A M)
@>>>
\roman H^{n-k}(L,M)
\tag2.11.1
$$
for all non-negative integers $k$ and all
left $(A,L)$-modules $M$ and, furthermore,
natural isomorphisms
$$
\psi
\colon
\roman H_k(L,N)
@>>>
\roman H^{n-k}(L,\roman{Hom}_A(C_L,N))
\tag2.11.2
$$
for all non-negative integers $k$ and all
right $(A,L)$-modules $N$.
\endproclaim
Here
$C_L\otimes_A M$
and
$\roman{Hom}_A(C_L,N)$
carry the corresponding
right- and left
$(A,L)$-module structures
explained in (2.4) and (2.3), respectively.
\smallskip
A duality isomorphism
of the kind (2.11.1)
may be found in \cite{\chemlone\ (5.2.2)}.
The additional information
provided for by Theorem 2.10 does not seem to 
be in the literature, though.
Occasionally we refer to duality isomorphisms
of the kind (2.11.1) and (2.11.2)
as {\it naive duality\/}.
\smallskip
\noindent
{\smc Example 2.12.}
Let $\fra g$ be a Lie algebra over $R$ which, as an $R$-module,
is supposed
to be finitely generated and projective of constant rank $n$ (say),
let $A$ be a
commutative $R$-algebra,
 and suppose that $\fra g$ acts on $A$
 by derivations.
In view of the defining properties (0.3)
for an $(R,A)$-Lie algebra,
the Lie bracket on $\fra g$ and the $\fra g$-action on $A$
induce a bracket
on $L= A \otimes_R \fra g$
which,
together with the obvious left $A$-module structure on $L$,
turns
$L$ into an  $(R,A)$-Lie algebra.
As an $A$-module, the dualizing module
$C_L$ of $L$ is plainly 
isomorphic to $A \otimes _R\Lambda^n\fra g^*$.
Moreover, the universal algebra $U(A,L)$ may be written
in the form $A \otimes _R U\fra g$,
and its algebra structure is 
given by
$$
(a \otimes _R 1)(1 \otimes_R x) = a \otimes_R x,
\quad
(1 \otimes_R x)(a \otimes _R 1)
=
a \otimes_R x + (x(a))\otimes_R 1,
$$
where $a \in A$ and $x \in \fra g$.
Further, the duality isomorphisms
of the kind (2.11.1) and (2.11.2)
obtain.
Now,
$(R,\fra g)$ being considered as a Lie-Rinehart algebra as well,
with trivial $\fra g$-action
on $R$ (which coincides with $A$ in this case),
we denote the corresponding dualizing module by $C_{\fra g}$ 
($\cong \Lambda^n\fra g^*$).
When $A$ is different from $R$
(and comes with a non-trivial
$\fra g$-action), 
as an $A$-module,
$C_L$ is  plainly isomorphic to $A \otimes _R C_{\fra g}$, 
and 
on this isomorphic image of $C_L$,
the 
right $(A,L)$-module structure 
is 
induced by the obvious right $\fra g$-module structure on
$A \otimes _R C_{\fra g}$
and
may thus be described by the formula
$$
(b \otimes_R \phi)(a \otimes_R x) =
-(a x(b))\otimes_R \phi - a \otimes_R(\lambda_x \phi) -
(x(a) b) \otimes_R \phi;
\tag2.12.1
$$
here
$a,b \in A,\ x \in \fra g,\  \phi \in 
C_{\fra g}\cong \Lambda^n\fra g^*$ and,
for $y \in 
\Lambda^n\fra g$, the value
$(\lambda_x \phi)(y)$
is given by
$
(\lambda_x \phi)(y)
=- \phi (x(y))$.
The formula (2.12.1) may look a bit odd but  is
forced by
(2.0.3) and (2.0.4).
The duality isomorphisms for $(A,L)$
now boil down to those for the ordinary
Lie algebra $\fra g$. 
In fact,
given left- and right
$(A,L)$-modules $M$ and $N$,
respectively,
the duality isomorphisms (2.11.1)
and (2.11.2)
come down to
$$
\roman H_k(L, C_L\otimes_A M)
\cong
\roman H_k(\fra g, C_{\fra g}\otimes_R M)
@>{\phi}>>
\roman H^{n-k}(\fra g,M)
\cong
\roman H^{n-k}(L,M)
\tag2.12.2
$$
and
$$
\roman H_k(L,N)
\cong
 \roman H_k(\fra g,N)
@>{\psi}>>
\roman H^{n-k}(\fra g,\roman{Hom}_R(C_{\fra g},N))
\cong
\roman H^{n-k}(L,\roman{Hom}_A(C_L,N)),
\tag2.12.3
$$
respectively.
Cf. Remark 1.5 above.
For later reference we mention that,
given a left $\fra g$-module $\fra m$,
the induced $A$-module $A \otimes_R \fra m$
inherits an obvious left $(A,L)$-module structure
given by the formula
$$
(a\otimes_Rx)(b\otimes_R y) = (a x(b)) \otimes_R y + (a b) \otimes_R (xy),
\quad a, b \in A,\ x \in \fra g,\ y\in \fra m,
\tag2.12.4
$$
which is in fact forced by (2.0.1) and (2.0.2).
Furthermore,
given any left $(A,L)$-module $M$, the formula
$$
y (a\otimes_Rx)= 
-(a\otimes_Rx) y
-(x(a)) y,
\quad a \in A,\ x \in \fra g,\ y\in M,
\tag2.12.5
$$
yields a right $(A,L)$-module structure on
$M$;
this is forced by (2.0.3) and (2.0.4).
For a general Lie-Rinehart algebra,
we cannot naively pass from left- to right $(A,L)$-module
structures in such a naive fashion, though.

\medskip\noindent{\bf 3. Multiplicative structures}\smallskip\noindent
In the definition of duality, we did not require that the duality isomorphisms
commute with connecting homomorphisms nor with maps
induced by morphisms in the $(R,A)$-Lie algebra argument.
However, 
Theorem 3.7 below
shows that duality isomorphisms can be given by a 
cap-product,
and these are, of course, natural in any reasonable sense.
\smallskip
As before, let $L$ be an $(R,A)$-Lie algebra.
Recall from Section 1 of \cite\poiscoho,
cf. also (1.6) in \cite\extensta,
that,
given a pairing
$M_1 \otimes_A M_2 \to M$
of left $(A,L)$-modules,
the standard shuffle multiplication of alternating maps
induces a pairing
$$
\cup\colon
\roman{Alt}_A(L,M_1)
\otimes_R
\roman{Alt}_A(L,M_2)
@>>>
\roman{Alt}_A(L,M)
\tag3.1
$$
of chain complexes over the ground ring $R$
(see (1.5') in \cite\poiscoho).
When $L$ is projective as an $A$-module so that
the homology of $\roman{Alt}_A(L,M)$
etc. gives $\roman H^*(L,M)$ etc.,
the pairing (3.1) induces a pairing
in the cohomology of $L$ generalizing the ordinary {\it cup\/}
pairing in Lie algebra cohomology
and we refer to the resulting pairing as
{\it cup pairing\/} as well,
written
$$
\cup\colon
\roman H^*(L,M_1)
\otimes_R
\roman H^*(L,M_2)
@>>>
\roman H^*(L,M)
\tag3.2
$$
(with an abuse of the notation $\cup$).
This pairing also generalizes the ordinary multiplicative
structure in de Rham
cohomology.
\smallskip
For intelligibility, we recall an explicit description;
it is not completely simple because
the distinction between
graded $A$-objects and differential graded $R$-modules
due to the in general non-trivial action of $L$ on $A$
persists throughout:
The ordinary diagonal map
$$
\Delta \colon
\Lambda_RL
\longrightarrow
\Lambda_RL\otimes \Lambda_RL
\tag3.1.1
$$
determined by
$$
\Delta(v) = v \otimes 1 + 1 \otimes v,\quad v \in L,
\tag3.1.2
$$
makes 
the graded exterior $R$-algebra
$\Lambda_RL$   
into a graded commutative and graded
cocommutative $R$-Hopf algebra and hence in particular endows
$\Lambda_RL$
with 
a graded cocommutative coalgebra structure;
see e.~g. {\smc Mac Lane}~\cite\maclaboo\ for details.
This diagonal map is also referred to as
{\it shuffle coproduct\/}
or {\it shuffle diagonal\/}.
Since 
the property of being a Hopf algebra implies in particular
that its diagonal map is multiplicative, the assignment (3.1.2)
in fact completely determines (3.1.1).
Explicitly, 
given
$x_1,\dots,x_{p+q} \in L$,
the value 
$
\Delta(x_1 \wedge\dots\wedge x_{p+q})$ in 
$\Lambda_RL$ 
is given by the formula
$$
\aligned
\Delta (x_1\wedge\dots\wedge x_{p+q})&
\\ 
=
\sum_{\sigma}\roman{sign}(\sigma)
&(x_{\sigma(1)}\wedge\dots\wedge x_{\sigma(p)})
\otimes(x_{\sigma(p+1)}\wedge\dots\wedge x_{\sigma(p+q)}),
\endaligned
\tag3.1.3
$$
where
$\sigma$
runs through $(p,q)$-shuffles
and where $\roman{sign}(\sigma)$ refers to the sign of $\sigma$.
Formally the same construction, but in the category of
$A$-modules rather than in that of $R$-modules,
yields a diagonal map
$$
\Delta_A \colon
\Lambda_AL
\longrightarrow
\Lambda_AL\otimes_A \Lambda_AL
\tag3.1.4
$$
on
(on the exterior $A$-algebra)
$\Lambda_AL$
determined by
$$
\Delta_A(v) = v \otimes_A 1 + 1 \otimes_A v,\quad v \in L,
\tag3.1.5
$$
and this diagonal map makes
the graded exterior $A$-algebra
$\Lambda_AL$ into
a graded commutative and graded
cocommutative $A$-Hopf algebra and hence in particular 
endows
$\Lambda_AL$ with 
a graded cocommutative $A$-coalgebra structure.
Moreover the
two diagonal maps (3.1.1) and (3.1.4) are 
compatible in the sense that the diagram
$$
\CD
\Lambda_RL
@>{\Delta}>>
\Lambda_RL\otimes \Lambda_RL
\\
@VVV
@VVV
\\
\Lambda_AL
@>{\Delta_A}>>
\Lambda_AL\otimes_A \Lambda_AL
\endCD
\tag3.1.6
$$
is commutative,
the unlabelled arrows being the obvious maps.
Dualization now yields 
the pairing (3.1):
Given 
$\alpha \in \roman{Alt}^p_A(L,M_1)$
and $\beta \in \roman{Alt}^q_A(L,M_2)$
and,
furthermore, arbitrary $x_1,\dots,x_{p+q} \in L$,
the value
$(\alpha \cup \beta)(x_1,\dots,x_{p+q}) $
in $\roman{Alt}^{p+q}_A(L,M)$
is given by
the  explicit
expression
$$
\aligned
(\alpha \cup \beta)&(x_1,\dots,x_{p+q}) 
\\
&= (-1)^{|\alpha||\beta|}
\sum_{\sigma}\roman{sign}(\sigma)
\mu_A(\alpha(x_{\sigma(1)},\dots,x_{\sigma(p)})
\otimes_A \beta(x_{\sigma(p+1)},\dots,x_{\sigma(p+q)})),
\endaligned
$$
where $\mu_A\colon M_1 \otimes_A M_2 \to M$
denotes the given $(A,L)$-module pairing
from $M_1 \otimes_A M_2$ to $M$.
Without the differentials, we would also obtain
a pairing
$$
\roman{Alt}_A(L,M_1)
\otimes_A
\roman{Alt}_A(L,M_2)
@>>>
\roman{Alt}_A(L,M)
$$
of graded $A$-modules
but for reasons already explained our pairing of
primary interest is (3.1).
\smallskip
Likewise, for left- and right $(A,L)$-modules
$M$ and $N$, respectively, the standard 
construction still yields a 
cap-pairing.
Notice that,
for reasons already hinted at,
we must be a bit circumspect here:
the standard construction of multiplicative structures
in e.~g. the cohomology of a Lie algebra
in the ordinary sense (over the ground ring $R$)
involves the diagonal map on its universal algebra;
under our more general circumstances
of a general $(R,A)$-Lie algebra $L$
the universal algebra
$U(A,L)$
will not in general inherit a  
Hopf algebra 
structure (over $R$)
unless $L$ is a Lie algebra over $A$
in the ordinary sense, that is, when $L$ acts trivially on $A$
and $A$ coincides with the ground ring $R$.
In fact, when the algebra $A$ is different
from the ground ring and $L$ acts non-trivially on $A$,
the question of Hopf algebra structure 
on $U(A,L)$
is only well posed 
when $A$ itself comes with a Hopf algebra structure
compatible with the $L$-action; this happens to be the case,
for example, when $A$ is 
the algebra of algebraic functions on
an algebraic group, e.~g. on a compact Lie group over the reals,
and when $L$ is the Lie algebra of derivations of such an $A$.
We refrain from spelling out details, since we shall not need them.

\proclaim{Theorem 3.3}
When $L$ is projective as an $A$-module,
for left- and right $(A,L)$-modules
$M$ and $N$, respectively, the customary formula
involving the requisite shuffle map yields a cap pairing
$$
\cap\colon 
\roman H_{\ell}(L,N) 
\otimes_R
\roman H^k(L,M) 
@>>>
\roman H_{\ell-k}(L,N\otimes_A M)
\tag3.3.1
$$
which is natural in terms of the data.
\endproclaim

To prepare for the proof, we recall the notion of cap product,
in a form tailored to our purposes:
Let $\Lambda$ be a (differential) graded $R$-coalgebra, and let
$M$ and $N$ be $R$-modules
(we could allow $M$ and $N$ to be $R$-chain complexes but we do
not need this greater generality).
Given $\phi \in \roman{Hom}_R(\Lambda,M)$, the morphism
$$
\cdot\cap \phi
\colon
N\otimes_R \Lambda
@>>>
N\otimes_R M\otimes_R \Lambda
$$
is defined as the composite
$$
N\otimes_R \Lambda
@>{N\otimes_R \Delta}>>
N\otimes_R \Lambda \otimes_R \Lambda
@>{N\otimes_R\phi \otimes_R \Lambda}>>
N\otimes_R M\otimes_R \Lambda,
$$
where the identity morphism on an object is denoted by the same
symbol as that object.
The resulting pairing
$$
\cap \colon
N\otimes_R \Lambda \otimes _R \roman{Hom}_R(\Lambda,M)
@>>>
N\otimes_R M\otimes_R \Lambda
\tag3.3.2
$$
is a version of the ordinary {\it cap pairing\/} in differential
homological algebra.
See e.~g. (2.3) in \cite\gugenhtw.
We note that, in order to get differentials right (whenever differentials
come into play),
the Eilenberg-Koszul convention should systematically be in force:
Whenever two graded objects $a$ and $b$ (say) are interchanged,
the sign $(-1)^{|a| |b|}$
should be added.
In particular,
for an ordinary $R$-Lie algebra $\fra g$,
given left- and right $\fra g$-modules
$M$ and $N$, respectively,
taking
$\Lambda = \Lambda_R \fra g$,
with the shuffle diagonal
(3.1.1),
the pairing (3.3.2), viewed as one 
of graded $R$-modules, may also be written
$$
\cap\colon 
(N\otimes_{U\fra g} K_R\fra g) 
\otimes_R
\roman {Alt}_R(\fra g,M)
@>>>
(N\otimes_RM)\otimes_{U\fra g} K_R\fra g
$$
where $K_R\fra g$
denotes the Koszul complex
for $\fra g$ in the category
of $R$-modules, cf. what has been said in the previous section
($K_R \fra g$ coincides with the Rinehart complex $K(R,\fra g)$, cf. (2.9.2)).
In this form, the pairing is actually compatible with the requisite
differentials
(coming from the Koszul complex  $K_R\fra g$).
When $\fra g$ is projective as an $R$-module,
this pairing computes the customary cap pairing
$$
\cap\colon 
\roman H_{\ell}(\fra g,N) 
\otimes_R
\roman H^k(\fra g,M) 
@>>>
\roman H_{\ell-k}(\fra g,N\otimes_R M)
$$
in the (co)homology of $\fra g$
which, in turn,
is induced by
the standard Hopf algebra structure 
on the universal algebra $U\fra g$ of $\fra g$
obtained when the elements of $\fra g$ are required to be primitive.
(When $\fra g$ is not projective, these statements are
true for the corresponding relative homological algebra.)
While
this Hopf algebra structure is lurking behind,
there is no need for us
to make it explicit,
and we refrain from spelling out details.
Our crucial observation is that,
though in general $U(A,L)$ no longer inherits a Hopf
algebra structure,
the relevant cap pairing  is still available.

\demo{Proof}
At first, we view
$L$ as
a Lie algebra over $R$
and write
$K_RL$ for its Koszul complex
in the category
of $R$-modules.
Consider the
cap pairing 
$$
\cap\colon 
(N\otimes_{UL} K_RL) 
\otimes_R
\roman {Alt}_R(L,M)
@>>>
(N\otimes_RM)\otimes_{UL} K_RL
\tag3.3.3
$$
of chain complexes in the category of $R$-modules.
As observed above, this pairing is induced by the shuffle diagonal map
(3.1.1)
on the exterior $R$-algebra $\Lambda_R^{} L$
and hence may  be written
$$
\cap\colon 
(N\otimes_{R} \Lambda_R^{} L)
\otimes_R
\roman {Hom}_R(\Lambda_R^{} L,M)  
@>>>
(N\otimes_RM)\otimes_{R} \Lambda_R^{} L.
\tag3.3.4
$$
Likewise,
taking now into account the $A$-module
structure on $L$ and
using
the 
shuffle diagonal map (3.1.4) on
the exterior $A$-algebra $\Lambda_A^{} L$,
we obtain a similar pairing
$$
\cap\colon
(N\otimes_A \Lambda_A^{} L) 
\otimes_R
 \roman {Hom}_A(\Lambda_A^{} L,M) 
@>>>
(N\otimes_AM)\otimes_A \Lambda_A^{} L
\tag3.3.5
$$
of graded $R$-modules.
Replacing
the tensor product \lq\lq$\otimes_R$\rq\rq\ 
with the tensor product
\lq\lq$\otimes_A$\rq\rq\ 
over $A$ we would even obtain a pairing
of graded $A$-modules
but below we shall 
stick to the pairing (3.3.5) of
graded $R$-modules.
(We remarked earlier that the distinction between
graded $A$-objects and differential graded $R$-modules
due to the in general non-trivial action of $L$ on $A$
persists throughout.)
In view of the isomorphisms
$$
\roman {Hom}_A(\Lambda_A^{} L,M) 
\cong
\roman {Alt}_A(L,M)
\cong
\roman{Hom}_{U(A,L)}(K(A,L),M)
\tag3.3.6
$$
and
$$
\aligned
N\otimes_A \Lambda_A^{} L
&\cong
N\otimes_{U(A,L)}K(A,L)
\\
(N\otimes_AM)\otimes_A \Lambda_A^{} L
&\cong
(N\otimes_AM)\otimes_{U(A,L)}K(A,L)
\endaligned
\tag3.3.7
$$
of graded $A$-modules,
 differentials being ignored,
the pairing (3.3.5) has the form
$$
\aligned
\cap\colon 
(N\otimes_{U(A,L)}K(A,L))
&\otimes_R
\roman{Hom}_{U(A,L)}(K(A,L),M)
\\
@>>>
&(N\otimes_AM)\otimes_{U(A,L)}K(A,L).
\\
\endaligned
\tag3.3.8
$$
We now assert that this pairing is in fact
compatible with the
differentials.
In order to see this we recall that
(i) the inclusion
$$
\roman {Alt}_A(L,M)
\subseteq
\roman {Alt}_R(L,M)
$$
is a morphism of chain complexes, and that (ii) the projection
maps
$$
N\otimes_R \Lambda_R^{} L 
@>>>
N\otimes_A \Lambda_A^{} L
$$
and
$$
(N\otimes_RM)\otimes_R \Lambda_R^{} L
@>>>
(N\otimes_AM)\otimes_A \Lambda_A^{} L
$$
are morphisms of chain complexes, too.
These two observations
are due to Rinehart \cite\rinehart.
Hence restricting to
$\roman {Alt}_A(L,M)=\roman {Hom}_A(\Lambda_A^{} L,M)$, from
(3.3.4), we obtain the pairing
$$
(N\otimes_{R} \Lambda_R^{} L) 
\otimes_R
\roman {Hom}_A(\Lambda_A^{} L,M) 
@>>>
(N\otimes_RM)\otimes_{R} \Lambda_R^{} 
\tag3.3.9
$$
of chain complexes.
The latter, in turn, fits into the commutative
diagram
$$
\CD
\phantom{AB}(N\otimes_R \Lambda_R^{} L) 
\otimes_R
\roman {Hom}_A(\Lambda_A^{} L,M) 
@>>>
(N\otimes_RM)\otimes_R \Lambda_R^{} L\phantom{AAA}
\\
@VVV
@VVV
\\
\phantom{AB}(N\otimes_A \Lambda_A^{} L) 
\otimes_R
\roman {Hom}_A(\Lambda_A^{} L,M) 
@>>>
(N\otimes_AM)\otimes_A \Lambda_A^{} L\phantom{AAA}
\endCD
\tag3.3.10
$$
having the upper row and the two vertical maps morphisms
of chain complexes.
Since the vertical maps are surjective,
this implies
that the lower row is also a morphism of chain complexes.
The latter induces the cap pairing
$$
\cap\colon 
\roman H_{\ell}(L,N)
 \otimes_R 
\roman H^k(L,M)
@>>>
\roman H_{\ell-k}(L,N\otimes_A M)
$$
which we are looking for. \qed
\enddemo

Suppose now that $L$ is a duality
$(R,A)$-Lie algebra
of rank $n$.
Let $C=\Lambda_A^n L^*$, the dualizing module for $L$.
The duality isomorphism (2.11.2), for $N = C$ and $k=n$, takes the form
$$
\psi
\colon
\roman H_{n}(L,C)
@>>>
\roman {Hom}_U(C,C).
\tag3.4
$$
Let 
$e \in \roman H_n(L,C)$
be the class which under this isomorphism goes to $\roman{Id}_C$.
We refer to it as the {\it fundamental class\/}
of the Lie-Rinehart algebra $(A,L)$.
We note that
$\roman {Hom}_U(C,C)
\cong \roman H^0(L,A)$.
When $A$ and $L$ are the algebra of smooth functions
and Lie algebra of smooth vector fields on a smooth 
(connected)
$n$-dimensional manifold
$W$,
$\roman H^0(L,A) \cong \roman H^0(W,\Bobb R)\cong \Bobb R$
(the constant functions)
and, for $W$ compact, Poincar\'e duality identifies
$\roman H^0(W,\Bobb R)$
with $\roman H_n(W,\Bobb R_t)$ 
where
$\Bobb R_t$ refers to the real valued local system
on $W$ 
arising from the orientation bundle
which is non-trivial if and only if $W$ is not orientable;
our fundamental class
$e$ can then be identified with the ordinary fundamental class
of $W$.
See Remark 4.11 below for details.
\smallskip
For a general Lie-Rinehart algebra $(A,L)$,
with $N = C$, 
the cap pairing (3.3.1) yields
$$
\cap\colon 
\roman H_n(L,C)
\otimes_R
\roman H^k(L,M)  
@>>>
\roman H_{n-k}(L,C\otimes_A M)
\tag3.5
$$
while with $M = \roman {Hom}_A(C,N)$ and the
ordinary evaluation pairing
\linebreak
$\roman {Hom}_A(C,N) \otimes_A C \to N$, we get
$$
\cap\colon 
\roman H_n(L,C) 
 \otimes_R
\roman H^k(L,\roman {Hom}_A(C,N))
@>>>
\roman H_{n-k}(L,N) .
\tag3.6
$$

\proclaim{Theorem 3.7}
Given a duality $(R,A)$-Lie algebra $L$ of rank $n$,
cap-product with its fundamental class $e$ produces isomorphisms
$$
(e \cap\,\cdot\,) \colon\roman H^k(L,M)
@>>>
\roman H_{n-k}(L,C \otimes_A M)
\tag3.7.1
$$
for all non-negative integers $k$ and all
left $(A,L)$-modules $M$ and, furthermore,
isomorphisms
$$
(e \cap\,\cdot\,)  \colon 
\roman H^{n-k}(L,\roman{Hom}_A(C,N))
@>>>
\roman H_{k}(L,N)
\tag3.7.2
$$
for all non-negative integers $k$ and all
right $(A,L)$-modules $N$.
\endproclaim

Consequently the duality isomorphisms may be taken to be natural
in any reasonable sense.
The proof to be given below owes much to the corresponding proof
of Theorem 9.5 in \cite\bieriboo\ but the
proof 
in \cite\bieriboo\ 
cannot  just be adapted to 
Lie-Rinehart algebras;
the basic difference is that, for a Lie-Rinehart
algebra $(A,L)$, the algebra $U(A,L)$
acts non-trivially on $A$ while in \cite\bieriboo\ 
(as always in group cohomology)
the action of the group ring on the ground ring is trivial.

\demo{Proof}
It proceeds in five steps.
\smallskip\noindent
{\it Step 1.}
Here we prove that,
in the top dimension $n$, the duality isomorphism
$$
\psi
\colon
\roman H_n(L,N)
@>>>
\roman {Hom}_U(C,N)                                                            
=\roman H^0(L,\roman {Hom}_A(C,N))
\tag3.7.3
$$
is given by the formula
$$
(\psi(u))(c) = u \cap c,\quad
u \in \roman H_n(L,N),
\  c \in C= \roman H^n(L,U).
\tag3.7.4
$$
Indeed, for the highest non-zero $A$-exterior power
$\Lambda^n_A L$ of $L$, the pairing (3.3.5) 
comes down to the evaluation pairing
$$
(N\otimes_A \Lambda^n_A L) 
\otimes_R
\roman {Hom}_A(\Lambda^n_A L,M) 
@>>>
N\otimes_AM
\tag3.7.5
$$
which, on homology, induces
the pairing
$$
\cap\colon 
\roman H_n(L,N) 
\otimes_R
\roman H^n(L,M) 
@>>>
\roman H_0(L,N\otimes_A M)
=
(N\otimes_A M)\otimes_U A.
$$
In particular,
with $M=U$, $U$ being viewed as a left $U$-module,
we obtain
the pairing
$$
\cap\colon 
\roman H_n(L,N) 
\otimes_R
\roman H^n(L,U) 
@>>>
(N\otimes_A U)\otimes_U A \cong N
$$
which, with reference to the
right
$U$-module structure on $\roman H^n(L,U)$
coming from the right $U$-module structure on itself
(which does not come into play in the construction of
$\roman H^n(L,U)$)
and that  on $N$, is compatible with the
right $U$-module structures.
The adjoint of this pairing
is the morphism of $R$-modules
$$
\roman H_n(L,N) 
@>>>
\roman{Hom}_U(\roman H^n(L,U), N),
\quad
u \mapsto (c \mapsto  u \cap c).
$$
The latter, in turn,
is plainly induced by the adjoint
$$
N\otimes_A \Lambda^n_A L
@>>>
\roman{Hom}_U
(\roman {Hom}_A(\Lambda^n_A L,U),N) 
\tag3.7.6
$$
of (3.7.5)
with $M=U$.
However, 
when the resolution $K$ in (1.4.2) is taken to be the
Rinehart complex (2.9.2),
which under these circumstances is 
a projective resolution of $A$ in the category of left $U(A,L)$-modules,
(3.7.6)  is just a rewrite
of the isomorphism (1.4.2)
in the top dimension.
Since (1.4.2) induces
(3.7.3),
where 
$C=\roman H^n(L,U)$, 
it follows that
(3.7.3) is given by (3.7.4) as asserted.
In particular, 
since $\psi(e) = \roman{Id}_C$,
taking $N=C$, we obtain
$$
e \cap  c = c,\quad
 c \in C= \roman H^n(L,U).
\tag3.7.7
$$
\smallskip\noindent
{\it Step 2.}
Now we claim that,
still in the top dimension,
the duality isomorphism
$$
\phi
\colon
C \otimes _UM
@>>>
\roman H^n(L,M)
\tag3.7.8
$$
is given by the following formula:
Let $\zeta \colon K_n \to U$
represent 
$c \in C= \roman H^n(L,U)$,
let $x \in M$, let
$\omega_x \colon U \to M$ be given by
$\omega_x(1) = x$,
and write
$(\omega_x)_* \colon \roman H^n(L,U) \to  \roman H^n(L,M)$
for the induced map; then 
$$
\phi(c \otimes x) = (\omega_x)_* c.
\tag3.7.9
$$
In fact, the duality isomorphism (3.7.8) is induced by the
canonical isomorphism
$\Phi$ from $K^* \otimes _U M$ onto 
$\roman {Hom}_U(K,M)$ where $K$ still refers to the Rinehart complex
(2.9.2),
cf. (1.4.1).
However, under this isomorphism, $\zeta \otimes x$ goes to the
composite of $\zeta$ with $\omega_x$ whence (3.7.9) holds.
\smallskip\noindent
{\it Step 3.}
Now we claim that,
still in the top dimension,
the special case
$$
(e \cap\,\cdot\,) \colon\roman H^n(L,M)
@>>>
\roman H_0(L,C \otimes_A M)
=
C \otimes_U M
\tag3.7.10
$$
of
(3.7.1)
is 
the inverse of (3.7.8) and hence, in particular,
is an isomorphism.
In order to see this, 
as before, let $x \in M$ and let
$\omega_x \colon U \to M$ be given by
$\omega_x(1) = x$.
Consider the  diagram
$$
\CD
C\otimes _U U
@>{\roman{Id}}>>
\roman H^n(L,U)
@>{(e \cap\,\cdot\,)}>>
C\otimes _U U
\\
@VV{(\omega_x)_*}V
@VV{(\omega_x)_*}V
@VV{(\omega_x)_*}V
\\
C \otimes _U M
@>{\phi}>>
\roman H^n(L,M)
@>{(e \cap\,\cdot\,)}>>
C\otimes _U M
\endCD
\tag3.7.11
$$
By virtue of what has been proved in Step 2,
this diagram is commutative and, in view of (3.7.7),
the composite of the top row is the identity map.
Since $x \in M$ was arbitrary,
this shows that the composite of the bottom row is the identity
of
$C\otimes _U M$,
hence
$(e \cap\,\cdot\,) = \phi^{-1}$
is an isomorphism.
\smallskip\noindent
{\it Step 4.}
Here we show that
(3.7.1) is an isomorphism for every 
non-negative integer
$k$. In order to see this, let
$0 \to M_1 \to P \to M \to 0$
be a short exact sequence of left $U$-modules
with $P$ projective.
Since $C$ is projective as an $A$-module,
the resulting sequence
$$
0 @>>> C \otimes _AM_1 @>>> C \otimes _AP @>>> C \otimes _AM @>>> 0
$$
of right $U$-modules,
with right 
$U$-module structures given by (2.4),
is still exact, and hence naturality of the cap product
yields the commutative diagram
$$
\CD
\roman H^{n-1}(L,P)
@>>>
\roman H^{n-1}(L,M)
@>>>
\roman H^{n}(L,M_1)
@>>>
\roman H^{n}(L,P)
\\
@VVV
@VVV
@VV{\cong}V
@VV{\cong}V
\\
0
@>>>
\roman H_1(L,C \otimes_A M)
@>>>
\roman H_0(L,C \otimes_A M_1)
@>>>
\roman H_0(L,C \otimes_A P),
\endCD
$$
each vertical morphism being given by
$(e \cap\,\cdot\,)$.
However, from (2.11.1),
with $M=P$, we deduce
$$
\roman H^{n-1}(L,P)
\cong
\roman H_1(L, C\otimes_A P)
\cong
\roman {Tor}_U(C,P) = 0
$$
since $P$ is a projective
$U$-module.
Hence
$$
(e \cap\,\cdot\,) \colon\roman H^{n-1}(L,M)
@>>>
\roman H_1(L,C \otimes_A M)
$$
is an isomorphism, too,
and the argument can be iterated. 
This establishes the isomorphisms (3.7.1).                                     
\smallskip\noindent
{\it Step 5.}
The isomorphisms (3.7.2) are now a formal  consequence 
of those already established.
Indeed, 
taking $M = \roman{Hom}_A(C,N)$, 
we get
$$
(e \cap\,\cdot\,) \colon\roman H^k(L,\roman{Hom}_A(C,N))
@>>>
\roman H_{n-k}(L,C \otimes_A \roman{Hom}_A(C,N)).
$$
However,
$C \otimes_A \roman{Hom}_A(C,N) \cong N$ as right
$(A,L)$-modules.
This establishes the isomorphisms (3.7.2) as well. \qed
\enddemo

\smallskip
Let $L$ be a duality $(R,A)$-Lie algebra of rank $n$.
Its duality properties
give rise to certain bilinear pairings in its cohomology,
in the following way:
Let
$M_1 \otimes_A M_2 \to M$
be a pairing of left $(A,L)$-modules.
The naturality of the cap pairings 
and the compatibility properties between
the cup- and cap pairings
imply that
the requisite duality isomorphisms combine to a commutative diagram
$$
\CD
\phantom{{A}^{n-} }\roman H^k(L,M_1) \otimes_R \roman H^{n-k}(L,M_2)\,\,
@>>>
\roman H^{n}(L,M)
\\
@VV{(e \cap\,\cdot\,)\otimes_R \roman{Id}}V
@VV{(e \cap\,\cdot\,)}V
\\
\roman H_{n-k}(L,C \otimes _AM_1) 
\otimes_R \roman H^{n-k}(L,M_2) \phantom{\otimes M}
@>>>
\,\phantom { {}_U }C \otimes _UM
\endCD
\tag3.8
$$
having 
as upper row the corresponding cup pairing (3.2),
as lower row the cap pairing (3.3.1)
(with the notation $k$, $\ell$, $M$, $N$ adjusted appropriately)
 and, furthermore,
all vertical maps isomorphisms;
here $U = U(A,L)$. From (3.8) we obtain the bilinear pairing
$$
\roman H^k(L,M_1) \otimes_R \roman H^{n-k}(L,M_2)
@>>>
C \otimes _UM .
\tag3.9
$$
The question we will study in the next section is whether, for suitable 
pairings 
$M_1 \otimes_A M_2 \to M$,
nondegenerate in a suitable way,
the pairing (3.9) may be nondegenerate,
in an appropriate sense.
\smallskip
\noindent
{\smc Remark 3.10.}
This question will, in general, not have a naive solution,
as the following example shows:
Under the circumstances of (2.10),
suppose that the ground ring is that of the reals,
$\Bobb R$,
let
$A$ be the algebra of
smooth real functions on a smooth real manifold $W$
and $\fra g$ a real Lie algebra acting infinitesimally on $W$,
and let $L = A \otimes _{\Bobb R} \fra g$,
the $(\Bobb R,A)$-Lie algebra introduced in (2.10).
Further, let $\Cal O$ be the space of sections of the orientation
bundle of $W$,
compactly supported
when $W$ is not compact,
with its canonical
left
$(A,L)$-module structure,
and denote by $\omega _A$
the dualizing module of
$\roman{Der}(A) = \roman{Vect}(W)$;
under the present circumstances,
$\omega _A$ is just the highest non-zero exterior power
of the space of sections of the cotangent bundle of $W$.
Real valued bilinear cohomology pairings 
can then be obtained by integration over $W$, in the following way:
Consider a bilinear pairing
$$
M_1 \otimes _A M_2 @>>> \roman{Hom}_A(C_L, \omega_A \otimes _A \Cal O)
$$
of left
$(A,L)$-modules
arising from a nondegenerate pairing of
smooth vector bundles on $W$,
the target
$\roman{Hom}_A(C_L, \omega_A \otimes _A \Cal O)$
being endowed with the obvious left $(A,L)$-module structure
(2.2).
The resulting pairing (3.9) 
then has its values in the space
$\omega_A \otimes _{U(A,L)}\Cal O$
which amounts to that of (compactly supported) densities on $W$,
with the $\fra g$-action being divided
out. Integration yields a map from this space to the reals
but  this map will not be an isomorphism
unless the structure map from $L$ to $\roman{Vect}(W)$ is surjective
whence 
the resulting real-valued pairing 
may be nondegenerate only if
the space
of densities with the $\fra g$-action being divided out
amounts to a single copy of the reals.
This confirms and explains the empirical observations
made in 
\cite\evluwein.
In particular,
the example given there shows that
$\omega_A \otimes _{U(A,L)}\Cal O$
may well be larger than a single copy of the reals.
The cure is provided by a notion of nondegeneracy
over more general rings:
When the structure map from $L$ to $\roman{Vect}(W)$ is not surjective,
in a sense, the resulting pairing (3.9) 
looks nondegenerate over
$\omega_A \otimes _{U(A,L)}\Cal O$
rather than
over $\Bobb R$ except that
this does not make sense as it stands since
$\omega_A \otimes _{U(A,L)}\Cal O$
is not a ring.
By means of 
an additional piece of structure
which may or may not exist,
that is to say, by means of the notion of a {\it trace\/},
we shall make precise 
this idea of nondegeneracy over more general rings
in the next section.

\medskip\noindent{\bf 4. Poincar\'e duality}\smallskip\noindent
In this section we shall show that,
for certain Lie-Rinehart algebras,
Poincar\'e duality holds in their cohomology,
the terminology \lq\lq Poincar\'e\rq\rq\ 
referring to the nondegeneracy of certain pairings
of the kind
(3.9) in a sense which we are about to make precise.
As a special case,
we obtain a new proof of Poincar\'e duality
in the de Rham cohomology of a smooth manifold.
We shall use the descriptions of
the duality isomorphisms in terms of the cap product with the fundamental 
class given in (3.7). This entails all the requisite naturality properties.
\smallskip
As before, let 
$A$ be a commutative algebra
and $L$  a duality $(R,A)$-Lie algebra
of rank $n$.
A left $(A,L)$-module $M$ which, as an $A$-module, 
is finitely generated and projective,
will henceforth be referred to as a
left {\it duality\/} $(A,L)$-module;
likewise a right $(A,L)$-module $N$ which, as an $A$-module, 
is finitely generated and projective,
will be called  a
right {\it duality\/} $(A,L)$-module.
We shall say that a {\it pretrace\/}
for $L$ consists of a left
$(A,L)$-module
$\Cal O$,
referred to as its {\it trace module\/},
together 
with an isomorphism
$$
t
\colon
\roman H^n(L,\Cal O)
@>>>
R
\tag4.1
$$
of $R$-modules;
a pretrace
$(\Cal O, t)$
for $L$
will be said to be a {\it trace\/}
for $L$ provided
a {\it weak\/} form of {\it Poincar\'e duality\/}
holds in the sense that,
for every
left duality $(A,L)$-module $M$,
$\roman H^0(L,M)=\roman{Hom}_U(A,M)$
and
$\roman H^n(L,\roman{Hom}_A(M,\Cal O))$
are dually paired, that is to say,
the canonical morphism
$$
\aligned
\roman{Hom}_U(A,M)
&@>>>
\roman {Hom}_R(\roman H^n(L,\roman{Hom}_A(M,\Cal O)),\roman H^n(L,\Cal O))
\\
&@>{t_*}>>
\roman {Hom}_R(\roman H^n(L,\roman{Hom}_A(M,\Cal O)),R)
\endaligned
\tag4.2.1
$$
of $R$-modules
is an isomorphism,
the morphism $t_*$ being the isomorphism induced by the trace.
\smallskip
A variant of this  notion of trace is obtained when the canonical morphism
(4.2.1) is only required to be an isomorphism
for every left duality $(A,L)$-module $M$ belonging to a suitable
class of
left duality $(A,L)$-modules. An example will be given in (4.15) below.
\smallskip
Our  notion of \lq trace\rq\  is similar to
a corresponding one
in the theory of Serre duality for the cohomology
of projective sheaves on a projective scheme,
cf. e.~g. \cite\hartshor,
but this terminology does not necessarily 
coincide with other usages thereof in de Rham and Poisson
cohomology.
\smallskip
Let $C_L$ denote the dualizing module of $L$,
cf. (2.10).
In view of the naturality of the duality isomorphism
$e \cap \cdot\,\, \colon \roman H^n(L,\Cal O) \to \roman C_L \otimes_U\Cal O$,
cf. the commutative diagram (3.8),
a trace may as well be described
as an isomorphism
$$
t
\colon
C_L \otimes_U\Cal O
@>>>
R
\tag4.3
$$
of $R$-modules such that, for every
left duality $(A,L)$-module $M$,
the canonical morphism
$$
\roman{Hom}_U(A,M)
@>>>
\roman {Hom}_R(C_L \otimes_U \roman{Hom}_A(M,\Cal O),C_L \otimes_U\Cal O)
@>{t_*}>>
\roman {Hom}_R(C_L \otimes_U M,R)
\tag4.4.1
$$
is an isomorphism of $R$-modules.
\smallskip\noindent
{\smc Example 4.5.}
Let 
$A$ be a commutative algebra and 
$\fra g$ an ordinary Lie algebra over $A$
which we suppose finitely generated and projective as an $A$-module,
of constant rank $k$.
To subsume
this example under our general
theory,
we take $A$ as the ground ring,
written $R$,
and distinguish 
$A$ and $R$ deliberately in notation
although
$A$ and $R$ coincide,
but this notational distinction will
enable us to abstract from the special situation
of this example.
We then view $(A, \fra g)$ as a duality
$(R,A)$-Lie algebra of rank $k$.
A typical case
arises from  a principal $G$-bundle
$\xi \colon P \to B$
with structure group $G$: its Lie
algebra $\fra g(\xi)$ of infinitesimal gauge transformations
is the space of sections of the adjoint bundle
and inherits a Lie algebra structure 
over the ring of smooth functions on $B$.
Returning to a general Lie algebra
$\fra g$ over $A$ of constant rank $k$,
let $\Cal O_{\fra g} = \Lambda^{\roman k}\fra g$,
the top exterior power of $\fra g$,
with its standard  left
$\fra g$-module structure.
With reference to the canonical right
$\fra g$-module structure on
$\roman{Hom}_A(U\fra g, A)$
induced by the 
left $\fra g$-structure on $U\fra g$,
we have 
the homology group
$$
\roman H_k(\fra g,\roman{Hom}_A(U\fra g, A) );
\tag4.5.1
$$
from
the right $\fra g$-structure on $U\fra g$
which remains free
under the construction
of $\roman H_k(\fra g,\roman{Hom}_A(U\fra g, A) )$,
this homology group  inherits
a  left 
$\fra g$-module structure  and,
with this structure,
$\roman H_k(\fra g,\roman{Hom}_A(U\fra g, A) )$
 is canonically isomorphic to
$\Cal O_{\fra g}$
as a left
$\fra g$-module; in fact,
$\roman H_k(\fra g,\roman{Hom}_A(U\fra g, A) )$
and
$\roman H^0(\fra g,\roman{Hom}_A(C_{\fra g},\roman{Hom}_A(U\fra g, A) ))$
are canonically isomorphic by duality,
and
$$
\roman H^0(\fra g,\roman{Hom}_A(C_{\fra g},\roman{Hom}_A(U\fra g, A) ))
\cong
\roman{Hom}_A(C_{\fra g}, A)
\cong \Cal O_{\fra g}.
$$
As a right
$(A,L)$-module, the dualizing module $C_{\fra g}$
is therefore naturally isomorphic to
$\roman{Hom}_A(\Cal O_{\fra g},A)=\roman{Hom}_A(\Lambda^k\fra g,A)$,
the 
right
$\fra g$-module structure on the latter being induced by the 
left
$\fra g$-module structure on 
$\Lambda^k\fra g$.
Then
$$
C_{\fra g}\otimes_A\Cal O_{\fra g}
\cong \roman{Hom}_A(C_{\fra g},C_{\fra g}) \cong A
\tag4.5.2
$$
inherits a  right
$\fra g$-module structure which has to be trivial,
since manifestly
$$
\roman{Hom}_{\fra g}(C_{\fra g},C_{\fra g}) \cong R =A.
$$
Consequently,
in the top dimension,
the duality isomorphism for 
the left $\fra g$-module
$\Cal O_{\fra g}$ takes the form
$$
\roman H^k(\fra g,\Cal O_{\fra g}) @>{e \cap\,\cdot\,\,}>> \roman H_0(\fra g, A)
\cong A
$$
and this in fact defines a trace
$$
t\colon
\roman H^k(\fra g,\Cal O_{\fra g})
@>>> A.
\tag4.5.3
$$
It is readily seen
that the property (4.2.1) holds;
we do not give the details.
Thus an ordinary Lie algebra $\fra g$ 
of finite constant rank
always comes with a trace,
having as trace module
the highest non-zero homology group
$\roman H_{\roman{top}}(\fra g,\roman{Hom}_A(U\fra g, A) )$,
and this homology group is canonically isomorphic to
the highest non-zero exterior power
$\Lambda^{\roman{top}}\fra g$
of $\fra g$;
here $A$ is to be viewed as a right $\fra g$-module
with trivial $\fra g$-action.
Moreover, given a left $\fra g$-module $M$,
the
canonical (or cap) pairing
$$
\roman H^{\ell}(\fra g,M) 
\otimes_A\roman H_{\ell}(\fra g,\roman{Hom}_A(M, A))
@>>>
\roman H_0(\fra g,A)
@>\cong>> A
\tag4.5.4
$$
is manifestly nondegenerate
because it arises from a nondegenerate pairing of
$A$-chain complexes.
In view of the obvious isomorphisms
$$
C_{\fra g} \otimes _A\roman{Hom}_A(M, \Cal O_{\fra g})
\cong
\roman{Hom}_A(M, C_{\fra g} \otimes _A\Cal O_{\fra g})
\cong
\roman{Hom}_A(M, A),
$$
the nondegeneracy of (4.5.4) implies at once that
of the cup pairing
$$
\roman H^{\ell}(\fra g,M) 
\otimes_A\roman H^{k-\ell}(\fra g,\roman{Hom}_A(M, \Cal O_{\fra g}))
@>>>
\roman H^k(\fra g,\Cal O_{\fra g})
@>t>> A .
\tag4.5.5
$$
In particular, the canonical duality pairing
between 
$C_{\fra g}$ and
$\Cal O_{\fra g}$
may be described as the canonical evaluation pairing
$$
\roman H^k(\fra g,U\fra g)
\otimes_A
\roman H_k(\fra g,\roman{Hom}_A(U\fra g, A) )
@>>>
\roman H_0(\fra g,A) 
\cong A.
\tag4.5.6
$$
The pairing (4.5.5) is the Poincar\'e duality pairing
in ordinary Lie algebra cohomology,
phrased here over an arbitrary commutative algebra $A$.
The corresponding Lie-Rinehart structure is trivial
and the duality properties are well known;
we therefore refer to this pairing as the {\it classical\/} 
duality pairing.
For example,
given a principal  bundle $\xi$ over a smooth manifold $B$,
the pairing (4.5.5) is available 
over the ring $A$ of smooth functions on $B$
for
the Lie algebra $\fra g(\xi)$ of infinitesimal gauge transformations.
\smallskip
For a general duality
Lie-Rinehart algebra $(A,L)$,
the highest non-zero exterior power
of $L$ over $A$
with its naive left $L$-module structure
does {\it not\/}
yield a left $(A,L)$-module and hence
cannot naively serve for the construction of a trace module,
and a trace is an additional piece of structure.
We now give a different description of a trace
involving the appropriate direct generalization of the trivial
right $\fra g$-module structure
on the copy of $A$ occurring in
the top homology of $\fra g$
with coefficients in
$\roman{Hom}_A(U\fra g, A)$ just explained,
cf. (4.5.1).
This description will be used in the next section to construct
traces for extensions of Lie-Rinehart algebras.

\proclaim{Proposition 4.6}
Let $(A,L)$ be a duality
Lie-Rinehart algebra.
\newline\noindent
{\rm (i)} Given a trace $(\Cal O, t)$ for
$(A,L)$,
the 
projective rank one  $A$-module
$V= C_L \otimes _A\Cal O$
with
right 
$(A,L)$-module structure {\rm (2.4)}
has the property that
$t$ induces an
isomorphism $\iota \colon V \otimes_UA \to R$
and that,
for every right duality $(A,L)$-module $N$,
the canonical map
$$
\roman{Hom}_U(N,V)
@>>>
\roman{Hom}_R(N \otimes _UA,
V \otimes _UA)
@>{\iota_*}>>
\roman{Hom}_R(\roman H_0(L,N),R)
\tag4.6.1
$$
is an isomorphism.
\newline\noindent
{\rm (ii)} 
Conversely,
given a 
right 
$(A,L)$-module 
$V$ which, as an $A$-module,
is projective of rank one,
together with an isomorphism
$\iota \colon V \otimes_UA \to R$,
let
$\Cal O = \roman{Hom}_A(C_L,V)$,
with its induced left $(A,L)$-module structure
{\rm (2.3)}  
so that $C_L \otimes _A\Cal O$ is canonically isomorphic to $V$;
then the composite
$$
\roman H^n(L,\Cal O)
@>>>
\roman H_0(L,C_L \otimes _A\Cal O)
\cong 
\roman H_0(L,V)
\cong
V \otimes _U A
@>{\iota}>> R
\tag4.6.2
$$
yields a trace provided that,
for every right duality $(A,L)$-module $N$,
the canonical map
$$
\roman{Hom}_U(N,V)
@>>>
\roman{Hom}_R(N \otimes _UA,
V \otimes _UA)
@>{\iota_*}>>
\roman{Hom}_R(\roman H_0(L,N),R)
\tag4.6.3
$$
is an isomorphism.
\endproclaim

\demo{Proof}
The argument for (i) is straightforward.
We briefly indicate one for (ii):
Given a left duality $(A,L)$-module $M$,
let $N =\roman{Hom}_A(M,V)$,
with
right 
$(A,L)$-module structure {\rm (2.4)}.
Then
$$
\roman{Hom}_U(A,M)
\cong
\roman{Hom}_U(\roman{Hom}_A(M,V),V)
\cong
\roman{Hom}_U(N,V)
$$
and
$\roman{Hom}_A(M, \Cal O) \cong
\roman{Hom}_A(C_L,N)$ whence
$$
\roman{Hom}_R(\roman H^n(L, \roman{Hom}_A(M, \Cal O)),R)
\cong
\roman{Hom}_R(\roman H^n(L, \roman{Hom}_A(C_L,N)),R)
$$
and, by duality,
$$
\roman{Hom}_R(\roman H^n(L, \roman{Hom}_A(C_L,N)),R)
\cong
\roman{Hom}_R(\roman H_0(L,N),R).
$$
Since (4.6.3) is an isomorphism, 
(4.2.1) is an isomorphism, too. \qed
\enddemo
The right $(A,L)$-module $V$ 
in (4.6) arises by abstraction from the isomorphism (4.5.2). 
In all the examples that we know,
as an $A$-module,
this module is actually free of rank 1, that is, just a copy
of $A$, so that the two modules $C_L$ and $\Cal O$ are
dual to each other
as $A$-modules.

\smallskip\noindent
{\smc Example 4.7.}
Let 
the ground ring be that of the reals, $\Bobb R$, let
$W$ be
a  smooth real $n$-dimensional manifold $W$,
$A$ its algebra of 
smooth functions,
$L$ the $(\Bobb R,A)$-Lie algebra $\roman{Vect}(W)$ of smooth vector fields
on $W$,
and let $\Cal O$ be the space of compactly supported sections of the
orientation bundle of $W$,
with its canonical 
flat connection and hence
left $(A,L)$-module structure.
The dualizing module $C_L$ of $L$ is just the space
$\omega_A$ of smooth $n$-forms on $W$, with  right
$(A,L)$-module structure given by (2.9.1), for $M=A$.
Then the tensor product
$\omega_A \otimes_A \Cal O$ is the space of compactly supported
densities on $W$,
and the integration map 
from $\omega_A \otimes_A \Cal O$
to $\Bobb R$
induces a trace 
$$
t \colon
\omega_A \otimes_U \Cal O
@>>> \Bobb R
\tag4.7.1
$$
for $L$. In fact, a left duality $(A,L)$-module $M$
arises as the space of sections of a flat vector bundle
$\zeta_M$ on $W$,
and (4.2.1) comes down to the fact that
$\roman H^0(W,\zeta_M)$
and
$\roman H_{\roman{cs}}^n(W,\zeta_M^*)$
are dually paired
where
$\roman H_{\roman{cs}}$ refers to cohomology with compact supports
and $\zeta_M^*$ to the dual bundle.
Notice that when $W$ is compact
$\Cal O$ is projective as an $A$-module
and there is no need to talk about compactly
supported forms;
furthermore,
in this case,
$\omega_A \otimes_A \Cal O$
is a free $A$-module of dimension one ---
that of densities on $W$ ---
and thus, endowed with the 
right $(A,L)$-module structure (2.4),
$A$ carries a  right
$(A,L)$-module structure; when $A$
is viewed as a right
$(A,L)$-module in this way,
we denote it henceforth by $A_r$.
Integration 
then induces
an isomorphism
$\iota\colon A_r \otimes _UA \to \Bobb R$
in such a way that
the right $(A,L)$-module
$V$ in (4.6) is just $A_r$ and that 
$\iota$ satisfies the conditions of (4.6)(ii).
Thus, with
$t$ being defined by (4.6.2),
$(\Cal O,t)$ is a trace
for $(A,L)$.
\smallskip
We now return to a general duality $(R,A)$-Lie algebra $L$,
for an arbitrary $R$-algebra $A$.
Given a trace
$(\Cal O, t)$ for $L$
having $\Cal O$ projective as an $A$-module,
in (4.2.1) we can replace 
$M$ with $\roman{Hom}_A(M,\Cal O)$; thus instead of (4.2.1)
we can then equivalently 
define a trace by
the more appealing
requirement 
that the canonical morphism
$$
\roman{Hom}_U(M,\Cal O)
@>>>
\roman {Hom}_R(\roman H^n(L,M),\roman H^n(L,\Cal O))
@>{t_*}>>
\roman {Hom}_R(\roman H^n(L,M),R)
\tag4.2.2
$$
be an isomorphism of $R$-modules
for {\it every\/}
left duality $(A,L)$-module $M$;
likewise, (4.4.1) 
to be an isomorphism is then equivalent
to the canonical morphism
$$
\roman{Hom}_U(M,\Cal O)
@>>>
\roman {Hom}_R(C_L \otimes_U M,C_L \otimes_U\Cal O)
@>{t_*}>>
\roman {Hom}_R(C_L \otimes_U M,R)
\tag4.4.2
$$
of $R$-modules
being an isomorphism.
Technically
the definition given at the beginning of the present section
allows for more freedom,
though,
for example when Lie-Rinehart algebras
defined over non-compact smooth manifolds
come into play; see e.~g. (4.10) below.
The usual abstract nonsense establishes  the following
the proof of which is left to the reader.

\proclaim{Lemma 4.8}
A trace for $L$
having trace module $\Cal O$
projective as an $A$-module is unique up to isomorphism.
More precisely,
given two such traces
$(\Cal O_1, t_1)$ and
$(\Cal O_2, t_2)$,
there is a unique isomorphism
of left $(A,L)$-modules between 
$\Cal O_1$ and $\Cal O_2$
identifying the two traces. \qed
\endproclaim

Let $L$ be a duality $(R,A)$-Lie algebra,
and suppose  that $L$ is endowed with a trace
$(\Cal O, t)$.
Given
a left
$(A,L)$-module $M$,
in view of what has been said in Section 3 before,
the canonical pairing
$M\otimes_A
\roman{Hom}_A(M, \Cal O)
@>>>
\Cal O
$
of $A$-modules 
induces 
the cap pairing
$$
\cap
\colon
\roman{Tor}_k^U(C_L, M)
\otimes_R
\roman{Ext}^k_U(M,\Cal O)
@>>>
C_L \otimes _U \Cal O @>\iota>> R
\tag4.9.1
$$
and the cup pairing
$$
\cup\colon
\roman H^k(L,M)
\otimes_R
\roman H^{n-k}(L,\roman{Hom}_A(M,\Cal O))
@>>>
\roman H^n(L,\Cal O) @>t>> R.
\tag4.9.2
$$
The cap pairing (4.9.1) may manifestly  be written
$$
\cap
\colon
\roman H_k(L,C_L\otimes _A M)
\otimes_R
\roman H^k(L,\roman{Hom}_A(M,\Cal O))
@>>>
C_L \otimes _U \Cal O @>\iota>> R
\tag4.9.3
$$
as well.
\smallskip
Likewise,
given
a right
$(A,L)$-module $N$,
taking $M=
\roman{Hom}_A(N, V)$
with left $(A,L)$-module structure (2.3),
where $V = C_L \otimes_A \Cal O$, cf. (4.6),
with right $(A,L)$-module structure (2.4),
we obtain the corresponding cap pairing
(3.3.1)
which now has the form
$$
\cap
\colon
\roman H_k(L,N)
\otimes_R
\roman H^k(L,\roman{Hom}_A(N,V))
@>>>
\roman H_0(L,N \otimes _A\roman{Hom}_A(N,V));
$$
combining it with the morphism on $\roman H_0$
induced by
the canonical morphism
$$
N\otimes_A
\roman{Hom}_A(N, V)
@>>>
V
$$
of right $(A,L)$-modules,
we obtain the pairing
$$
\cap
\colon
\roman H_k(L,N)
\otimes_R
\roman H^k(L,\roman{Hom}_A(N,V))
@>>>
\roman H_0(L ,V) @>\iota>> R 
\tag4.9.4
$$
which we refer to as
{\it cap pairing\/} as well.
Given a left
$(A,L)$-module $M$,
the pairing (4.9.3)
coincides with the pairing (4.9.4)
for the right $(A,L)$-module $N= C_L\otimes _A M$,
with right $(A,L)$-module structure (2.4);
given a right
$(A,L)$-module $N$,
the pairing (4.9.4)
coincides with the pairing (4.9.3)
for the left $(A,L)$-module $M= \roman{Hom}_A(C_L, N)$,
with left $(A,L)$-module structure (2.3).
The naturality of the duality isomorphisms via the cap product
with the fundamental class, cf. (3.7), establishes at once the following.

\proclaim{Lemma 4.9} 
Given a left
$(A,L)$-module $M$,
the following are equivalent.
\newline\noindent
{\rm (i)} The cap pairing {\rm (4.9.1)} is nondegenerate 
in the sense that its adjoint
$$
\cap^{\sharp}
\colon
\roman{Ext}^k_U(M,\Cal O)
@>>>
\roman{Hom}_R(\roman{Tor}_k^U(C_L, M), R)
\tag4.9.5
$$
is an isomorphism of $R$-modules.
\newline\noindent
{\rm (ii)} The cup pairing {\rm (4.9.2)} is nondegenerate 
in the sense that its adjoint
$$
\cup^{\sharp}
\colon
\roman H^{n-k}(L,\roman{Hom}_A(M,\Cal O))
@>>>
\roman{Hom}_R(\roman H^k(L,M),R)
\tag4.9.6
$$
is an isomorphism of $R$-modules. \qed
\endproclaim

We note that the adjoint (4.9.5)
of (4.9.1)
may as well be written
$$
\cap^{\sharp}
\colon
\roman H^k(L,\roman{Hom}_A(M,\Cal O))
@>>>
\roman{Hom}_R(\roman H_k(L,C_L\otimes _A M), R).
\tag4.9.7
$$
Likewise
the adjoint of
the cap pairing (4.9.4) takes the form
$$
\cap^{\sharp}
\colon
\roman H^k(L,\roman{Hom}_A(N,V))
@>>>
\roman{Hom}_R(\roman H_k(L,N),R).
\tag4.9.8
$$
\smallskip
Given a 
duality $(R,A)$-Lie algebra with a trace
$(\Cal O,t)$ and, furthermore, given a
left $(A,L)$-module $M$,
we shall say that $L$
satisfies {\it Poincar\'e duality 
for\/} $M$  if
the adjoint (4.9.5) of the
cap pairing (4.9.1) is an isomorphism or, equivalently,
if the adjoint (4.9.6) of the
cup pairing (4.9.2) is an isomorphism of $R$-modules;
given a right $(A,L)$-module $N$,
we shall say that $L$
satisfies {\it Poincar\'e duality 
for\/} $N$  if
the adjoint (4.9.8) of the
cap pairing (4.9.4) is an isomorphism.
We note that
$L$
satisfies 
Poincar\'e duality 
for a left $(A,L)$-module $M$  if
and only if it satisfies
Poincar\'e duality 
for the right $(A,L)$-module $C_L\otimes _AM$;
likewise $L$
satisfies 
Poincar\'e duality 
for a right $(A,L)$-module $N$  if
and only if it satisfies
Poincar\'e duality 
for the left $(A,L)$-module $\roman{Hom}_A(C_L,N)$.

\smallskip
Let $W$ be a  connected real smooth manifold.
We maintain the notation in (4.7) and do not repeat it.
There are various known ways of establishing
Poincar\'e duality in this case:
One may introduce a Riemannian metric on $W$.
This induces inner product structures on the
constituents of the various de Rham complexes
and, via the appropriate Sobolev completions,
the canonical pairing
$C_L \otimes _{\Bobb R}\Cal O \to \Bobb R$
passes to a perfect pairing of Hilbert spaces.
In the orientable case this is a version of the standard
$L_2$-pairing between functions and densities.
From this observation,
the nondegeneracy of the cap  pairing
(4.9.3)
may be deduced:
The adjoint
of the canonical pairing
$$
\left(\omega_A \otimes _U (K \otimes_A M)\right)\otimes _{\Bobb R}
\roman{Hom}_U(K \otimes_A M,\Cal O)
@>>>
\omega_A \otimes _U\Cal O @>\iota>>
\Bobb R
$$
of chain complexes
has the form
$$
\omega_A \otimes _U (K \otimes_A M)
@>>>
\roman{Hom}_{\Bobb R}(\roman{Hom}_U(K \otimes_A M,\Cal O), \Bobb R);
$$
the corresponding space of de Rham currents $C$
is a subspace of the target
and the canonical map
from
$\omega_A \otimes _U (K \otimes_A M)$ to $C$
has dense image and
is a chain equivalence.
This is similar to the reasoning in de Rham's book
\cite\derhaboo.
We do not spell out the details
since this would not provide any new insight.
\smallskip
Another argument establishing
nondegeneracy proceeds by integration against
suitable dual cell decompositions of $W$.
Again this would not provide anything new and we refrain from spelling out 
details. 
\smallskip
A third way of establishing
nondegeneracy
is by means of a Mayer-Vietoris
argument.
Since
this will allow for later generalization
we now give a proof  
for the ordinary de Rham cohomology of manifolds
along these lines,
within the present framework.

\proclaim{Proposition 4.10}
Under the circumstances of {\rm (4.7)},
suppose that $W$ can be covered by finitely many open contractible sets
$V_1,\dots, V_\ell$
such that each non-empty intersection 
$V_{j_1} \cap \dots \cap V_{j_k}$
is itself contractible.
Then, given a left
duality $(A,L)$-module $M$, 
the adjoint
{\rm (4.9.5)}
of the cap pairing 
{\rm (4.9.1)}
(both pairings over $R=\Bobb R$)
is an isomorphism of real vector spaces.
\endproclaim

For example, every compact manifold
may be covered by finitely many open sets in such a way that the hypothesis
spelled out above is satisfied.
We remind the reader that
here
a left
duality $(A,L)$-module
is just the space of sections of a flat vector bundle on $W$.

\demo{Proof}
On an open contractible subset of $W$,
Poincar\'e duality comes essentially down to the defining property
of a trace.
The idea of the proof is to 
reduce the general case to that of an open contractible
manifold
by means of a Mayer-Vietoris argument.
\smallskip
The details are as follows.
Let $K=K(A,L)$ be the 
Rinehart complex (2.9.2); it is a projective
resolution of $A$ in the category of left
$U(A,L)$-modules,
cf. what is said in Section 2.
Then $K \otimes_A M$, endowed with the 
left $(A,L)$-module structure (2.1) (in each degree),
is a projective
resolution of $M$ in the category of left
$(A,L)$-modules, and the adjoint
(4.9.4) is induced by the canonical homomorphism 
$$
\roman{Hom}_U(K \otimes_A M,\Cal O)
@>>>
\roman{Hom}_{\Bobb R}(\omega_A \otimes _U (K \otimes_A M),
\omega_A \otimes _U\Cal O)
\tag4.10.1
$$
of chain complexes of
real vector spaces.
\smallskip
Let $V \subset W$ be an open subspace.
The two sides of (4.10.1) are then defined over
$W$ and over $V$; we indicate this by the notation
$\roman{Hom}_U(K \otimes_A M,\Cal O)|_W$,
$\roman{Hom}_U(K \otimes_A M,\Cal O)|_V$,
etc.
Restriction induces a canonical 
homomorphism
$$
\omega_A \otimes _U (K \otimes_A M)|_W
@>>>
\omega_A \otimes _U (K \otimes_A M)|_V;
$$
furthermore, since
$\Cal O$
consists of compactly supported sections,
there are canonical injection homomorphisms
$$
\roman{Hom}_U(K \otimes_A M,\Cal O)|_V
@>>>
\roman{Hom}_U(K \otimes_A M,\Cal O)|_W
$$
and
$$
\omega_A \otimes _U\Cal O|_V
@>>>
\omega_A \otimes _U\Cal O|_W,
$$
and these combine to a commutative diagram
$$
\CD
\roman{Hom}_U(K \otimes_A M,\Cal O)|_V
@>>>
\roman{Hom}_{\Bobb R}(\omega_A \otimes _U (K \otimes_A M),
\omega_A \otimes _U\Cal O)|_V
\\
@VVV
@VVV
\\
\roman{Hom}_U(K \otimes_A M,\Cal O)|_W
@>>>
\roman{Hom}_{\Bobb R}(\omega_A \otimes _U (K \otimes_A M),
\omega_A \otimes _U\Cal O)|_W
\endCD
\tag4.10.2
$$
of real chain complexes.
\smallskip
For any open subset $V$ of $W$, we now write
$$
C(V) = \roman{Hom}_U(K \otimes_A M,\Cal O)|_V,\ 
\Omega(V)=\omega_A \otimes _U (K \otimes_A M)|_V,\ 
\Omega^*(V)=\roman{Hom}_{\Bobb R}(\Omega(V),\Bobb R).
$$
Given two open sets $V_1$ and $V_2$, 
we obtain  two exact sequences of chain complexes
$$
0
@>>>
C(V_1 \cap V_2)
@>>>
C(V_1)\oplus C(V_2)
@>>>
C(V_1 \cup V_2)
@>>>
0
\tag4.10.3
$$
and
$$
\Omega(V_1 \cap V_2)
@<<<
\Omega(V_1)\oplus \Omega(V_2)
@<<<
\Omega(V_1 \cup V_2)
@<<<
0.
\tag4.10.4
$$
The sequence (4.10.4) is the ordinary one 
used to derive the Mayer-Vietoris sequence
in de Rham theory; 
even though the morphism
from $\Omega(V_1)\oplus \Omega(V_2)$
to
$\Omega(V_1 \cap V_2)$
in (4.10.4) 
is not surjective,
(4.10.4) 
is well known to induce a Mayer-Vietoris sequence
$$
\cdots
@<<<
\roman H_j\Omega(V_1 \cap V_2)
@<<<
\roman H_j\Omega(V_1) \oplus \roman H_j  \Omega(V_2)
@<<<
\roman H_j\Omega(V_1 \cup V_2)
@<<<
\roman H_{j+1}\Omega(V_1 \cap V_2)
@<<<
\cdots
$$
In fact, the more usual description of this Mayer-Vietoris sequence
looks like
$$
@>>>
\roman H^{n-j-1}(V_1 \cap V_2)
@>>>
\roman H^{n-j}(V_1 \cup V_2)
@>>>
\roman H^{n-j}(V_1) \oplus \roman H^{n-j}(V_2)
@>>>
\roman H^{n-j}(V_1 \cap V_2)
@>>>
\cdots
$$
where the coefficients are not indicated in notation.
The exactness of the Mayer-Vietoris sequence
may e.~g. be derived 
from the corresponding sheaf version 
by an application of the standard device relating 
\v Cech cohomology
with ordinary (singular) cohomology;
cf. e.~g. (II.5.6)  in \cite{\godebook} (p. 219).
It follows that
(4.10.4) 
induces as well a Mayer-Vietoris sequence
$$
\cdots
@>>>
\roman H^j\Omega^*(V_1 \cap V_2)
@>>>
\roman H^j\Omega^*(V_1) \oplus \roman H^j  \Omega^*(V_2)
@>>>
\roman H^j\Omega^*(V_1 \cup V_2)
@>>>
\roman H^{j+1}\Omega^*(V_1 \cap V_2)
@>>>
$$
On the other hand,
it is manifest that 
(4.10.3)
induces a Mayer-Vietoris sequence
$$
\cdots
@>>>
\roman H^j C(V_1 \cap V_2)
@>>>
\roman H^j C(V_1) \oplus \roman H^j C(V_2)
@>>>
\roman H^j C(V_1 \cup V_2)
@>>>
\roman H^{j+1}C(V_1 \cap V_2)
@>>>
\cdots
$$
This is in fact the Mayer-Vietoris sequence
$$
\cdots
@>>>
\roman H_{\roman{cs}}^j (V_1 \cap V_2)
@>>>
\roman H_{\roman{cs}}^j (V_1) \oplus \roman H_{\roman{cs}}^j (V_2)
@>>>
\roman H_{\roman{cs}}^j (V_1 \cup V_2)
@>>>
\roman H_{\roman{cs}}^{j+1}(V_1 \cap V_2)
@>>>
\cdots
$$
in compactly supported 
de Rham cohomology $\roman H_{\roman{cs}}$
where again the coefficients are not indicated in notation.
The morphism
(4.10.1) 
yields
in fact a transformation
of functors 
from $C(\cdot)$ to $\Omega^*(\cdot)$
which, in turn,
induces a morphism of Mayer-Vietoris
sequences
from that involving the functor  $C$ to that involving
the functor $\Omega^*$.
Thus it suffices to prove that
$$
C(W) @>>> \Omega^*(W)
\tag4.10.5
$$
is an isomorphism on cohomology for contractible $W$.
However, for contractible $W$,
the projective $A$-module $M$ is necessarily free, and 
it suffices to consider $M=A=C^{\infty}(W)$.
Now
$\roman{Ext}^k(A,\Cal O) \cong \roman H_{\roman{cs}}^k(W)$
is zero for $0 \leq k <n$
and a copy of the reals for $k=n$.
Further,
$\roman{Tor}_k^U(\omega_A,A) \cong \roman H^{n-k}(L,A)
\cong \roman H^{n-k}(W)
$
which is again
zero
for $0 \leq k <n$
and a copy of the reals for $k=n$,
and so is
$\roman H^k \Omega^*(W)$.
The chain map
(4.10.5) identifies the two. 
This is, in fact, just property (4.4.1) of a trace. \qed
\enddemo
\smallskip\noindent
{\smc Remark 4.11.}
By means of the  functor
$\roman{Tor}^U_k(\omega_A, \cdot)$,
the argument for (4.10)
reduces the proof of nondegeneracy
of the Poincar\'e duality pairing
to that of
nondegeneracy of the canonical pairing
$$
\roman{Tor}_k^U(C_L, M)
\otimes_{\Bobb R}
\roman{Ext}^k_U(M,\Cal O)
@>>>
\roman{Tor}_0^U(C_L,\Cal O) \cong \Bobb R
$$
between the indicated Tor- and Ext-groups.
Under the circumstances of (4.5),
the nondegeneracy 
of the corresponding pairing (4.5.4)
is immediate
but, 
given an arbitrary duality
$(R,A)$-Lie algebra $L$,
in view of the subtle distinction between the ground ring
$R$ and the algebra $A$ on which $L$ acts in general non-trivially,
nondegeneracy 
of pairings of the kind (4.9.1) and (4.9.2)
is not immediate and perhaps not even always true
though we haven't found a counterexample.
Moreover,
our proof shows that,
under the circumstances of (4.10),
given a left duality $(A,L)$-module $M$,
the group
$\roman{Tor}^U_k(\omega_A, M)$
may be viewed as one of
classes of de Rham $k$-{\it cycles\/}
with values in the flat vector
bundle corresponding to $M$.
Thus the vector spaces
$\roman{Tor}^U_*(\omega_A, M)$
play the role of a kind of de  Rham {\it homology vector spaces\/} for $W$.
For $W$ compact they in fact coincide with the ordinary homology groups
with real coefficients.
For a general (smooth $n$-dimensional) manifold $W$,
$\roman{Tor}^U_n(\omega_A, A) = \roman H_n(L,\omega _A)$,
and 
$\roman H_n(L,\omega _A)$ is one-dimensional, generated
by the fundamental class $e$ (introduced after (3.4),
whether or not $W$ is compact.
This notion of fundamental class comes down to the ordinary one 
for compact $W$ 
but, for non-compact $W$, our theory thus still provides
a fundamental class, as does Borel-Moore homology theory.

\smallskip\noindent
{\smc Remark 4.12.}
The argument for (4.10),
suitably formalized, 
can presumably be used to prove Poincar\'e duality
for other Lie-Rinehart algebras $(A,L)$,
for example for $A$ a regular affine algebra
over a field $k$ and
$L=\roman{Der}(A)$.
 
\smallskip\noindent
{\smc Example 4.13.}
Under the circumstances of (2.12),
let $R= \Bobb R$, and
suppose that $\fra g$ is the Lie algebra of a compact (connected) Lie group
$G$ of dimension $n$ and that
$A$ is the algebra of smooth functions on a smooth manifold $W$
in such a way that the $\fra g$-action on $A$ is induced by
a $G$-action on $W$.
Then, cf. (2.12), the cohomology 
$\roman H^*(L,M)$
of the $(\Bobb R,A)$-Lie algebra
$L = A \otimes _{\Bobb R} \fra g$
with coefficients 
in any left $(A,L)$-module
$M$ is just the ordinary Lie algebra cohomology
$\roman H^*(\fra g,M)$.
Furthermore,
cf. \cite{\ginzwein\ (3.5)}, the Lie algebra cohomology
$\roman H^*(\fra g,A)$
is isomorphic to
$\roman H^*(\fra g,\Bobb R) \otimes_{\Bobb R} A^{\fra g}$.
Thus $\Cal O = A$ and 
the canonical map
$$
t
\colon
\roman H^n(\fra g,A)
\cong \roman H^n(\fra g,\Bobb R) \otimes_{\Bobb R} A^{\fra g}
@>>>
A^{\fra g}
$$
to the $\fra g$-invariants  $A^{\fra g}$
is an isomorphism.
Presumably
$(\Cal O, t)$
yields a trace for $L$,
with
$A^{\fra g}$
as ground ring.
Details have not been verified yet.
Moreover,
over
$A^{\fra g}$,
the cup pairing
$$
\roman H^j(L,A)
\otimes_{A^{\fra g}}
\roman H^{n-j}(L,A)
@>>>
\roman H^n(L,A) \cong
A^{\fra g}
\tag4.13.1
$$
manifestly comes down,
in the classical 
Poincar\'e duality
pairing
$$
\roman H^j(\fra g,\Bobb R)
\otimes_{\Bobb R}
\roman H^{n-j}(\fra g,\Bobb R)
@>>>
\roman H^n(\fra g,\Bobb R) \cong
\Bobb R
\tag4.13.2
$$
in ordinary Lie algebra cohomology,
to 
an extension of scalars from
the reals 
to $A^{\fra g}$.
Since
(4.13.2) is a nondegenerate pairing
of finitely generated real vector spaces,
(4.13.1) is a nondegenerate pairing
of finitely generated free
$A^{\fra g}$-modules.
Thus $L$ satisfies Poincar\'e duality
for $M=A$,
with the algebra  $A^{\fra g}$
of invariants  as ground ring.
Presumably $L$ satisfies Poincar\'e duality
for any duality $(A,L)$-module.
We note that, unless $W$ reduces to a point,
the {\it nondegeneracy\/} of (4.13.2) {\it cannot even be phrased
over the reals\/}.
We also note that,
when the $\fra g$-action is not induced by a
$G$-action,
$\roman H^*(\fra g,A)$
need no longer be isomorphic to
$\roman H^*(\fra g,\Bobb R) \otimes_{\Bobb R} A^{\fra g}$.
A counterexample
is given in \cite\ginzwein.
Perhaps a counterexample to Poincar\'e duality
may be found along these lines.
\smallskip
Our next aim is to illustrate Poincar\'e duality
with another class of examples.
For this purpose we shall need
the following.

\proclaim{Lemma 4.14}
Let $(A,L)$ be a duality Lie-Rinehart algebra of rank $n$
over a commutative ring $R$, let $(\Cal O, t)$ be a trace for
$L$, and let $M_1$ and $M_2$ be 
left duality $(A,L)$-modules (i.e. left $(A,L)$-modules which
are finitely generated and projective as $A$-modules).
Then $L$ satisfies Poincar\'e duality for $M=M_1 \oplus M_2$
if and only if it satisfies Poincar\'e duality for
$M_1$ and for $M_2$.
\endproclaim

\demo{Proof}
For $M=M_1 \oplus M_2$,
the $R$-module morphism (4.9.6) may be written in the form
$$
\aligned
\roman H^{n-k}(L,\roman{Hom}_A(M_1,\Cal O))
&\oplus
\roman H^{n-k}(L,\roman{Hom}_A(M_2,\Cal O))
\\
@>>>
&\roman{Hom}_R(\roman H^k(L,M_1),R)
\oplus
\roman{Hom}_R(\roman H^k(L,M_2),R).
\endaligned
\tag4.14.1
$$
Plainly, (4.14.1) is an isomorphism of $R$-modules
if and only if its constituents
$$
\roman H^{n-k}(L,\roman{Hom}_A(M_1,\Cal O))
@>>>
\roman{Hom}_R(\roman H^k(L,M_1),R)
$$
and
$$
\roman H^{n-k}(L,\roman{Hom}_A(M_2,\Cal O))
@>>>
\roman{Hom}_R(\roman H^k(L,M_2),R)
$$
are isomorphisms of $R$-modules. \qed
\enddemo

\smallskip\noindent
{\smc Example 4.15.}
Let again $W$ be a (real) smooth $n$-dimensional manifold,
write $A = C^{\infty}(W)$,
let $\Cal F$ be a  foliation on $W$
of codimension $n-k$
with compact leaves,
let $L_{\Cal F}\subseteq \roman{Vect}(W)$ be the $(\Bobb R, A)$-Lie algebra of
vector fields tangent to the foliation,
and let $\Cal O_{\Cal F}$ be the $A$-module of 
sections of the
orientation bundle of $\Cal F$. 
Being a real line bundle,
the orientation bundle of $\Cal F$ 
inherits a canonical 
flat connection; the latter, in turn,
yields a left $(A,\roman{Vect}(W))$-module 
and hence $(A,L_{\Cal F})$-module structure
on $\Cal O_{\Cal F}$.
Write $B$ for the space of leaves of $\Cal F$ and $C^{\infty}(B)$ 
for its algebra
of smooth functions, that is, $C^{\infty}(B)$ is the algebra of smooth
functions on $W$ which are constant on the leaves.
As an $A$-module,
the dualizing module $C_{L_{\Cal F}}=\roman{Hom}_A(\Lambda_A^k L_{\Cal F},A)$ 
of $L_{\Cal F}$ (cf. 2.10) is 
the space of sections of a line bundle on $W$,
and the tensor product
$C_{L_{\Cal F}} \otimes_A \Cal O_{\Cal F}$ 
is the space of  densities 
along the leaves;
hence the integration map 
$$
C_{L_{\Cal F}} \otimes_A \Cal O_{\Cal F}
@>>> C^{\infty}(B),
\quad
\omega \mapsto \int_{F_b} \omega,
$$
induces a map
$$
\tau\colon
\roman H^k(L_{\Cal F}, \Cal O_{\Cal F})
@>{e \cap \,\cdot\ }>>
C_{L_{\Cal F}} \otimes_U \Cal O_{\Cal F}
@>>> \roman{Map}(B,\Bobb R),
$$
$e\in \roman H_k(L_{\Cal F}, C_{\Cal F})$ being the fundamental class 
of $L_{\Cal F}$, cf. (3.4).
Given a density along the leaves $\omega$,
the resulting real valued function $h$ on $B$
given by $h(b) = \int_{F_b} \omega$
need not be smooth, though.
This happens for example when, roughly speaking,
the \lq\lq volumes\rq\rq\ of the
leaves \lq\lq jump\rq\rq, e.~g. for a Moebius band with a short leaf.
(This observation is due to the referee.)
\smallskip
We now explain an important special case where $\tau$ yields 
a trace (or rather a variant thereof) and
where, furthermore,
$L_{\Cal F}$
satisfies Poincar\'e duality
with appropriate coefficients.
We take
as ground ring 
$R$ the algebra
$C^{\infty}(B)$
rather than the reals.
Suppose that the foliation $\Cal F$ constitutes
a fiber bundle $\eta\colon W  \to B$ with (compact) fiber $F$,
and write
$L_{\eta}= L_{\Cal F}$
and
$\Cal O_{\eta}= \Cal O_{\Cal F}$.
Then
$\tau$ yields  an isomorphism
$$
\tau\colon
\roman H^k(L_{\Cal F}, \Cal O_{\Cal F})
@>>> 
C^{\infty}(B)
\tag4.15.1
$$
of $(C^{\infty}(B))$-modules.
For
$C_{L_{\eta}} \otimes_A \Cal O_{\eta}$ is the space of sections
of a trivial line bundle on $W$;
a 
nowhere vanishing section
$\widetilde \omega$ is a density along the leaves which, on each 
fiber (or leaf) $F_b$,
satisfies $\int_{F_b} \widetilde\omega \ne 0$,
and $f_{\widetilde\omega}: b \mapsto \int_{F_b} \widetilde\omega$
is a smooth nowhere vanishing function on $B$
(this is not necessarily true when the foliation is not
induced from a fiber bundle);
when we multiply $\widetilde \omega$ by the pull back of 
$f_{\widetilde \omega}$,
we obtain a density along the leaves $\omega$ which has as image under $\tau$
the constant function 1 on $B$.
This shows that $\tau$ maps onto $C^{\infty}(B)$.
Injectivity of $\tau$ follows from the observation that the customary 
argument
which,
given a $k$-form $\alpha$ on a smooth $k$-dimensional manifold
$F$ whose integral over
$F$ vanishes,
by integration against suitable paths
yields a $(k-1)$-form
$\beta$
such that
$d\beta = \alpha$
carries over since integration is compatible with parameters.
Thus
given a 
density along the leaves $\sigma$ 
such that
$\int_{F_b} \widetilde\sigma = 0$
for every leaf $F_b$,
there is a
$(k-1)$-form $\beta$ on
$L_{\eta}$ with values in $\Cal O_{\eta}$
such that $d \beta = \sigma$.
Hence $\tau$ is injective.
Consequently
multiplication of $\omega$ by the pull back of a smooth function 
on $B$ yields the inverse mapping of $\tau$.
\smallskip
Let $\pi \colon P \to B$ be 
a principal bundle for $\eta$, having
compact structure group $G$,
let $\zeta\colon V \to F$ be a flat smooth $G$-vector bundle on $F$
or, equivalently, a $G$-local system on $F$,
and let $\zeta_W \colon V_W=P \times _G V \to P \times _G F = W$
be its extension to $W$;
this is a fibered vector bundle on $W$,
endowed with a flat connection defined only for smooth vector 
fields tangent to the fibers, that is, for elements of $L_\eta$.
Write $M=\Gamma(\zeta_W) $ for its space of sections;
it inherits  a left $(A,L_{\eta})$-module structure.
In particular,  $\zeta$ might be the trivial vector bundle on $F$;
in this case,
$M$ is a finitely generated free $A$-module,
with the obvious 
left $(A,L_{\eta})$-module structure.

\proclaim{Theorem 4.15.3}
The $(\Bobb R,A)$-Lie algebra $L_{\eta}$,
endowed with the pretrace $(\Cal O_{\eta}, \tau)$,
satisfies Poincar\'e duality for every such $M$.
More precisely,
$$
\cup
\colon
\roman H^\ell(L_{\eta},M) \otimes_{C^{\infty}(B)} 
\roman H^{k-\ell}(L_{\eta},\roman{Hom}_A(M,\Cal O_{\eta}))
@>>>
\roman H^k(L_{\eta},\Cal O_{\eta}) \cong C^{\infty}(B)
\tag4.15.4
$$
is a nondegenerate pairing of finitely generated
projective $(C^{\infty}(B))$-modules.
In particular, 
$$
\cup
\colon
\roman H^\ell(L_{\eta},A) \otimes_{C^{\infty}(B)} 
\roman H^{k-\ell}(L_{\eta},\Cal O_{\eta})
@>>>
\roman H^k(L_{\eta},\Cal O_{\eta}) \cong C^{\infty}(B)
\tag4.15.5
$$
is a nondegenerate pairing of finitely generated
projective $(C^{\infty}(B))$-modules.
\endproclaim

We now prepare for the proof of this theorem.
To describe the chain complex 
$\roman{Alt}_A(L_{\eta},M)$,
consider the tangent bundle 
$\tau_{\Cal F}\colon \roman T\Cal F
@>>> W$
of $\Cal F$
(which is assumed to arise from a fibre bundle),
let 
$\Lambda \tau_{\Cal F} = \{
\Lambda^0 \tau_{\Cal F},
\Lambda^1 \tau_{\Cal F},
\dots,
\Lambda^k \tau_{\Cal F} \}
$ 
be the collection of its graded exterior powers,
and consider the system 
$\roman{HOM}(\Lambda \tau_{\Cal F}, \zeta_W)$
of smooth vector bundles on $W$ where
$\roman{HOM}$ stands for morphisms of vector bundles.
Write
$\tau_F \colon \roman TF \to F$ for the tangent bundle of $F$
and let
$A_F = C^{\infty}(F)$,
$L_F = \roman{Vect}(F)$, 
$M_F = \Gamma (\zeta)$;
each $\roman{HOM}(\Lambda^j \tau_{\Cal F}, \zeta_W)$
may clearly be written
$$
P \times _G\roman{HOM}(\Lambda^j \tau_F, \zeta)\colon
P \times _G\roman{HOM}(\Lambda^j \roman TF, V)
@>>> 
P \times _G F = W,
$$
and
the chain complex $\roman{Alt}_A(L_{\eta},M)$
is obtained from
$
(C^{\infty}(P)) \otimes_{\Bobb R} \roman{Alt}_{A_F}(L_F,M_F)
$
when $G$-invariants are taken.
On the other hand,
$\roman{Alt}_{A_F}(L_F,M_F)$
is just the de Rham complex of the fibre $F$
with coefficients in the flat vector bundle $\zeta$,
and its cohomology  
$\roman H^*(F,\zeta)$ inherits  a $G$-module structure; consider
the
associated smooth vector bundle
$$
h^*(\pi,\zeta)\colon P \times _G \roman H^*(F,\zeta) @>>> B.
$$
Its space of sections
$\Gamma (h^*(\pi,\zeta))$
boils down to the space
of $G$-invariants 
of
$$
(C^{\infty}(P)) \otimes_{\Bobb R} \roman H^*(F,\zeta)
\cong \roman H^*((C^{\infty}(P)) \otimes_{\Bobb R} \roman{Alt}_{A_F}(L_F,M_F)).
$$
Further, it is manifest that,
as
$(C^{\infty}(B))$-modules,
the cohomology groups
$\roman H^*(L_{\eta},M)$
and
$\roman H^*(L_{\eta},\roman{Hom}_A(M,\Cal O_{\eta}))$
are projective.
The proof of (4.15.3) now relies on the following.

\proclaim{Lemma 4.15.6}
The canonical map from
$$
\roman H^*(L_{\eta},M) \cong\roman H^*\left( ( (C^{\infty}(P)) 
\otimes_{\Bobb R} 
\roman{Alt}_{A_F}(L_F,M_F))^G\right)
$$
to the space of sections of
$h^*(\pi,\zeta)$
is an isomorphism of
$(C^{\infty}(B))$-modules.
\endproclaim

\demo{Proof}
Write
$C=(C^{\infty}(P)) \otimes_{\Bobb R} \roman{Alt}_{A_F}(L_F,M_F)$
for short;
the group $G$ being compact,
the standard argument
involving invariant integration on $G$
shows that
the canonical map
from
$\roman H^*(C^G)$
to
$(\roman H^*C)^G$
is an isomorphism, that is to say,
the canonical map
$$
\roman H^*(L_{\eta},M)
@>>>
\left((C^{\infty}(P)) \otimes_{\Bobb R} \roman H^*(F,\zeta)\right)^G
=\Gamma (h^*(\pi,\zeta))
$$
is an isomorphism
of $(C^{\infty}(B))$-modules. \qed
\enddemo

\demo{Proof of Theorem {\rm (4.15.3)}}
Lemma 4.15.6 reduces
the question whether
the pairing (4.15.4) is nondegenerate
to ordinary Poincar\'e duality for $F$
with coefficients in $\zeta$.
Alternatively,
we may cover $F$ by finitely many $G$-invariant open 
contractible subsets
$V_1,\dots, V_\ell$
such that each non-empty intersection 
$V_{j_1} \cap \dots \cap V_{j_k}$
is itself contractible
and,
on any such non-empty intersection, 
work with  those sections 
of the orientation bundle
of the corresponding foliation
which, on each leaf,
are compactly supported.
Poincar\'e duality may then be established
by an   argument generalizing that for
(4.10) above. \qed
\enddemo

Let $(A,L)$ be a Lie-Rinehart algebra over a general commutative ring $R$.
A left $(A,L)$-module
$M$ will be said to be an {\it induced\/}
left $(A,L)$-module
provided there is an $A^L$-module $M'$ such that,
as a left $(A,L)$-module,
$M$ coincides with
$A \otimes_{A^L} M'$
where
$A \otimes_{A^L} M'$
is endowed with  the obvious
left $(A,L)$-module structure induced by that on $A$.

\proclaim{Corollary 4.15.7}
The $(\Bobb R,A)$-Lie algebra $L=L_{\eta}$,
endowed with the pretrace $(\Cal O_{\eta}, \tau)$,
satisfies Poincar\'e duality for every
induced left $(A,L_{\eta})$-module
$M$ of the kind $M=A \otimes_{A^L} M'$
where $M'$ is a finitely generated projective
$A^L$-module.
That is to say,
$$
\cup
\colon
\roman H^\ell(L_{\eta},M) \otimes_{C^{\infty}(B)} 
\roman H^{k-\ell}(L_{\eta},\roman{Hom}_A(M,\Cal O_{\eta}))
@>>>
\roman H^k(L_{\eta},\Cal O_{\eta}) \cong C^{\infty}(B)
\tag4.15.8
$$
is a nondegenerate pairing of finitely generated
projective $(C^{\infty}(B))$-modules.
\endproclaim

\demo{Proof}
By (4.15.3),
the pairing (4.15.4)
is nondegenerate when $M$ is a finitely generated
free $A$-module, endowed with the obvious
left $(A,L)$-module structure.
The general case is then a consequence of
Lemma 4.14. \qed
\enddemo

\medskip\noindent
{\bf 5. Extensions of Lie-Rinehart algebras and Poincar\'e duality}
\smallskip\noindent
In this section we shall show that Poincar\'e duality
is preserved under extensions 
of Lie-Rinehart algebras.
In Section 7 below, this will enable us to show that, in certain cases,
Poisson (co)homology
satisfies Poincar\'e duality.
\smallskip
Let
$$
0
@>>>
L'
@>>>
L
@>>>
L''
@>>>
0
\tag5.1
$$
be an extension of duality $(R,A)$-Lie algebras;
let $k= \roman{rank}(L'),
\, m = \roman{rank}(L''),
\, n = k+m= \roman{rank}(L)$.
Then $L'$ is
an ordinary Lie algebra in the category of $A$-modules
which,
as an $A$-module, is finitely generated and projective 
of constant rank $k$.
Hence, cf. (4.5), it has a canonical trace
$(\Cal O', t')$
and satisfies {\it classical\/} Poincar\'e duality for every left
$L'$-module, cf. (4.5) above.

\proclaim{Proposition 5.2}
A trace for $L''$
induces 
a trace for $L$.
\endproclaim
The precise relationships between the traces for
$L$ and $L''$ will be given in (5.2.1) and (5.2.2) below.

\demo{Proof} 
Let
$(\Cal O_{L''},t'')$
be a trace for $L''$,
with 
$t''\colon C_{L''} \otimes _{U''} \Cal O_{L''} @>{\cong}>> R$,
and, cf. (4.6), let
$V'' = C_{L''} \otimes _A \Cal O_{L''}$,
endowed with the  right $(A,L'')$-module structure (2.3).
Furthermore,
let
$$
\Cal O_L = \roman{Hom}_A(C_L, C_{L''}) \otimes_A \Cal O_{L''}
\tag5.2.1
$$
with left $(A,L)$-module structure on
$\roman{Hom}_A(C_L, C_{L''})$
given by (2.3) and with the ordinary tensor product
left $(A,L)$-module structure on
$\Cal O_L$,
and write
$V$ for $V''$, viewed as a right
$(A,L)$-module via the projection from $L$ to $L''$.
Then the canonical maps
$$
C_L \otimes _U (\roman{Hom}_A(C_L, C_{L''}) \otimes_A \Cal O_{L''}) 
@>>>
(C_L \otimes _ A\roman{Hom}_A(C_L, C_{L''})) \otimes_U \Cal O_{L''} 
$$
and
$$
(C_L \otimes _ A\roman{Hom}_A(C_L, C_{L''})) \otimes_U \Cal O_{L''} 
@>>>
C_{L''} \otimes_{U''} \Cal O_{L''}
$$
are isomorphisms (of $R$-modules)
and the composite
$$
t\colon C_L \otimes _U \Cal O_L 
@>>>
C_{L''} \otimes_{U''} \Cal O_{L''} 
@>{t''}>> R
\tag5.2.2
$$
is an isomorphism of $R$-modules
which in fact yields a trace for $L$.
To verify the defining property
of a trace,
we use (4.6)(ii):
Let $N$ be
a right
duality $(A,L)$-module.
Then
$$
\roman{Hom}_U(N,V)
\cong 
\roman{Hom}_{U''}(N\otimes_{U'}A,V'')
$$
and
$$
\roman{Hom}_R(N\otimes_U A,V \otimes _U A)
\cong 
\roman{Hom}_R((N\otimes_{U'}A)\otimes_{U''}A,V''\otimes_{U''}A).
$$
Since
$(\Cal O_{L''},t'')$
is a trace for $L''$
and since
$N$ is a right duality $(A,L)$-module,
by (4.6) (i), the canonical morphism
$$
\roman{Hom}_{U''}(N\otimes_{U'}A,V'')
@>>>
\roman{Hom}_R((N\otimes_{U'}A)\otimes_{U''}A,V''\otimes_{U''}A)
$$
is an isomorphism.
Consequently the canonical morphism
$$
\roman{Hom}_U(N,V)
@>>>
\roman{Hom}_R(N\otimes_U A,V \otimes _U A)
$$
is an isomorphism.
Since $N$ was arbitrary,
by (4.6) (ii),
$(\Cal O, t)$ indeed yields a trace for $(A,L)$.
\qed
\enddemo

\smallskip 

It is readily seen that,
for a
left $(A,L)$-module $M$,
the canonical $L$-action
on the chain complex
$\roman{Alt}_A(L',M)$
(over the ground ring $R$)
makes $\roman{Alt}_A(L',M)$ into
a chain complex in the category of
left
$(A,L)$-modules,
and 
the cohomology
$\roman H^*(L',M)$
thus inherits 
a left $(A,L'')$-module structure.
Likewise,
for a
right $(A,L)$-module $N$,
the  $(A,L)$-action
on the chain complex
$N \otimes _A K(A,L')$
(over the ground ring $R$)
given by (2.4) in each degree
makes
$N \otimes _A K(A,L')$
into
a chain complex in the category of
right
$(A,L)$-modules,
and 
the homology
$\roman H_*(L',N)$
inherits a right
$(A,L'')$-module 
structure.

\proclaim{Theorem 5.3}
Let $(\Cal O'', t'')$
be a trace for $L''$, and
let $(\Cal O, t)$
be the corresponding trace for $L$ given by {\rm (5.2)}.
Given a right $(A,L)$-module $N$,
if the $(R,A)$-Lie algebra $L''$ satisfies Poincar\'e duality for 
the right $(A,L'')$-modules
$$
\roman H_0(L',N),\roman H_1(L',N),\dots, \roman H_k(L',N),
$$
then the $(R,A)$-Lie algebra $L$ satisfies Poincar\'e duality for the right
$(A,L)$-module $N$.
\endproclaim

We note that the theorem may also be phrased for
a left $(A,L)$-module (instead of the right $(A,L)$-module $N$)
but the wording would be technically more involved.
The translation is by means of the standard device:
Given a left
$(A,L)$-module $M$, let $N = C_L \otimes _A M$,
with 
right $(A,L)$-module structure (2.4);
given a right
$(A,L)$-module $N$, let $M = \roman{Hom}_A(C_L, N)$,
with 
left $(A,L)$-module structure (2.3).
We leave the details of the translation of the wording of (5.3) 
into a statement involving left modules rather than right ones
to the reader.

\demo{Proof}
Let $V = C_L \otimes _A O$,
with  right $(A,L)$-module
structure (2.4).
Consider the cohomology change of rings spectral sequence 
$(\roman E^{p,q}_r(\roman{Hom}_A(N,V)),d_r)$
for the canonical surjection
$U(A,L) \to U(A,L'')$,
with coefficients in the left $(A,L)$-module 
$\roman{Hom}_A(N,V)$,
with
left $(A,L)$-module structure (2.3).
See $(2)_4$ in (XVI.5) (p. 349) of \cite\cartanei.
This spectral sequence has
$$
\roman E^{p,q}_2(\roman{Hom}_A(N,V)) =
\roman H^p(L'', \roman H^q(L',\roman{Hom}_A(N,V))),
\quad p,q \geq 0.
$$
Likewise,
with coefficients in the right $(A,L)$-module $N$,
the homology change of rings spectral sequence 
$(\roman E_{p,q}^r(N),d^r)$
for
the canonical surjection
$U(A,L) \to U(A,L'')$
($(2)_2$ in (XVI.5) (p. 348) of \cite\cartanei)
has
$$
\roman E_{p,q}^2(N) =
\roman H_p(L'', \roman H_q(L',N)),\quad p,q \geq 0.
$$
Moreover,
for every $q \geq 0$,
the cap pairing
$$
\roman H_q(L',N)
\otimes_A
\roman H^q(L',\roman{Hom}_A(N,V))
@>>>
V
\tag5.3.1
$$
for the ordinary $A$-Lie algebra $L'$
is a pairing 
of $A$-modules and,
by classical duality for $L'$,
cf. (4.5),
the adjoint
$$
\roman H^q(L',\roman{Hom}_A(N,V))
@>>>
\roman{Hom}_A(\roman H_q(L',N),V)
\tag5.3.2
$$
of (5.3.1)
is an isomorphism of $A$-modules.
Furthermore,
for $q \geq 0$,
$\roman H_q(L',N)$
and
$\roman H^q(L',\roman{Hom}_A(N,V))$
inherit  
right- and
left $(A,L'')$-module structures, respectively,
in such a way that,
when
$\roman{Hom}_A(\roman H_q(L',N),V)$
is endowed with the left
$(A,L'')$-module
structure (2.3),
(5.3.2)
is in fact an isomorphism of
left
$(A,L'')$-modules.
Consequently, given $p,q \geq 0$, the 
cap pairing
$$
\roman H_p(L'',\roman H_q(L',N))
\otimes_R
\roman H^p(L'',\roman H^q(L',\roman{Hom}_A(N,V)))
@>>>
\roman H_0(L'', V'') \cong R
\tag5.3.3
$$
for $L''$, with reference to
the bilinear pairing (5.3.1),
amounts to the corresponding Poincar\'e duality pairing
(4.9.4)
(with $L''$ and
$\roman H_q(L',N)$ instead of
$L$ and $N$, respectively); it 
has the form
$$
\roman H_p(L'',\roman H_q(L',N))
\otimes_R
\roman H^p(L'', \roman{Hom}_A(\roman H_q(L',N),V))
@>>>
\roman H_0(L'', V'') \cong R.
\tag5.3.4
$$
By hypothesis,
$L''$ satisfies Poincar\'e duality
for the right
$(A,L'')$-modules
$$
\roman H_0(L',N),\dots,\roman H_k(L',N).
$$
Hence, given $p,q \geq 0$,
the adjoint
$$
\roman H^p(L'', \roman{Hom}_A(\roman H_q(L',N),V))
@>>>
\roman{Hom}_R(\roman H_p(L'',\roman H_q(L',N)),R)
$$
of (5.3.4) is an isomorphism of $R$-modules whence
the adjoint
$$
\roman H^p(L'',\roman H^q(L',\roman{Hom}_A(N,V)))
@>>>
\roman{Hom}_R(\roman H_p(L'',\roman H_q(L',N)),R)
\tag5.3.5
$$
of the
cap pairing 
(5.3.3)
is an isomorphism of $R$-modules.
\smallskip
The description of the change of rings spectral sequences
in \cite\cartanei\ 
involves suitable bicomplexes arising from appropriate resolutions.
For the present circumstances,
we now indicate briefly
a different description
which is more appropriate
for our purposes:
The projection from $L$ to $L''$ induces
a surjection 
from $\Lambda_A L$ onto
$\Lambda_A L''$,
whence the degree filtration
of
$\Lambda_A L''$
induces a filtration $F_0 \subseteq F_1 \subseteq ... $ of
$\Lambda_A L$ where
an element $\alpha_1 \wedge \alpha_2 \wedge \dots \wedge \alpha_k$
of
$\Lambda_A L$
lies in $F_p$ if 
among the $\alpha_1, \alpha_2, \dots, \alpha_k$ in $L$
at most $p$
have non-zero image in $L''$.
This filtration induces a filtration of the Rinehart complex
$K(A,L)$ and, therefore,
filtrations of the chain  complexes
$$
N \otimes_U K(A,L)
\quad
\text{and}
\quad
\roman{Hom}_U(K(A,L), \roman{Hom}_A(N,V)).
$$
The spectral sequences
$
(\roman E_{p,q}^r(N),d^r)$
and
$(\roman E^{p,q}_r(\roman{Hom}_A(N,V)),d_r)$
arise 
from these filtrations.
The cap pairing 
between the relevant chain complexes
is compatible with the
filtrations and hence,
in view of standard  spectral sequence
multiplicative properties, induces
a pairing
$$
(\roman E_{p,q}^r(N),d^r)
\otimes_R
(\roman E^{p,q}_r(\roman{Hom}_A(N,V)),d_r)
@>>>
\roman E_{0,0}^r(V),
\quad p,q \geq 0,\ r \geq 2,
\tag5.3.6
$$
of spectral sequences.
However,
since $t$ is a trace for $L$,
the canonical morphisms
$$
\roman H_0(L'',\roman H_0(L',V))
=
\roman E_{0,0}^2(V)
@>>>
\dots
@>>>
\roman E_{0,0}^r(V)
$$
are all isomorphisms and the trace induces canonical isomorphisms
between the $\roman E_{0,0}^r(V)$'s
and the ground ring $R$.
Consequently
the pairing
(5.3.6) of spectral sequences may be written
$$
(\roman E_{p,q}^r(N),d^r)
\otimes_R
(\roman E^{p,q}_r(\roman{Hom}_A(N,V)),d_r)
@>>>
R,\quad p,q \geq 0,\ r \geq 2.
\tag5.3.7
$$
Furthermore,
since (5.3.5) is an isomorphism,
the pairing (5.3.5)
has the property that, for $r=2$,
the adjoint
$$
\roman E^{p,q}_2(\roman{Hom}_A(N,V))
@>>>
\roman{Hom}_R(\roman E_{p,q}^2(N),R)
$$
is an isomorphism of bigraded $R$-modules.
Hence the adjoint
$$
(\roman E^{p,q}_r(\roman{Hom}_A(N,V)),d_r)
@>>>
\roman{Hom}_R((\roman E_{p,q}^r(N),d^r),R)
=
\left(\roman{Hom}_R((\roman E_{p,q}^r(N),R),(d^r)^*\right)
$$
of
the pairing (5.3.7) 
is an isomorphism of spectral sequences in the category of $R$-modules.
Thus
the adjoint
$$
\roman E^{p,q}_{\infty}(\roman{Hom}_A(N,V))
@>>>
\roman{Hom}_R(\roman E_{p,q}^{\infty}(N),R),
\quad p,q \geq 0,
$$
of the induced pairing
$$
\roman E_{p,q}^{\infty}(N)
\otimes_R
\roman E^{p,q}_{\infty}(\roman{Hom}_A(N,V))
@>>>
R
\tag5.3.8
$$
is an isomorphism of bigraded $R$-modules.
Consequently, for every $\ell\geq 0$,
the 
adjoint
$$
\roman H^{\ell}(L,\roman{Hom}_A(N,V))
@>>>
\roman{Hom}_R(\roman H_{\ell}(L,N), R)
\tag5.3.9
$$
of the cap pairing
$$
\roman H_{\ell}(L,N)
\otimes_R
\roman H^{\ell}(L,\roman{Hom}_A(N,V))
@>>>
\roman H_0(L, V) \cong R
\tag5.3.10
$$
is an isomorphism of $R$-modules.
This proves the claim. \qed
\enddemo

We now illustrate the last result with a special class of examples.
Let
the ground ring
be that of the reals, $\Bobb R$,
and let
$A$ be the algebra of smooth
functions on
a smooth manifold $W$.
Recall that 
a Lie algebroid on $W$,
cf. e.~g. \cite\mackbook,
determines a
duality $(\Bobb R,A)$-Lie algebra,
and
every
duality $(\Bobb R,A)$-Lie algebra
arises in this way.
By (5.2),
given an 
$(\Bobb R,A)$-Lie algebra $L$
which arises from a transitive
Lie algebroid,
when 
$(\Cal O_A, t_A)$
refers to the  canonical trace
for $\roman{Der}(A) = \roman{Vect}(W)$,
cf. (4.7),
the canonical isomorphism
between $C_L \otimes _A(Q_L \otimes _A \Cal O_A)$
and
$\omega_A \otimes _A \Cal O_A$
yields a trace
$$
(\Cal O_L, t_L) =(Q_L \otimes _A \Cal O_A, t_L)
$$
for $L$.

\proclaim{Theorem 5.4}
Every 
(duality) Lie-Rinehart algebra
$(A,L)$
over the reals
arising from a transitive
Lie algebroid on a smooth manifold
satisfies
Poincar\'e duality
for any
duality $(A,L)$-module.
\endproclaim

\demo{Proof}
This is an immediate consequence of (5.3),
combined with 
Poincar\'e duality for the
$(\Bobb R,A)$-Lie algebra of smooth vector fields
on the manifold which comes into play, i.~e. with
standard Poincar\'e duality
in ordinary de Rham cohomology
of smooth manifolds,
cf. (4.10). \qed
\enddemo

For an {\it integrable\/}
transitive Lie algebroid, i.~e. one arising from a principal bundle,
this result follows of course from ordinary Poincar\'e duality for
the total space of the principal bundle.
The statement of the theorem is more general, though, since there
are transitive Lie algebroids which do not arise from a principal bundle
\cite\almemoli.                                        

\medskip\noindent
{\bf 6. The Picard group and modular class of a Lie-Rinehart algebra}
\smallskip\noindent
Let $L$ be an $(R,A)$-Lie algebra.
We shall denote the $A$-module 
of derivations of $L$ in $A$, that is, that of
$A$-valued 1-cocycles, by
$\roman{Der}(L,A)$.
We remind the reader that $\phi\in \roman{Hom}_A(L,A)$
is called a {\it derivation\/}
provided
$$
\phi([\alpha,\beta])
=
\alpha \phi(\beta) - \beta \phi(\alpha),
\quad \alpha, \beta \in L.
$$
Let $M$ be an $A$-module.
Recall that an
$L$-{\it connection\/}
on $M$
may be described as a pairing
$L\otimes_R M \to M$
of $R$-modules, written
$(\alpha,x) \mapsto \alpha(x)$,
such that 
$$
\align
\alpha (ax) &= \alpha(a) x + a \alpha(x),
\\
(a\alpha)(x) &= a(\alpha (x)),
\endalign
$$
where $a \in A,\ x \in M, \ \alpha \in L$;
an $L$-connection on $M$
which is actually
a {\it left\/}
$(A,L)$-module structure
is also said to be {\it flat\/}.
See \cite\poiscoho\ (2.16)
for historical remarks on these algebraic notions of connection etc.
Recall \cite\poiscoho\ (Section 2) that
any projective $A$-module has an $L$-connection.
Let now $M$ be a projective rank one $A$-module
admitting a flat 
$L$-connection, that is, a
left $(A,L)$-module structure.
We fix such a structure 
$L \to \roman{End}_R(M)$
and denote 
by $\alpha_M$
the operator
on $M$ which is
the image of $\alpha \in L$.                           
The next statement generalizes
the well  known fact that,
for   a flat
line bundle
on a smooth manifold $W$,
the group of de Rham 1-cocycles
$Z^1(W, \Bobb R)$
acts 
faithfully and transitively 
on the space of flat connections.

\proclaim{Proposition 6.1}
The $A$-module $M$ being projective of rank one,
the assignment to
$\phi \in \roman{Der}(L,A)$
and $\alpha \in L$
of
$$
\alpha_{M,\phi} \colon M \to M,\quad
\alpha_{M,\phi}(m) = \alpha_M (m)+ \phi(\alpha) m,
\quad m \in M,
$$
yields a faithful and transitive action of
$\roman{Der}(L,A)$
on the set of 
left $(A,L)$-module structures on $M$.
\endproclaim

\demo{Proof} This is straightforward and left to the reader. \qed
\enddemo

For $M=A$, endowed with the $L$-action which is part
of the  $(R,A)$-Lie algebra structure (of $L$),
the statement of (6.1) comes down to the well known fact that
every derivation $\phi$ of $L$ in $A$ induces a  
left $(A,L)$-module structure on $A$ which is entirely characterized by
$\alpha (1) = \phi (\alpha)$, where $\alpha \in L$.
Given a derivation $\phi$ of $L$ in $A$,
the general statement of (6.1) is then obtained
when  $M$ is canonically identified with $M \otimes_A A$,
the tensor product $M \otimes_A A$ 
being endowed with the left $(A,L)$-module structure (2.1).
\smallskip
What corresponds to the group of gauge transformations
for a line bundle is now the group
of units $A^{\times}$ in $A$;
indeed, when $A$ is the algebra of smooth functions on a smooth
manifold,
$A^{\times}$ is precisely the multiplicative group of nowhere-zero functions,
that is, that of functions with values in the multiplicative group
$\Bobb R^{\times}$
of non-zero real numbers.
In our general situation,
the algebra $\roman{End}_A(M)$
is canonically isomorphic to $A$, whence the group
$\roman{Aut}_A(M)$ amounts to $A^{\times}$.
It acts on 
$\roman{Der}(L,A)$
by the association
$$
\phi \mapsto {}^u\phi,\quad
{}^u\phi = \phi - u^{-1} du,
\quad \phi \in \roman{Der}(L,A),
\ u \in A^{\times},
$$
where as usual $du\in \roman{Der}(L,A)$ refers to the
1-coboundary $du(\alpha) = \alpha (u)$, for
$\alpha \in L$.
We denote the space of $A^{\times}$-orbits
in
$\roman{Der}(L,A)$
by $\Cal H^1(L,A)$.

\proclaim{Proposition 6.2}
The 
space of $A^{\times}$-orbits
$\Cal H^1(L,A)$
inherits an abelian group
structure
from that of $\roman{Der}(L,A)$.
\endproclaim
\smallskip
\noindent
{\smc Remark 6.3.}
In the gauge theory situation, 
$\Cal H^1(L,A)$ is just the abelian group
of flat line bundles which are trivial as line bundles.
Given a smooth real manifold $W$,
this group amounts to the first 
cohomology
group
$\roman H^1(W,\Bobb R_+^{\times})$
with coefficients in the multiplicative group
of positive real numbers
$\Bobb R_+^{\times}$.
In fact, the holonomy induces a map from the  space
of flat connections on a trivial line bundle on $W$
to the space 
$$
\roman{Hom}(\pi_1(W),\Bobb R_+^{\times})
\cong
\roman H^1(W,\Bobb R_+^{\times}),
$$
and this map induces an isomorphism
from
$\Cal H^1(L,A)$
onto
$\roman H^1(W,\Bobb R_+^{\times})$,
$A$ and $L$ being the algebra of smooth functions and Lie algebra
of smooth vector fields, respectively, on $W$.
Since the structure group,
$\Bobb R_+^{\times}$, is abelian,
this map is just given by the assignment to a smooth 1-form
$\alpha$ of 
the homomorphism
$\phi_\alpha \colon \pi_1(W) \to \Bobb R_+^{\times}$
given by
$$
\phi_\alpha [c] = \roman{exp}\int_c \alpha,
$$
perhaps with a minus sign (depending on the choice
of convention),
where $c$ is a closed curve in $W$ representing an
element of $\pi_1(W)$.
Thus,
for an general Lie-Rinehart algebra $(A,L)$,
 the group
$\Cal H^1(L,A)$ 
generalizes
the cohomology group
$\roman H^1(W,\Bobb R_+^{\times})$.
Via the logarithm, this group amounts to the more traditional
$\roman H^1(W, \Bobb R)$.

\demo{Proof of {\rm (6.2)}}
Let $\phi,\phi', \psi, \psi'\in\roman{Der}(L,A)$,
and suppose that
there are $u,v \in A^{\times}$
such that
$$
\phi' = \phi-u^{-1}du,
\quad
\psi' = \psi-v^{-1}dv. 
$$
Then
$$
\phi' +\psi'= \phi +\psi-(uv)^{-1}d(uv). 
$$
Hence
the space of $A^{\times}$-orbits constitutes an abelian group. \qed
\enddemo

Next we realize this group, in fact, construct a larger group
containing it as a subgroup.
Recall that $\roman{Pic}(A)$ is the abelian group
of isomorphism classes of projective rank one $A$-modules,
the group structure being induced by the tensor product.
We denote by $\roman{Pic}^{\roman{flat}}(L,A)$
the set of isomorphism classes of left
$(A,L)$-modules which are
projective of rank one as $A$-modules.
We mention in passing
that 
we have 
written
$\roman{Pic}^{\roman{flat}}(L,A)$
rather than just
$\roman{Pic}(L,A)$
since, 
when $L$ has suitable additional structure, there is presumably
also a group
of projective of rank one $A$-modules
with \lq\lq harmonic\rq\rq\ $L$-connections
generalizing the group of line bundles with harmonic connection;
this group is then to be denoted by
$\roman{Pic}(L,A)$.
The operation of tensor product endows 
$\roman{Pic}^{\roman{flat}}(L,A)$
with an abelian group structure.
The inverse of
the class of such an $(A,L)$-module $M$ is given by
the class of 
the $(A,L)$-module $\roman{Hom}_A(M,A)$,
with left
$(A,L)$-module structure given by
(2.2).
The zero element of this group
is 
the class of the free rank one $A$-module,
with $L$-action given by
$\alpha (b) = 0$, where
$\alpha \in L$ and $b$ is a basis element.
Notice in particular that an automorphism
of a
left $(A,L)$-module which is
projective of rank one as an $A$-module
may be viewed as a gauge transformation.

\proclaim{Proposition 6.4}
For a free rank one $A$-module $M$, with basis $b$,
the assignment to
$\phi \in \roman{Der}(L,A)$
and $\alpha \in L$
of
$$
\alpha_{M,\phi} \colon M \to M,\quad
\alpha_{M,\phi}(b) = \phi(\alpha) b
$$
identifies the
kernel of the obvious forgetful  map
from $\roman{Pic}^{\roman{flat}}(L,A)$
to
$\roman{Pic}(A)$ with the group
$\Cal H^1(L,A)$.
\endproclaim

\demo{Proof}
Given
$\phi \in \roman{Der}(L,A)$,
write $M_\phi$ for 
the free rank one $A$-module with basis $b$
and left $(A,L)$-structure given by
$$
\alpha_{M,\phi} \colon M \to M,\quad \alpha \in L.
$$
For
$\phi, \psi \in \roman{Der}(L,A)$,
on $M_\phi \otimes_AM_\psi$,
the corresponding 
left $(A,L)$-module structure
(2.1)
is given by
$$
\alpha(b \otimes_Ab)
=
(\alpha_{M,\phi}b) \otimes_Ab
+
b \otimes_A(\alpha_{M,\psi} b)
=
(\phi(\alpha) b) \otimes_Ab
+
b \otimes_A (\psi(\alpha)b)
$$
which
amounts to
$$
\alpha_{M\otimes_AM,\phi+\psi}\colon
(b \otimes_A b) \longmapsto
(\phi(\alpha) +\psi(\alpha))(b \otimes_Ab)
$$
whence the assertion. \qed
\enddemo

Summing up and combining with 
\cite\poiscoho\  (2.15.1), we obtain the following.

\proclaim{Theorem 6.5}
The assignment
to a projective rank one $A$-module 
of its characteristic class in 
$\roman H^2(L,A)$
yields a map from
$\roman{Pic}(A)$ to $\roman H^2(L,A)$ which
fits into the exact sequence
$$
0
@>>>
\Cal H^1(L,A)
@>>>
\roman{Pic}^{\roman{flat}}(L,A)
@>>>
\roman{Pic}(A)
@>>>
\roman H^2(L,A).
\tag6.5.1
$$
\endproclaim

We now suppose that, as an $A$-module,
the $(R,A)$-Lie algebra
$\roman{Der}(A)$
is finitely generated and projective
of constant rank.
For simplicity, we then refer to $A$ as being
{\it regular\/}.
Extending notation already used in (4.7),
we denote the dualizing module
$\roman {Hom}_A(\Lambda_A^{\roman{top}} \roman{Der}(A), A)$
of $\roman{Der}(A)$
by $\omega_A$.
Here
$\Lambda_A^{\roman{top}}$
refers to the highest non-zero exterior power of
$\roman{Der}(A)$;
its degree equals
the rank of $\roman{Der}(A)$.
Our notion of regularity
is consistent with standard terminology
in algebraic geometry:
When the module $\roman D_A$ of formal differentials
is finitely generated and projective
of constant rank as an $A$-module,
$\roman{Der}(A)$
will likewise have this property,
and the requirement that
$\roman D_A$ be finitely generated and projective
of constant rank essentially 
amounts to the usual notion of regularity in algebraic geometry.
Further, we can then identify
$\omega_A$ with the highest non-zero exterior power 
$\Lambda_A^{\roman{top}} \roman D_A$ of $D_A$;
in algebraic geometry, 
a closely related
object is called
the \lq\lq canonical sheaf\rq\rq\ 
and denoted by $\omega$ with an appropriate subscript
whence the notation.
Likewise, cf. (4.7), when
$A$ is the algebra of smooth functions
on a smooth real manifold $W$,
$\omega_A$ amounts to
the highest non-zero exterior power of the 
space of sections of the
cotangent bundle of $W$. As observed in (2.8),
the operation of Lie derivative (2.6.2) endows
the projective rank one $A$-module
$\omega_A$
with  a right $(A,\roman{Der}(A))$-module structure.
\smallskip
Let $L$ be a duality
$(R,A)$-Lie algebra
of rank $n$,
and let
$C_L$  
($=\Lambda_A^n L^*$)
be its dualizing module.
Let
$$
Q_L = \roman{Hom}_A(C_L,\omega_A);
\tag6.6
$$
this is a projective rank one $A$-module
which has a canonical left $(A,L)$-module structure
given by (2.3).
We refer to
the class $[Q_L]$ of
$Q_L$
in 
$\roman{Pic}^{\roman{flat}}(L,A)$
as the {\it modular class\/}
of $L$.
A version of the module $Q_L$ occurred already in Section 5:
when the structure map $L \to  \roman{Der}(A)$
is surjective,
and when $\roman{Der}(A)$ has a trace,
written $(\Cal O_A, t_A)$,
the module written
$\roman{Hom}_A(C_L,C_{L''})$ in (5.2.1)
is precisely of the kind $Q_L$,
where now $L'' = \roman{Der}(A)$;
letting $\Cal O_L= Q_L \otimes _A \Cal O_A$,
with
left $(A,L)$-module structure (2.1),
we then obtain a trace for
$L$ by means of
the canonical isomorphism
$$
C_L \otimes _U \Cal O_L
@>>>
\omega_A \otimes _{U''} \Cal O_A
@>{t_A}>>
R,
$$
cf. (5.2.2).
\smallskip
When $A$ and
$L$ 
are
the 
algebra of  smooth real functions
and
$(\Bobb R,A)$-Lie algebra
of smooth vector fields, respectively, on a smooth real manifold
$W$, 
as an $A$-module,
$$
Q_L = 
\roman{Hom}_A(C_L,\omega_A)
=\roman{Hom}_A(\omega_A,\omega_A),
$$
that is,
$Q_L$ is free of rank one,
generated by the identity map.
Furthermore,
its
left $(A,L)$-module structure is \lq\lq trivial\rq\rq\ 
in the sense that it has a basis, here the identity map,
which remains invariant under the
$L$-action;
in fact, this means that
$Q_L$
and $A$ are isomorphic
even as
left $(A,L)$-modules,
and the modular class $[Q_L]$ is trivial.

\smallskip
\noindent
{\smc Remark 6.7.}
Let $R=\Bobb R$,
let $W$ be a smooth real manifold,
let $A$ be its algebra of  smooth real functions,
and let
$L$ be 
the $(\Bobb R,A)$-Lie algebra
coming from a Lie algebroid on $W$.
The group
$\Cal H^1(L,A)$ 
is  the space
of $A^{\times}$-orbits
in $\roman {Der}(L,A)$, 
and the passage to
$A^{\times}$-orbits amounts to the identification of
two 1-cocycles $\phi$ and $\psi$ whenever
$\phi-\psi = u^{-1}du$ for some $u \in A^{\times}$.
Under the present circumstances,
a unit $u$ boils down to a nowhere vanishing function $f$,
and
$$
u^{-1}du
=f^{-1}df = d \roman{log}|f|
$$
is a coboundary.
Hence the identity map of
$\roman {Der}(L,A)$
induces a canonical map from
$\Cal H^1(L,A)$ 
to
$\roman H^1(L,A)$.
Moreover, the group $\roman{Pic}(A)$ is just  the group
$\roman{Hom}(\pi_1(W),\Bobb Z /2)\cong \roman H^1(W,\Bobb Z /2)$. 
Consequently,
the square $[Q_L]^2$ of
the modular class
in $\roman{Pic}^{\roman{flat}}(L,A)$
lies in
$\Cal H^1(L,A)$ and hence maps 
to
$\roman H^1(L,A)$.
In \cite\evluwein,
this image, divided by 2, is taken as the definition of the modular class
(of the corresponding Lie algebroid).
One of the disadvantages of such a definition is this:
If $M$ is a projective rank one $A$-module
with a left $(A,L)$-structure,
as an $A$-module, its tensor square
$M \otimes_A M$ is free, and hence
the 
tensor square
left $(A,L)$-structure
is given by a derivation $\phi$ of $L$ in $A$ as in (6.4);
when we then endow the free rank one $A$-module 
with the left $(A,L)$-structure
which corresponds to the 
derivation $\frac 12 \phi$ of $L$ in $A$,
$M$ and the free rank one $A$-module with the
left $(A,L)$-structure
given by
$\frac 12 \phi$
would determine the same element
in
$\roman H^1(L,A)$;
the two 
left $(A,L)$-modules are, of course, distinguished in
$\roman{Pic}^{\roman{flat}}(L,A)$.
Another disadvantage is that such a definition will not work
over other ground rings
like e.~g. that of the complex numbers,
unless the 
projective rank one $A$-module 
under consideration
has finite order in the Picard group.
\smallskip
\noindent
{\smc Remark 6.8.}
Under the circumstances of (6.7),
when $L$ is the
$(\Bobb R,A)$-Lie algebra of smooth vector fields on $W$
the exact sequence
(6.5.1)
boils down to the exact sequence
$$
0
@>>>
\roman H^1(W, \Bobb R_+^{\times})
@>>>
\roman H^1(W,\Bobb R^{\times}) 
@>>>
\roman H^1(W,\Bobb Z /2) 
@>>>
0
$$
associated with the split exact coefficient sequence
$$
0
@>>>
\Bobb R_+^{\times}
@>>>
\Bobb R^{\times} 
@>>>
\Bobb Z /2 
@>>>
0
$$
of abelian groups.
For 
a general 
$(\Bobb R,A)$-Lie algebra
$L$
of the kind coming into play in (6.7),
by naturality,
the corresponding morphism
$L \to \roman{Vect}(W)$
of 
$(\Bobb R,A)$-Lie algebras
induces a commutative diagram
$$
\CD
0
@>>>
\roman H^1(W, \Bobb R_+^{\times})
@>>>
\roman H^1(W,\Bobb R^{\times}) 
@>>>
\roman H^1(W,\Bobb Z /2) 
@>>>
0
\\
@.
@VVV
@VVV
@VV{\cong}V
@.
\\
0
@>>>
\Cal H^1(L,A)
@>>>
\roman{Pic}^{\roman{flat}}(L,A)
@>>>
\roman{Pic}(A)
@>>>
0 
\endCD
\tag6.8.1
$$
of abelian groups.
Thus
the group
$\roman{Pic}^{\roman{flat}}(L,A)$
may be constructed from
$\roman H^1(W,\Bobb R^{\times})$
and the group 
$\Cal H^1(L,A)$
of 
line bundles with a flat
$L$-connection
which are trivial as line bundles.

\smallskip
\noindent
{\smc Remark 6.9.}
Still under the circumstances of (6.7),
when $A$ is the algebra of smooth complex functions
on $W$ and
$L$  the
$(\Bobb C,A)$-Lie algebra 
arising from complexification of 
the Lie algebra
of smooth vector fields on $W$,
the exact sequence
(6.5.1)
boils down to the exact sequence
$$
0
@>>>
\roman H^1(W,\Bobb C^{\times})_0
@>>>
\roman H^1(W,\Bobb C^{\times}) 
@>>>
\roman H^2(W,\Bobb Z) 
@>>>
\roman H^2(W, \Bobb C); 
$$
here
$\roman H^1(W,\Bobb C^{\times})\cong
\roman{Hom}(\pi_1(W),\Bobb C^{\times})$,
the group of flat complex line bundles
and $
\roman H^1(W,\Bobb C^{\times})_0$
is the group
of flat complex line bundles
which are topologically trivial;
this group may also be described 
as the connected component
of the trivial homomorphism in
$\roman{Hom}(\pi_1(W),\Bobb C^{\times})$
or as the quotient
$\roman H^1(W,\Bobb C)\big/\roman H^1(W,\Bobb Z)$.
Notice that
the second cohomology group
$
\roman H^2(W,\Bobb Z) $
with integer coefficients
amounts to the Picard group of the ring $A$
of smooth complex functions on $W$.

\smallskip
\noindent
{\smc Example  6.10.}
Under the circumstances
of (6.7),
suppose  that $L$ arises from 
the tangent bundle
of a foliation $\Cal F$
on $W$.
Then the line bundle which corresponds to 
$Q_L$ amounts to the top exterior power of the conormal bundle
to $\Cal F$.
The left $(A,L)$-module structure on $Q_L$ is then the 
{\it Bott connection\/}
and,
provided the normal bundle
of $\Cal F$ is orientable, the elements of
$Q_L$
are transverse measures to $\Cal F$.
With the definition of modular class given 
in \cite\evluwein,
the modular class of $Q_L$ is there called the {\it modular class
of the foliation\/}.
See also \cite\weinsfte.
\smallskip
\noindent
{\smc Example  6.11.}
Let $L=\fra g$ be an ordinary  Lie algebra 
of finite dimension $n$ (say)
over a field
$k$,
considered as the Lie-Rinehart algebra $(k,\fra g)$,
with $A=k$.
Then the group $\Cal H^1(L,A)$
boils down to
$\roman H^1(\fra g,k) = \roman {Hom}(\fra g,k)$,
and the exact sequence (6.5.1)
yields an isomorphism from
$\roman H^1(\fra g,k)$ onto
$\roman{Pic}^{\roman{flat}}(L,A)$.
Moreover,
$C_{\fra g} = \Lambda_R^n \fra g^*$ and
$Q_{\fra g} = \roman{Hom}_k(C_{\fra g},k) =\Lambda_R^n \fra g$,
with the obvious respective right- and left $\fra g$-module structures,
cf. (1.5) and (4.5). Thus,
as an element of
$\roman H^1(\fra g,k)$,
the modular class  then comes down to
the adjoint character $\xi_0\colon \fra g \to k$
given by
$$
\xi_0(x) = \roman{tr}(\roman{ad}_x),\quad x \in \fra g.
$$
\smallskip
Given a general Lie-Rinehart algebra $(A,L)$
over a commutative ring $R$ and a free $A$-module
$M$ of rank one and being
endowed with a 
left $(A,L)$-module structure, 
with reference to a basis element $b$ of $M$,
for $\alpha \in L$,
we define the {\it divergence\/}
$\roman{div}_b\alpha$
of $\alpha$ by
$$
\alpha_M(b) =  (\roman{div}_b\alpha) b.
\tag6.12
$$

\smallskip
\noindent
{\smc Example  6.13.}
Let $L = A \otimes _R \fra g$ be an $(R,A)$-Lie algebra
of the kind considered in (2.12).
In view of (2.12),
the dualizing module $C_L$ of $L$ may be written 
$C_L= A\otimes_R\Lambda^n \fra g^*$ and,
in this description,
the right $(A,L)$-module structure on $C_L$ 
is given by (2.12.1). 
Hence, as an $A$-module,
$
Q_L = \roman{Hom}_A(C_L, \omega_A)
$
is plainly isomorphic to
$$
\roman{Hom}_R(C_{\fra g}, \omega_A)
\cong \Lambda^n\fra g \otimes_R \omega_A.
$$ 
Now the $A$-dual 
$\roman{Hom}_A(C_L,A)$
of the dualizing module
$C_L$
manifestly has the form
$A \otimes _R \Lambda^n \fra g$
and from this description
it is obvious
that
(2.12.4)
turns
this
$A$-dual 
$\roman{Hom}_A(C_L,A)$
into a left $(A,L)$-module;
likewise
$\omega_A$
inherits a left $(A,L)$-module 
structure and, with reference to these structures,
the canonical $A$-module isomorphism
$$
\roman{Hom}_A(C_L,\omega_A)
\cong 
(A \otimes _R \Lambda^n \fra g) \otimes _A \omega_A
\tag6.13.1
$$
is one of
left $(A,L)$-modules,
the left
$(A,L)$-module structure on the right-hand side being given by
(2.1).
This description of 
$Q_L =
\roman{Hom}_A(C_L,\omega_A)$
is somewhat simpler than
that given in (6.6) above for a {\it general\/} 
$(R,A)$-Lie algebra.
If, moreover, as an $A$-module,
$\omega_A$ is free, with basis $b$ (say),
the modular class of
$Q_L \cong \Lambda^n \fra g \otimes _R \omega_A$
lies in
the subgroup $\Cal H^1(L,A)$
of $\roman{Pic}^{\roman{flat}}(L,A)$
and is the class of 
$$
\xi_0 + \roman{div}_b
\in  \roman{Der}(\fra g,A) \cong \roman{Der}(L,A)
\tag6.13.2
$$
$\xi_0$ being the adjoint character given in (6.11).
This follows from the description
(6.13.1)
of the left $(A,L)$-module structure
on $Q_L$ for the special case under consideration.
See also 
what is said in \cite\evluwein.

\medskip\noindent{\bf 7. Poisson algebras}\smallskip\noindent
Let $A$ be a Poisson algebra, with Poisson structure
$\{\cdot,\cdot\}$,
and let $D_{\{\cdot,\cdot\}}$
be its $A$-module of formal differentials $D_A$
or a suitable
quotient $\overline D_A$ thereof
so that the canonical map from
$\roman {Hom}_A(\overline D_A, A)$ to
$\roman {Hom}_A(D_A, A) (\cong \roman{Der}(A))$
is an isomorphism.
Suppose
that $D_{\{\cdot,\cdot\}}$ is endowed with the
$(R,A)$-Lie algebra structure
induced by the Poisson structure;
see \cite\poiscoho\ (3.8) for details.
We further suppose that,
as an $A$-module, $D_A$ 
or the appropriate quotient $\overline D_A$ thereof
is finitely
generated and
projective of constant rank $n$;
then
the $A$-module
$\roman{Der}(A)$
is finitely
generated and projective
of constant rank $n$ as well, and $A$ is regular
in a sense explained in Section 3 above.
For example, when $A$ is the algebra of smooth functions
on a smooth real Poisson manifold $P$ of dimension $n$,
the appropriate object to be taken
is the $A$-module
of 1-forms $\overline D_A =\Omega^1(P)$, 
with its
$(\Bobb R,A)$-Lie algebra structure.
For a general Poisson algebra
$A$ of the kind under discussion,
given left- and
right $(A,D_{\{\cdot,\cdot\}})$-modules $M$ and $N$,
{\it Poisson cohomology\/}
$\roman H^*_{\roman{Poisson}}(A,M)$
and {\it Poisson homology\/} 
$\roman H_*^{\roman{Poisson}}(A,N)$
are then given by  (or, indeed, may then be defined by)
$$
\roman H^*_{\roman{Poisson}}(A,M)
=
\roman H^*(D_{\{\cdot,\cdot\}},M),
\quad
\roman H_*^{\roman{Poisson}}(A,N)
=
\roman H_*(D_{\{\cdot,\cdot\}},N).
$$
(When 
$A$ is not regular in our sense,
this may not be the appropriate definition;
see \cite\poiscoho\  for details.)
The dualizing module for
$D_{\{\cdot,\cdot\}}$,
which we denote by $C_{\{\cdot,\cdot\}}$,
is just
$$
C_{\{\cdot,\cdot\}}=\Lambda_A^{n} D_{\{\cdot,\cdot\}}^*
=\roman{Hom}_A(\Lambda_A^{n} D_{\{\cdot,\cdot\}},A)
\cong\Lambda_A^{n} \roman{Der}(A),
\tag7.1
$$
the last isomorphism being one of $A$-modules, so that
the right $(A,D_{\{\cdot,\cdot\}})$-module structure
on
$\Lambda_A^{n} \roman{Der}(A)$
is induced by that 
on $\roman{Hom}_A(\Lambda_A^{n} D_{\{\cdot,\cdot\}},A)$
via this isomorphism.
Given a right
$(A,D_{\{\cdot,\cdot\}})$-module $N$,
the inverse duality isomorphism
(2.11.2)
now manifestly has the form
$$
\roman H_k^{\roman{Poisson}}(A,N)
@>>>
\roman H^{n-k}_{\roman{Poisson}}(A,\roman{Hom}_A(C_{\{\cdot,\cdot\}},N)),
\tag7.2
$$
for every non-negative integer $k$.
By Theorem 3.7,
(7.2) may be obtained as the cap product with a suitable
fundamental class. Moreover,
the proof of Proposition 1.4
shows that an isomorphism of the
kind (7.2) may  even be taken to be induced by an isomorphism of complexes:
Let $K=K(A,D_{\{\cdot,\cdot\}})$, 
the Rinehart complex reproduced in (2.9.2),
which is in fact a projective resolution 
of $A$ in the category of left
$(A,D_{\{\cdot,\cdot\}})$-modules,
since $D_{\{\cdot,\cdot\}}$ is assumed to be projective as an $A$-module;
with this choice of $K$,
the left-hand side of the isomorphism (1.4.2)
is the ordinary complex calculating Poisson homology
with coefficients in $N$
while the right-hand side calculates the corresponding
Poisson cohomology with coefficients in 
$\roman{Hom}_A(C_{\{\cdot,\cdot\}},N)$,
cf. \cite\poiscoho.
In particular,
besides its obvious   left
$(A,D_{\{\cdot,\cdot\}})$-module structure,
the algebra $A$ inherits also  a {\it right\/}
$(A,D_{\{\cdot,\cdot\}})$-module structure;
it is given
the formula
$$
a(b(du)) = \{ab,u\},\quad a,b,u \in A,
\tag7.3
$$
where $du\in D_A$ refers to the formal differential of $u \in A$.
We denote the resulting
right
$(A,D_{\{\cdot,\cdot\}})$-module 
by $A_{\{\cdot,\cdot\}}$.
As a special case of (7.2), with $N=A_{\{\cdot,\cdot\}}$,
we now obtain
the  isomorphism
$$
\roman H_k^{\roman{Poisson}}(A,A_{\{\cdot,\cdot\}})
@>>>
\roman H^{n-k}_{\roman{Poisson}}(A,
\roman{Hom}_A(C_{\{\cdot,\cdot\}},A_{\{\cdot,\cdot\}})).
\tag7.4
$$
Here
$
\roman{Hom}_A(C_{\{\cdot,\cdot\}},A_{\{\cdot,\cdot\}}) 
$
is endowed with the  left
$(A,D_{\{\cdot,\cdot\}})$-module structure
given by
(2.3).
We note that
$
\roman{Hom}_A(C_{\{\cdot,\cdot\}},A_{\{\cdot,\cdot\}}) 
$
is plainly
isomorphic to
$\Lambda_A^{n} D_{\{\cdot,\cdot\}}$
as an $A$-module. Further,
in view of Proposition 4.6 (ii),
$\roman{Hom}_A(C_{\{\cdot,\cdot\}},A_{\{\cdot,\cdot\}})$
looks like a trace module,
with $V=A_{\{\cdot,\cdot\}}$.
We shall come back to this observation later.
When $A$ is the Poisson algebra of smooth functions
on a smooth real Poisson manifold $P$
and
$D_{\{\cdot,\cdot\}}$ the $A$-module of
1-forms $\Omega^1(P)$, with its  $(\Bobb R,A)$-Lie algebra
structure mentioned at the beginning of this section,
the isomorphism (7.4) is precisely the one
given in \cite\evluwein\ 
and written 
$$
\roman H^k(P,\Lambda^{\roman{top}}\roman T^*P)
\cong
\roman H_{n-k}(P)
$$
there.
A special case thereof may be found in Section 6 of \cite\chemlone.
\smallskip
Next we consider the left 
$(A,D_{\{\cdot,\cdot\}})$-module $Q_{D_{\{\cdot,\cdot\}}}$,
cf. (6.6) for its description.
Under the present circumstances, by definition, this module is of the form
$$
Q_{D_{\{\cdot,\cdot\}}} = \roman{Hom}_A(C_{\{\cdot,\cdot\}}, \omega_A)
\tag7.5.1
$$
where the
left $(A,D_{\{\cdot,\cdot\}})$-module structure on
$\roman{Hom}_A(C_{\{\cdot,\cdot\}}, \omega_A)$
is given by
(2.3), with reference to the right
$(A,D_{\{\cdot,\cdot\}})$-module structures on
$C_{\{\cdot,\cdot\}}$ and $\omega_A$.
On the other hand,  as $A$-modules,
$$
\roman{Hom}_A(C_{\{\cdot,\cdot\}},A) \otimes _A \omega_A
\cong 
\Lambda_A^{n} D_{\{\cdot,\cdot\}} \otimes _A \omega_A,
\tag7.5.2
$$
and 
the corresponding  left
$(A,D_{\{\cdot,\cdot\}})$-module structure on
${\Lambda_A^{n} D_{\{\cdot,\cdot\}} \otimes _A \omega_A}$
may be described
as being given by the Lie derivative on both
$\Lambda_AD_{\{\cdot,\cdot\}}$ and 
$$
\omega_A = 
\roman{Hom}_A(\Lambda_A^{n} \roman{Der}(A),A)
=\Lambda_A^{n} D_A.
\tag7.6
$$
We shall say a bit more about this description in the proof of
Theorem 7.9 below.
However, as $A$-modules,
$$
\Lambda_A^{n} D_{\{\cdot,\cdot\}}
\cong \roman{Hom}_A(C_{\{\cdot,\cdot\}},A),
\quad
\omega_A \cong \roman{Hom}_A(C_{\{\cdot,\cdot\}},A)
\tag7.7
$$
and, 
as we have  already observed,
with reference to the right
$(A,D_{\{\cdot,\cdot\}})$-module structure on
$C_{\{\cdot,\cdot\}}$ and 
that on $A$
given by (7.3) above
(indicated by the notation $A_{\{\cdot,\cdot\}}$),
(2.3) endows
$\roman{Hom}_A(C_{\{\cdot,\cdot\}},A_{\{\cdot,\cdot\}})$
with a left $(A,D_{\{\cdot,\cdot\}})$-module structure.

\proclaim{Lemma 7.8}
This 
left
$(A,D_{\{\cdot,\cdot\}})$-module structure on
$\roman{Hom}_A(C_{\{\cdot,\cdot\}},A_{\{\cdot,\cdot\}}) 
\cong \Lambda_A^{n} D_{\{\cdot,\cdot\}}$
is given by the formula
$$
(du)(\beta) = \lambda_{du}(\beta),\quad
u \in A,\ \beta \in \Lambda_A^{n} D_{\{\cdot,\cdot\}}.
$$
\endproclaim

Here
\lq\lq $\lambda$\rq\rq\ 
refers to the operation of Lie derivative (2.6.2);
notice that we do {\it not\/} assert that the
left
$(A,D_{\{\cdot,\cdot\}})$-module structure on
$\roman{Hom}_A(C_{\{\cdot,\cdot\}},A_{\{\cdot,\cdot\}})$ 
is given by the operation of Lie derivative
on $\Lambda_A^{n} D_{\{\cdot,\cdot\}}$;
only the result of the operation
with the $A$-module {\it generators\/}
$du$ ($u \in A$) of $D_{\{\cdot,\cdot\}}$ 
is given in this way.
We note that, when
$\Lambda_A^{n} D_{\{\cdot,\cdot\}}$
is actually free 
(necessarily of rank 1)
as an $A$-module, a choice
of basis element
$\beta$ of 
$\Lambda_A^{n} D_{\{\cdot,\cdot\}}$
determines a derivation
$\Phi_\beta$ of $A$
such that,
for $u \in A$, the result $\lambda_{du}(\beta)$
is given by the formula
$$
\lambda_{du}(\beta) = (\Phi_\beta u)\beta.
$$
When $A$ is the algebra of smooth functions
on a smooth (orientable) Poisson manifold $P$,
this derivation $\Phi_{\beta}$
is precisely the
\lq\lq modular vector field\rq\rq\ 
or \lq\lq curl\rq\rq\ 
(\lq\lq rotationnel\rq\rq\ in French), that is,
the infinitesimal
generator of the modular automorphism
group of $P$;
this automorphism group
was introduced in \cite\weinsfte,
and the modular vector field is written there in the form
$\Cal L_{H_f}\mu /\mu$.
The modular vector field
occurs already in
\cite\koszulon; see \cite\weinsfte\ 
for its significance and for additional references.

\demo{Proof}
Write $\Lambda =\Lambda_A^{n} D_{\{\cdot,\cdot\}}$.
Then $C_{\{\cdot,\cdot\}}\cong \roman{Hom}_A(\Lambda,A)$ as $A$-modules
and,
for $\beta\in \Lambda$ and $\phi\in \roman{Hom}_A(\Lambda,A)$,
the assignment
$\beta(\phi) = \phi(\beta)$
identifies
$\Lambda$ with
$$
\roman{Hom}_A(C_{\{\cdot,\cdot\}},A)
\cong \roman{Hom}_A(\roman{Hom}_A(\Lambda,A),A)
$$ 
as $A$-modules.
Let $\alpha\in D_{\{\cdot,\cdot\}}$. Then
$$
\align
(\alpha \beta) (\phi) 
&=
\beta(\phi \alpha) -(\beta\phi) \alpha 
\\
&=
(\phi \alpha)(\beta) -(\phi(\beta)) \alpha 
\\
&=
\phi (\lambda_\alpha(\beta))
-\alpha(\phi(\beta))
-(\phi(\beta)) \alpha, 
\endalign
$$
that is to say,
$$
\phi(\alpha (\beta)) = \phi(\lambda_\alpha (\beta)) 
-\alpha(\phi(\beta))
-(\phi(\beta)) \alpha.
$$
Let
$u \in A$ and $\alpha = du$; we then
have
$$
\alpha(\phi(\beta)) = \{u,\phi(\beta)\},
\quad
(\phi(\beta))\alpha = \{\phi(\beta),u\},
$$
whence
$
\phi(\alpha (\beta))=\phi(\lambda_\alpha (\beta))$.
Since $\phi$ is arbitrary, we conclude
that $\alpha (\beta)=\lambda_\alpha (\beta)$.
This proves the assertion. \qed
\enddemo

As a consequence we obtain the following.

\proclaim{Theorem 7.9}
With reference to the 
right 
$(A,D_{\{\cdot,\cdot\}})$-module
$A_{\{\cdot,\cdot\}}$
and the
left
$(A,D_{\{\cdot,\cdot\}})$-module structure on
$\roman{Hom}_A(C_{\{\cdot,\cdot\}},A_{\{\cdot,\cdot\}})$
(given by {\rm (2.3)}), we have
$$
Q_{D_{\{\cdot,\cdot\}}}
\cong \roman{Hom}_A(C_{\{\cdot,\cdot\}},A_{\{\cdot,\cdot\}}) \otimes_A 
\roman{Hom}_A(C_{\{\cdot,\cdot\}},A_{\{\cdot,\cdot\}})
\tag7.9.1
$$
as left 
$(A,D_{\{\cdot,\cdot\}})$-modules.
\endproclaim

\demo{Proof}
For a general duality $(R,A)$-Lie algebra $L$ of rank $n$,
it is readily seen that,
when 
the left $(A,L)$-module 
$Q_L = \roman{Hom}_A(C_L,\omega_A)$ (cf. (6.6))
is written in the form
$$
Q_L= \Lambda_A^n L \otimes_A \omega_A,
$$
the left $L$-module structure underlying
the left $(A,L)$-module structure
of $Q_L$
may be described 
as the 
ordinary tensor product
$L$-action 
arising
from the 
operations of Lie derivative
on the tensor factors
$\Lambda_A^n L$ and $\omega_A$;
neither of these Lie derivative actions defines a  
left $(A,L)$-module structure separately but their combination
indeed yields 
such a structure on the tensor product.
This remark applies in particular to
$L=D_{\{\cdot,\cdot\}}$
and
$\Lambda_A^n D_{\{\cdot,\cdot\}} \otimes _A \omega_A$.
By Lemma 7.8,
the result of the operation
on $\roman{Hom}_A(C_{\{\cdot,\cdot\}},A_{\{\cdot,\cdot\}})$
with the $A$-module generators
$du$ ($u \in A$) of $D_{\{\cdot,\cdot\}}$ 
is given by the operation of Lie derivative.
In view of the axioms for
a general $(R,A)$-Lie algebra $L$
and of those for a left $(A,L)$-module,
the left 
$(A,D_{\{\cdot,\cdot\}})$-module structure
on 
$\roman{Hom}_A(C_{\{\cdot,\cdot\}},A_{\{\cdot,\cdot\}})$
is then completely determined,
and the claim follows. \qed
\enddemo

\smallskip
\noindent
{\smc Remark 7.10.}
Theorem 7.9 generalizes
a corresponding result in \cite\evluwein,
and our argument 
involving the right
$(A,D_{\{\cdot,\cdot\}})$-module $A_{\{\cdot,\cdot\}}$
simplifies the proof thereof.
This result
has also been generalized in
\cite\xuone,
where other issues related to the present
paper are discussed, too.
See as well our follow up paper
\cite\bv.
\smallskip
We now return to our general Poisson algebra $A$
and write
$$
Q_{\{\cdot,\cdot\}} =
\roman{Hom}_A(C_{\{\cdot,\cdot\}},A_{\{\cdot,\cdot\}});
$$
the 
class
$[Q_{\{\cdot,\cdot\}}]$
in $\roman{Pic}^{\roman{flat}}(D_{\{\cdot,\cdot\}},A)$
of $Q_{\{\cdot,\cdot\}}$
is plainly characteristic for 
the Poisson algebra
$A$, and
we refer 
to it 
as the {\it modular class\/}
of $A$.
In view of Theorem 7.9,
$$
[Q_{\{\cdot,\cdot\}}]^2
=
[Q_{D_{\{\cdot,\cdot\}}}]
\in \roman{Pic}^{\roman{flat}}(D_{\{\cdot,\cdot\}},A),
$$
that is,
the modular class of the Lie-Rinehart algebra
$(A,D_{\{\cdot,\cdot\}})$
is twice the modular class 
(or its square when we think of the group structure
as being multiplicative)
of the Poisson algebra $A$.
We note that, when
$A$ is the Poisson algebra of smooth real functions on a
smooth compact Poisson manifold $P$,
the modular class is zero if and only if $P$ has a global volume form
which is annihilated by all Hamiltonian vector fields.

\smallskip
\noindent
{\smc Remark 7.11.}
When
$A$ is the Poisson algebra of smooth real functions on a
smooth Poisson manifold $P$,
our definition of the modular class of
$A$ is related to the definition of the modular
class of $P$
in \cite\weinsfte\ 
in the same way as
our definition of the modular class
$[Q_L]$ for a general Lie-Rinehart algebra
$(A,L)$ is related to the definition
of the modular class of a Lie algebroid
given in \cite\evluwein,
cf. Remark 6.7 above.
\smallskip\noindent
{\smc Illustration 7.12.}
Let $A$ be a Poisson algebra of functions on the dual of an $R$-Lie algebra
$\fra g$, the Poisson bracket being induced by the Lie bracket
on $\fra g$ in the customary way;
so $A$ could be the algebra of
smooth functions on the dual $\fra g^*$
of an ordinary  real (or complex) Lie algebra
$\fra g$,
or $A$ could be the algebra of polynomials
on an $R$-Lie algebra
$\fra g$ for a general ground ring $R$.
We suppose that, as an $R$-module,
$\fra g$ is finitely generated and projective of constant rank
$n$ (say); when $R$ is a field this is of course just the dimension
of $\fra g$. Now
$\fra g$
acts on
$A$ in an obvious fashion in such a way that
$D_{\{\cdot,\cdot\}}$
is isomorphic
to the corresponding $(R,A)$-Lie algebra
$A\otimes_R \fra g$
explained in (2.12) above.
See also \cite\poiscoho\ (3.18).
This implies that
$
C_{\{\cdot,\cdot\}}$
is isomorphic to
$A \otimes_R  \Lambda^n \fra g^*$
whence
$
\roman{Hom}_A(C_{\{\cdot,\cdot\}},A)$
is isomorphic to
$A \otimes_R  \Lambda^n \fra g$,
the $\fra g$-structures being the obvious ones.
Moreover,
$
Q_{\{\cdot,\cdot\}}$
is isomorphic to
$A \otimes_R  \Lambda^n \fra g \otimes_R  \Lambda^n \fra g $,
and the isomorphism (7.9.1) now comes down to the obvious isomorphism
between
$A \otimes_R  \Lambda^n \fra g \otimes_R  \Lambda^n \fra g $
and
$(A \otimes_R  \Lambda^n \fra g) \otimes_A (A \otimes_R  \Lambda^n \fra g)$.

\proclaim{Corollary 7.13}
The pairing {\rm (3.11)} induces a bilinear pairing
$$
\roman H_k^{\roman{Poisson}}(A,A_{\{\cdot,\cdot\}}) \otimes_{\Bobb R} 
\roman H_{n-k}^{\roman{Poisson}}(A,A_{\{\cdot,\cdot\}}) 
@>>>
\omega_A \otimes _UA.
\tag7.13.1
$$
\endproclaim

\demo{Proof}
The cohomology pairing (3.9) with reference
to the resulting pairing
$$
\roman{Hom}_A(C_{\{\cdot,\cdot\}},A_{\{\cdot,\cdot\}}) 
\otimes_A \roman{Hom}_A(C_{\{\cdot,\cdot\}},A_{\{\cdot,\cdot\}})
@>>>
Q_{D_{\{\cdot,\cdot\}}}
=\roman{Hom}_A(C_{\{\cdot,\cdot\}},\omega_A)
$$
of left $(A,D_{\{\cdot,\cdot\}})$-modules
takes the form
$$
\align
\roman H_{\roman{Poisson}}^k&
(A,\roman{Hom}_A(C_{\{\cdot,\cdot\}},A_{\{\cdot,\cdot\}})) 
\otimes_R 
\roman H^{n-k}_{\roman{Poisson}}
(A,\roman{Hom}_A(C_{\{\cdot,\cdot\}},A_{\{\cdot,\cdot\}}))
\\
&@>>>
C_{\{\cdot,\cdot\}} \otimes _U\roman{Hom}_A(C_{\{\cdot,\cdot\}},\omega_A)
\cong
\omega_A\otimes _U A
\endalign
$$
where $U= U(A,D_{\{\cdot,\cdot\}})$ 
denotes the corresponding universal algebra.
Naive duality (cf. (2.11) above) carries this pairing into (7.13.1). \qed
\enddemo

When
$A$ is the Poisson algebra of smooth real functions on an
orientable Poisson manifold $P$,
integration yields a map
from
$\omega_A \otimes _UA$
to the reals
which is non-zero when the manifold is compact or when we are working
with compactly supported functions and forms
and hence we then obtain a bilinear
real-valued pairing
$$
\roman H_k^{\roman{Poisson}}(P) \otimes_{\Bobb R} 
\roman H_{n-k}^{\roman{Poisson}}(P) 
@>>>
\Bobb R.
\tag7.14
$$
This is precisely the pairing given in 
of \cite\evluwein.
For reasons similar to those explained
in (3.10), this  pairing
will
not be nondegenerate
unless the Poisson structure is a symplectic one.
In the light of the notion of Poincar\'e duality
introduced in Section 4,
we now point a way towards
truely 
nondegenerate 
(co)homology pairings in Poisson (co)homology.

\smallskip \noindent 
{\smc Example 7.15.} 
Let $A$ be the algebra of smooth 
real functions on a smooth $n$-dimensional manifold $B$, endowed with the 
trivial Poisson structure.  
The corresponding $(\Bobb R,A)$-Lie algebra is just 
the space $D$ of sections of the cotangent bundle of $B$ with trivial Lie 
bracket, and the Poisson cohomology
$\roman H^*_{\roman{Poisson}}(A,A)$ is the 
graded $A$-algebra $\Lambda_A V$ 
of multi vector fields on $B$, where we have 
written $V = \roman{Vect}(B)$, the Lie bracket on $V$ being ignored.  In 
particular, letting $\Cal O = \Lambda^n_A D$, we have 
$$ \roman H^n_{\roman{Poisson}}(A,\Cal O) = \roman H^n(D,\Cal O) = 
\roman{Hom}_A(\Lambda^n_A D,\Lambda^n_A D) = A, 
$$ 
and this defines a trace $t \colon \roman H^n(D,\Cal O) \to A$.  
Poincar\'e duality now comes down to the 
fact that the canonical pairing 
$$ 
\Lambda^*_A D \otimes_A \Lambda^{n-*}_A D @>>> \Lambda^*_A D 
$$ 
is a perfect pairing of $A$-modules.  In general,
this notion of 
duality cannot even be phrased as one of real vector spaces!  Thus to 
understand Poincar\'e duality for Lie-Rinehart algebras one is forced to admit 
ground rings more general than the naive ones.  
\smallskip
We now switch to non-trivial Poisson structures.
Thus let $\{\cdot,\cdot\}$ be a non-trivial Poisson
structure on $A$.
In view of what has been said in Section 4, see in particular
(4.15),
the appropriate ground ring to be taken is the subalgebra
$R =\roman H^0_{\roman{Poisson}}(A,A) \subseteq A$
of Casimir elements in $A$.
Under suitable circumstances,
a trace can now be constructed by means of (4.6) (ii),
with $V= A_{\{\cdot,\cdot\}}$ (with the right
$(A,D_{\{\cdot,\cdot\}})$-module
structure
(7.3)):
The isomorphism $\iota$
(cf. (4.6))
then comes down to an isomorphism
$$
\iota
\colon
A_{\{\cdot,\cdot\}} \otimes _U A
\cong
\roman H_0^{\roman{Poisson}}(A,A)
@>{\iota}>>
\roman H^0_{\roman{Poisson}}(A,A) =R,
$$
that is, a trace boils essentially down
to an isomorphism from
$
A_{\{\cdot,\cdot\}} \otimes _U A
$
onto the algebra of Casimir
elements
such that,
for every
right duality
$(A,D_{\{\cdot,\cdot\}})$-module $N$,
the canonical map (4.6.3)
with
$L=D_{\{\cdot,\cdot\}}$
is an isomorphism. We now explain a situation where
indeed a trace is obtained in this way and where
Poincar\'e duality holds,
with the algebra of Casimir elements as ground ring.
\smallskip
\noindent
{\smc Example 7.16.}
Let 
$W$ be a smooth real $n$-dimensional Poisson manifold, 
$A$ its
Poisson algebra of smooth functions,
and suppose that
the symplectic foliation $\Cal F$ of $W$ constitutes a fiber bundle
with compact fiber $F$, of dimension $m$,
and base $B$.
For example,
$W$ could be the total space of a fiber bundle
whose fiber has a symplectic structure preserved by the
structure group.
An example
in which the cohomology class of the
symplectic structure \lq\lq varies\rq\rq\  from fiber to fiber
arises from
the regular part of the dual of a semisimple Lie algebra of compact type:
the fibers  are flag manifolds with a fixed complex
structure but with all \lq\lq different\rq\rq\  symplectic structures;
this follows from results in \cite\borhiron.
I am indebted to A. Weinstein for having pointed out these
examples to me.
The algebra
$\roman H^0_{\roman{Poisson}}(A,A)$ of Casimir elements
now amounts to the algebra $C^{\infty}(B)$
of smooth functions on $B$.
Further,  the image of
 the structure map from 
$D_{\{\cdot,\cdot\}}$ to
$\roman{Vect}(W)$
is the
$(\Bobb R,A)$-Lie algebra
$L_{\Cal F}$
which arises from the Lie algebroid
giving the infinitesimal structure of the foliation $\Cal F$,
and the resulting surjection from
$D_{\{\cdot,\cdot\}}$ to
$L_{\Cal F}$
fits into an extension
$$
0
@>>>
L'
@>>>
D_{\{\cdot,\cdot\}}
@>>>
L_{\Cal F}
@>>>
0
\tag7.16.1
$$
of
$(\Bobb R,A)$-Lie algebras
which, in turn, comes from an extension of Lie algebroids;
in particular,
besides
$D_{\{\cdot,\cdot\}}$
which has as underlying vector bundle
the cotangent bundle of $W$,
the $(\Bobb R,A)$-Lie algebras
$L'$ and $L_{\Cal F}$
have underlying vector bundles, too.
Let
$$
\tau\colon
\roman H^k(L_{\Cal F}, \Cal O_{\Cal F})
@>{e \cap \,\cdot\ }>>
C_{L_{\Cal F}} \otimes_U \Cal O_{\Cal F}
@>>> C^{\infty}(B)
$$
be the corresponding trace (4.15.1)
for $L_{\Cal F}$, and let
$(\Cal O,t)$
be the resulting trace for
$D_{\{\cdot,\cdot\}}$
given by (5.2).
We record the isomorphisms
$$
\Cal O = \roman{Hom}_A(C_{\{\cdot,\cdot\}}, C_{\Cal F}) 
\otimes _A \Cal O_{\Cal F}
\cong
\roman{Hom}_A(C_{\{\cdot,\cdot\}}, C_{\Cal F}\otimes _A \Cal O_{\Cal F}) 
\cong
\roman{Hom}_A(C_{\{\cdot,\cdot\}}, A_{\{\cdot,\cdot\}})
$$
of left
$(A,D_{\{\cdot,\cdot\}})$-modules
where $\roman{Hom}_A(C_{\{\cdot,\cdot\}}, A_{\{\cdot,\cdot\}})$
is endowed with the
left $(A,D_{\{\cdot,\cdot\}})$-module structure {\rm (2.3)};
now,
as an $A$-module,
$\roman{Hom}_A(C_{\{\cdot,\cdot\}}, A_{\{\cdot,\cdot\}})$
is
isomorphic to
$\Lambda^n_AD_{\{\cdot,\cdot\}}$,
that is,
to the space of sections of the highest non-zero
exterior power of the cotangent bundle of $W$.
Thus, as an $A$-module,
$\Cal O$ is just a copy of $A$ if and only if $W$ is orientable.

\proclaim{Theorem 7.16.2}
Under the present circumstances, the 
Poisson cohomology
cup pairing
$$
\roman H^k_{\roman{Poisson}}(A,A) \otimes_{C^{\infty}(B)}  
\roman H^{n-k}_{\roman{Poisson}}(A,\Cal O) 
@>>>
\roman H^n_{\roman{Poisson}}(A,\Cal O) 
@>t>>
C^{\infty}(B)
\tag7.16.3
$$
is a nondegenerate bilinear pairing
of finitely generated projective
$(C^{\infty}(B))$-modules.
\endproclaim

In view of the naturality
of the duality isomorphism (3.7.1),
naive duality (cf. (2.11) above) carries the pairing (7.16.3) into
a bilinear pairing in Poisson homology, and we obtain at once
the following.  

\proclaim{Corollary 7.16.4}
Under the circumstances of {\rm (7.16.2)}, the 
pairing {\rm (7.16.3)}
induces 
the nondegenerate bilinear pairing
$$
\roman H_{n-k}^{\roman{Poisson}}(A,C_{\{\cdot,\cdot\}}) 
\otimes_{C^{\infty}(B)}  
\roman H_k^{\roman{Poisson}}(A,A_{\{\cdot,\cdot\}}) 
@>>>
C^{\infty}(B)
\tag7.16.5
$$
of finitely generated projective
$(C^{\infty}(B))$-modules. \qed
\endproclaim

Thus, in Poisson (co)homology,
pairings of the kind (7.16.3) and (7.16.5) are more appropriate
objects of study
than those of the kind (7.14),
which are rarely nondegenerate, as already pointed out.
Furthermore,
the pairings
(7.16.3) and (7.16.5)
do not pose any problem
with the interpretation of the term
\lq\lq nondegenerate\rq\rq\ 
since the modules involved are finitely
generated over the corresponding ground ring,
the subalgebra of Casimir functions
(but in general not finitely generated over the reals).
\smallskip
We now prepare for the proof of Theorem 7.16.2.
To elucidate the structure of the extension (7.16.1),
consider the extension
$$
0
@>>>
\kappa_F
@>>>
\tau_W
@>>>
W\times_B \tau_B
@>>>
0
\tag7.16.6
$$
of vector bundles on $W$ arising from the projection
of the tangent bundle $\tau_W$ of $W$ to 
the tangent bundle $\tau_B$ of $B$.
The induced extension of the dual bundles may be written
$$
0
@>>>
W\times_B \tau^*_B
@>>>
\tau^*_W
@>>>
\kappa^*_F
@>>>
0.
\tag7.16.7
$$
The spaces of sections of this extension yield the extension
of $A$-modules which underlies (7.16.1);
in particular, we deduce that
$L'$ 
arises from the
projective $(C^{\infty}(B))$-module
$\Gamma (\tau^*_B)$
by extension of scalars, that is,
$L' \cong A \otimes _{C^{\infty}(B)}  \Gamma (\tau^*_B)$,
and this is an isomorphism of
$A$-Lie algebras, when
$\Gamma (\tau^*_B)$
and hence
$A \otimes _{C^{\infty}(B)}  \Gamma (\tau^*_B)$
are endowed with the trivial Lie bracket.
Now 
$\roman H^*(L',A)$ is just the exterior $A$-algebra
of $L'$ and, since
$L'$ is an induced module,
$$
\roman H^*(L',A)\cong 
A \otimes_{C^{\infty}(B)}  \left(\Lambda^*_{C^{\infty}(B)}\Gamma 
(\tau^*_B)\right).
\tag7.16.8
$$
\smallskip
Recall from (2.10) that the dualizing module $C_L$ 
for $L=D_{\{\cdot,\cdot\}}$ has the form
$C_L = \roman{Hom}_A(\Lambda^n_AL, A)$, endowed with the right
$(A,L)$-module structure (2.8.1), for $M=A$.
Inspection shows that, under the present circumstances,
the operation of Lie derivative of $L$ on $C_L$ factors through
the quotient $L_{\Cal F}$, and hence the
right $(A,L)$-module structure 
on $C_L$ 
factors through a right $(A,L_{\Cal F})$-module structure.
Consequently,
$$
\roman H_*(L',C_L)\cong 
 \left(\Lambda^*_{C^{\infty}(B)}\Gamma 
(\tau^*_B)\right)
\otimes_{C^{\infty}(B)}  C_L
\tag7.16.9
$$
as right 
$(A,L_{\Cal F})$-modules,
the right-hand side being endowed with the obvious
right 
$(A,L_{\Cal F})$-module structure
coming from that on $C_L$.

\demo{Proof of Theorem {\rm 7.16.2}}
As a 
$(C^{\infty}(B))$-module,
each module
$E^{p,q}_r(A)$
occurring in the cohomology spectral sequence
$(E^{p,q}_r(A),d_r)$
for the extension (7.16.1) with coefficients in $A$
is projective.
Hence 
each module
$E^{p,q}_{\infty}(A)$
is projective
as a 
$(C^{\infty}(B))$-module,
whence so 
are the cohomology modules
$\roman H^*_{\roman{Poisson}}(A,A)$
and $\roman H^*_{\roman{Poisson}}(A,\Cal O)$.
\smallskip
By Corollary 4.15.7,
the $(\Bobb R,A)$-Lie algebra $L_{\Cal F}$,
endowed with the trace
$(\Cal O_{\Cal F},\tau)$,
satisfies Poincar\'e duality for every finitely generated 
projective
induced left $(A,L)$-module of the kind
$M= A \otimes_{C^{\infty}(B)} M'$, 
$M'$ denoting an arbitrary
finitely generated projective
$(C^{\infty}(B))$-module.
Consequently $L_{\Cal F}$ satisfies
Poincar\'e duality for 
the right
$(A,L_{\Cal F})$-modules
$$
\roman H_0(L',C_L),
\roman H_1(L',C_L),
\dots,
\roman H_k(L',C_L).
$$
We now apply Theorem 5.3
and deduce that
the $(\Bobb R,A)$-Lie algebra $D_{\{\cdot,\cdot\}}$,
endowed with the trace
$(\Cal O,t)$,
satisfies Poincar\'e duality for 
the right $(A,L)$-module $C_L$.
Consequently
$D_{\{\cdot,\cdot\}}$
satisfies Poincar\'e duality for
the left $(A,L)$-module 
$\roman{Hom}_A(C_L,N)$
where $N$ is the right
$(A,L)$-module $C_L$, that is to say,
$D_{\{\cdot,\cdot\}}$
satisfies Poincar\'e duality for
just $A$, with its canonical 
left$(A,L)$-module structure,
as asserted.
This proves the claim. \qed
\enddemo

\smallskip
\noindent
{\smc Example 7.17.}
Let $\fra g$ be a real $n$-dimensional Lie algebra, 
consider its dual $\fra g^*$,
endowed with the ordinary Lie-Poisson structure,
and write $A$ for its Poisson algebra of smooth functions.
The subalgebra of Casimir functions
is the algebra $A^{\fra g}$ of invariants,
and we have the cup pairing
$$
\roman H^k_{\roman{Poisson}}(A,A) \otimes_{A^{\fra g}} 
\roman H^{n-k}_{\roman{Poisson}}(A,A) 
@>>>
\roman H^n_{\roman{Poisson}}(A,A).
\tag7.17.1
$$
Further, the Poisson
cohomology 
$\roman H^*_{\roman{Poisson}}(A,A)$
is well known to be isomorphic to
the Lie algebra cohomology
$\roman H^*(\fra g,A)$,
the algebra $A$ being endowed with the obvious $\fra g$-module structure,
cf. (3.18.4) in \cite\poiscoho,
and the pairing
(7.17.1)
boils down to the bilinear pairing
$$
\roman H^k(\fra g,A) \otimes_{A^{\fra g}} 
\roman H^{n-k}(\fra g,A) 
@>>>
\roman H^n(\fra g,A).
\tag7.17.2
$$
In view of what has been said in (4.13), when
$\fra g$ is of compact type,
the pairing
(7.17.2)
and hence
(7.17.1)
yields the
nondegenerate bilinear pairing
$$
\roman H^k_{\roman{Poisson}}(A,A) \otimes_{A^{\fra g}} 
\roman H^{n-k}_{\roman{Poisson}}(A,A) 
@>>>
A^{\fra g}
\tag7.17.3
$$
of finitely generated free
$A^{\fra g}$-modules.

\medskip 
\widestnumber\key{999} 
\centerline{References} \smallskip\noindent

\ref \no \almemoli
\by R. Almeida and P. Molino
\paper Suites d'Atiyah et feuilletages transversalement completes
\jour C. R. Acad. Sci. Paris I 
\vol 300
\yr 1985
\pages 13--15
\endref

\ref \no \bieriboo
\by R. Bieri
\book Homological dimension of discrete groups
\publ Queen Mary College, Department of Pure Mathematics
\publaddr London, U. K.
\yr 1976
\moreref 2nd edition: 1981
\finalinfo MR 57\#6224, 84h:20047
\endref

\ref \no \bierieck
\by R. Bieri and B. Eckmann
\paper Groups with homological duality generalizing
Poincar\'e duality
\jour Inventiones Math.
\vol 20
\yr 1973
\pages 103-124
\endref

\ref \no \borhiron
\by A. Borel and F. Hirzebruch
\paper Characteristic classes and homogeneous spaces. I
\jour Amer. J. Math.
\vol 80
\yr 1958
\pages 458--538
\moreref 
\paper II
\jour Amer. J. Math.
\vol 81
\yr 1959
\pages 315--382
\moreref
\paper III
\jour Amer. J. Math.
\vol 82
\yr 1960
\pages 491--504
\endref

\ref \no \boredmod
\by A. Borel  et al.
\book D-modules
\bookinfo Perspectives in Mathematics, vol. 2,
J. Coates and S. Helgason, eds.
\publ Academic Press
\publaddr Boston, Orlando, San Diego, New York, Austin,
London, Sydney, Tokyo, Toronto
\yr 1987
\endref

\ref \no \kbrownfo
\by K. S. Brown
\paper Homological criteria for finiteness
\jour Comment. Math. Helvetici
\vol 50
\yr 1975
\pages 129--135
\endref

\ref \no \cartanei
\by H. Cartan and S. Eilenberg
\book Homological Algebra
\publ Princeton University Press
\publaddr Princeton
\yr 1956
\endref

\ref \no \chemlone
\by S. Chemla
\paper Poincar\'e duality for $k$-$A$ Lie superalgebras
\jour Bull. Soc. math. France
\vol 122
\yr 1994
\pages 371--397
\endref
\ref \no \derhaboo
\by G. de Rham
\book Vari\'et\'es diff\'erentiables
\publ Hermann
\publaddr Paris
\yr 1955
\endref

\ref \no \evluwein
\by S. Evens, J.-H. Lu, and A. Weinstein
\paper Transverse measures, the modular class. and a cohomology pairing
for Lie algebroids
\paperinfo preprint,
{\tt dg-ga/9610008}
\endref

\ref \no \ginzwein
\by V. L. Ginzburg and A. Weinstein
\paper Lie-Poisson structure on some Poisson Lie groups
\jour J. Amer. Math. Soc.
\vol 5
\yr 1992 
\pages 445--453
\endref

\ref \no \godebook
\by R. Godement
\book Topologie alg\'ebrique et th\'eorie des faisceaux
\publ Hermann
\publaddr Paris
\yr 1964
\endref

\ref \no \gugenhtw
\by V.K.A.M. Gugenheim
\paper On the chain complex of a fibration
\jour Illinois J. of Mathematics
\vol 16
\yr 1972
\pages 398--414
\endref

\ref \no \hartshor
\by  R. Hartshorne
\book Algebraic Geometry
\bookinfo Graduate texts in Mathematics
 No. 52
\publ Springer
\publaddr Berlin-G\"ottingen-Heidelberg
\yr 1977
\endref

\ref \no \poiscoho
\by J. Huebschmann
\paper Poisson cohomology and quantization
\jour J. f\"ur die Reine und Angew. Math.
\vol 408
\yr 1990
\pages 57--113
\endref

\ref \no \bv
\by J. Huebschmann
\paper Lie-Rinehart algebras, Gerstenhaber algebras, and Batalin-
Vilkovisky algebras
\jour Annales de l'Institut Fourier
\vol 48
\yr 1998
\pages 425--440
\endref

\ref \no \extensta
\by J. Huebschmann
\paper 
Extensions of Lie-Rinehart algebras and the Chern-Weil construction
\paperinfo in: Festschrift in honor of J. Stasheff's 60th birthday
\jour Cont. Math. 
\vol 227
\yr 1999
\pages 145--176
\finalinfo{\tt dg-ga/9706002}
\endref

\ref \no \koszulon
\by J. L. Koszul
\paper Crochet de Schouten-Nijenhuis et cohomologie
\jour Ast\'erisque,
\vol hors-s\'erie,
\yr 1985
\pages 251--271
\paperinfo in: E. Cartan et les Math\'ematiciens d'aujourd'hui, 
Lyon, 25--29 Juin, 1984
\endref
\ref \no \mackbook
\by K. Mackenzie
\book Lie groupoids and Lie algebroids in differential geometry
\bookinfo London Math. Soc. Lecture Note Series, vol. 124
\publ Cambridge University Press
\publaddr Cambridge, England
\yr 1987
\endref

\ref \no \maclaboo
\by S. Mac Lane
\book Homology
\bookinfo Die Grundlehren der mathematischen Wissenschaften
 No. 114
\publ Springer
\publaddr Berlin $\cdot$ G\"ottingen $\cdot$ Heidelberg
\yr 1963
\endref

\ref \no \rinehart
\by G. Rinehart
\paper Differential forms for general commutative algebras
\jour  Trans. Amer. Math. Soc.
\vol 108
\yr 1963
\pages 195--222
\endref

\ref \no  \weinsfte
\by A. Weinstein
\paper The modular automorphism group of a Poisson manifold
\paperinfo in: Special volume in honor of A. Lichnerowicz
\jour J. of Geometry and Physics
\vol 23
\yr 1997
\pages 379--394
\endref

\ref \no \xuone
\by P. Xu
\paper 
Gerstenhaber algebras and BV-algebras
in Poisson geometry
{\tt dg-ga/9703001}
\endref
\enddocument